\documentstyle[12pt,amssymb,amscd]{amsart} 
\theoremstyle{definition}

\theoremstyle{remark}

\newcommand{\ov}{\overline}

\newcommand{\PP}{{\Bbb P}}
\newcommand{\CC}{{\Bbb C}}
\newcommand{\RR}{{\Bbb R}}
\newcommand{\TT}{{\Bbb T}}
\newcommand{\ZZ}{{\Bbb Z}}

\newcommand{\DD}{{\Bbb D}}
\newcommand{\iso}{\stackrel{\simeq}{\longrightarrow}}
\newcommand{\sm}{\smallsetminus}

\renewcommand{\marginpar}[1]{}
\catcode`\@=12

\def\Empty{}
\newcommand\oplabel[1]{
  \def\OpArg{#1} \ifx \OpArg\Empty {} \else
  	\label{#1}
  \fi}
		
\long\def\realfig#1#2#3#4{
\begin{figure}[htp]
\centerline{\psfig{figure=#2,width=#4}}
\caption[#1]{#3}
\oplabel{#1}
\end{figure}}

\newcommand{\comm}[1]{}
    
\input{psfig}
\begin{document}
\title{Dynamics of Singular Holomorphic Foliations on the Complex Projective Plane} 
\author{Saeed Zakeri}
\address{Department of Mathematics, University of Tehran, Iran, and Institute for Studies in Theoretical Physics and Mathematics (IPM), Tehran, Iran}
\curraddr{Department of Mathematics, SUNY at Stony Brook, NY 11794-3651, USA}
\email{zakeri@math.sunysb.edu}
\maketitle
\def\SBIMSMark#1#2#3{
 \font\SBF=cmss10 at 10 true pt
 \font\SBI=cmssi10 at 10 true pt
 \setbox0=\hbox{\SBF Stony Brook IMS Preprint \##1}
 \setbox2=\hbox to \wd0{\hfil \SBI #2}
 \setbox4=\hbox to \wd0{\hfil \SBI #3}
 \setbox6=\hbox to \wd0{\hss
             \vbox{\hsize=\wd0 \parskip=0pt \baselineskip=10 true pt
                   \copy0 \break%
                   \copy2 \break%
                   \copy4 \break}}
 \dimen0=\ht6   \advance\dimen0 by \vsize \advance\dimen0 by 8 true pt
                \advance\dimen0 by -\pagetotal
 \dimen2=\hsize \advance\dimen2 by .25 true in
%
%
  \openin2=publishd.tex
  \ifeof2\setbox0=\hbox to 0pt{}
  \else 
     \setbox0=\hbox to 3.1 true in{
                \vbox to \ht6{\hsize=3 true in \parskip=0pt  \noindent  
                {\SBI Published in modified form:}\hfil\break
                \input publishd.tex 
                \vfill}}
  \fi
  \closein2
  \ht0=0pt \dp0=0pt
 \ht6=0pt \dp6=0pt
 \setbox8=\vbox to \dimen0{\vfill \hbox to \dimen2{\copy0 \hss \copy6}}
 \ht8=0pt \dp8=0pt \wd8=0pt
 \copy8
 \message{*** Stony Brook IMS Preprint #1, #2 ***}
}

\SBIMSMark{1998/10}{September 1998}{}
\thispagestyle{empty}
\markboth{\sc Saeed Zakeri}
         {\sc Dynamics of Singular Holomorphic Foliations on $\CC\PP^2$}
\begin{center}
{\bf Table of Contents}
\begin{tabular}{ll}
{\small \bf Preface} & {\small \bf 2}\\
{\small \bf Introduction} & {\small \bf 3}\\
{\small \bf Chapter 1. Singular Holomorphic Foliations by Curves} & {\small \bf 8 }\\
{\footnotesize Holomorphic Foliations on Complex Manifolds} & {\footnotesize 8}\\
{\footnotesize Singular Holomorphic Foliations by Curves}\ \ \  & {\footnotesize 9}\\
{\footnotesize Polynomial SHFC's on $\CC \PP^2$} & {\footnotesize 11}\\
{\footnotesize Rigidity of SHFC's on $\CC \PP^2$} & {\footnotesize 16}\\
{\footnotesize Geometric Degree of a SHFC on $\CC \PP^2$} & {\footnotesize 22}\\
{\footnotesize Line at Infinity as a Leaf} & {\footnotesize 24}\\ 
{\small \bf Chapter 2. The Monodromy Group of a Leaf} & {\small \bf 30}\\
{\footnotesize Holonomy Mapping and the Monodromy Group} & {\footnotesize 30}\\
{\footnotesize Monodromy Pseudo-Group of a Leaf} & {\footnotesize 33}\\
{\footnotesize Multiplier of a Monodromy Mapping} & {\footnotesize 34}\\
{\footnotesize Monodromy Group of the Leaf at Infinity} & {\footnotesize 36}\\
{\footnotesize Equivalence of Foliations and Subgroups of Bih$_0(\CC)$} & {\footnotesize 40}\\
{\small \bf Chapter 3. Density and Ergodicity Theorems} & {\small \bf 44}\\
{\footnotesize Linearization of Elements in Bih$_0(\CC)$} & {\footnotesize 44}\\
{\footnotesize Approximation by Elements of a Pseudo-Group} & {\footnotesize 46}\\
{\footnotesize Ergodicity in Subgroups of Bih$_0(\CC)$} & {\footnotesize 47}\\
{\footnotesize Density of Leaves of SHFC's on $\CC \PP^2$} & {\footnotesize 51}\\
{\footnotesize Ergodicity of SHFC's on $\CC \PP^2$} & {\footnotesize 57}\\
{\small \bf Chapter 4. Non-Trivial Minimal Sets} & {\small \bf 59}\\
{\footnotesize Uniqueness of Minimal Sets} & {\footnotesize 60}\\
{\footnotesize Hyperbolicity in Minimal Sets} & {\footnotesize 62}\\
{\footnotesize Algebraic Leaves and Minimal Sets} & {\footnotesize 66}\\
{\small \bf References} & {\small \bf 69}\\
{\small \bf Symbol Index} & {\small \bf 71}
\end{tabular}
\end{center}
\newpage
\vspace*{4mm}
\begin{center}
{\bf Preface}
\end{center}
\vspace*{0.1in}
\thispagestyle{plain}
This manuscript is a revised version of my Master's thesis which was originally written in 1992 and was presented to the Mathematics Department of University of Tehran. My initial goal was to give, in a language accessible to non-experts, highlights of the 1978 influential paper of Il'yashenko on singular holomorphic foliations on $\CC \PP^2$ \cite{I4}, providing short, self-contained proofs. Parts of the exposition in chapters 1 and 3 were greatly influenced by the beautiful work of G\'omez-Mont and Ortiz-Bobadilla \cite{G-O} in Spanish, which contains more material, different from what we discuss here. It must be noted that much progress has been made in this area since 1992, especially in local theory (see for instance the collection \cite{I6} and the references cited there). However, Hilbert's $16$th Problem and the Minimal Set Problem are still unsolved. 

There is a well-known connection between holomorphic foliations in dimension $2$ and dynamics of iterations of holomorphic maps in dimension $1$, but many believe that this connection has not been fully exploited. It seems that some experts in each area keep an eye on progress in the other, but so far there have been rather few examples of a fruitful interaction. The conference on Laminations and Foliations held in May 1998 at Stony Brook was a successful attempt to bring both groups together. As a result, many people in dynamics expressed their interest in learning about holomorphic foliations. I hope the present manuscript will give them a flavor of the subject and will help initiate a stronger link between the young researchers in both areas.
     
My special thanks go to S. Shahshahani who encouraged me to study this subject years ago and supervised my Master's thesis back in Iran. I gratefully acknowledge the financial support of IPM (Tehran) during preparation of this manuscript. A series of lectures given by C. Camacho and Yu. Il'yashenko in Trieste, Italy, in 1991 and 1992 were a source of inspiration to me. I am indebted to X. G\'omez-Mont and J. Milnor who read parts of this revised version and made very useful comments. Finally, I would like to thank the Dynamical Systems group at Stony Brook, especially A. Epstein, M. Lyubich and J. Milnor. Without the interest they showed, I would have never found enough motivation to revise what I had written 6 years ago.\\ \\
\hspace*{10cm} Saeed Zakeri\\
\hspace*{9cm} Stony Brook, July 1998

\newpage
\vspace*{4mm}
\noindent
{\Large {\bf Introduction}}\\ \\ \\ \\ \\ \\ \\ \\ \\ \\ \\ \\
\thispagestyle{plain}
Consider the differential equation
$$(\ast) \left \{
\begin{array}{l}\vspace{3 mm}
\displaystyle{\frac{dx}{dt}=P(x,y)}\\
\displaystyle{\frac{dy}{dt}=Q(x,y)}
\end{array}
\right .$$
in the real plane $(x,y)\in {\RR}^2$, where $P$ and $Q$ are relatively prime
polynomials with $\max \{ {\rm deg}P, {\rm deg}Q \}=n$. What can be said, Hilbert
asked, about the number and the location of limit cycles of ($\ast$)? In particular,
is it true that there are only finitely many limit cycles for ($\ast$)? If so,
does there exist an upper bound $H(n)$, depending only on $n$,
for the number of limit cycles of an equation of the form ($\ast$)?
Surprisingly, the finiteness problem has been settled only in recent
years, and the existence and a value of $H(n)$ is still unknown, even when $n=2$!

In an attempt to answer the Hilbert's question, H. Dulac ``proved'' the finiteness
theorem in 1926 \cite{D}. Many years later, however, his ``proof'' turned out to be
wrong. In fact, in 1982 Russian mathematician Yu. Il'yashenko found a fundamental mistake in Dulac's
argument \cite{I3}, and gave a correct proof for the finiteness theorem later in 1987 \cite{I1}.

The second major attempt along this line was started in 1956 by a seminal paper
of two Russian mathematicians I. Petrovski\u\i\  and E. Landis \cite{P-L1}. They had
a completely different and perhaps more radical approach. They considered
($\ast$) as a differential equation in the complex plane $(x,y)\in \CC^2$,
with $t$ being a {\it complex time} parameter. The integral curves of the vector
field are now either singular points which correspond to the common zeros of
$P$ and $Q$, or complex curves tangent to the vector field which are holomorphically immersed in $\CC^2$. This gives rise to a holomorphic foliation by
complex curves with a finite number of singular points. One can easily see
that this foliation extends to the complex projective plane $\CC \PP^2$,
which is obtained by adding a line at infinity to the plane $\CC^2$. The
trajectories of ($\ast$) in the real plane are then the intersection of these
complex curves with the plane ${\rm Im}x={\rm Im}y=0$.

What makes this approach particularly useful is the possibility of applying
methods of several complex variables and algebraic geometry over an
algebraically closed field, not available in the real case. Intuitively,
the complexified equation provides ``enough space'' to go around and observe how the integral curves behave, whereas the real-plane dynamics of the trajectories is only the tip of a huge iceberg.

Viewing ($\ast$) as a complex differential equation, Petrovski\u\i\  and Landis
``proved'' that $H(2)$ exists and in fact one can take $H(2)= 3$ \cite{P-L1}. Later on, they ``proved'' the estimates
$$ H(n)\leq \left \{
\begin{array}{ll}\vspace{3 mm}
\displaystyle{\frac{6n^3-7n^2+n+4}{2}} & {\rm \ \ \ \ if\ } n {\rm \ is\ even},\\
\displaystyle{\frac{6n^3-7n^2-11n+16}{2}} & {\rm \ \ \ \ if\ } n {\rm \ is\ odd},
\end{array}
\right.$$
thus answering the Hilbert's question \cite{P-L2}. The result was regarded as a great achievement: Not only did they solve a difficult problem, but they introduced a truly novel method in the geometric theory of ordinary differential equations. However, it did not take long until a serious gap was discovered in their proof. Although the gap made their estimate
on $H(n)$ invalid, their powerful method provided a route to further studies in this direction.

In 1978, Il'yashenko made a fundamental contribution
to the problem. Following the general idea of Petrovski\u\i\  and Landis, he studied equations ($\ast$) with complex polynomials $P,Q$ from a topological standpoint without particular attention to the Hilbert's question. In his famous paper \cite{I4}, he
showed several peculiar properties of the integral curves of such equations.

From the point of view of foliation theory, it may seem that
foliations induced by equations like ($\ast$) form a rather small sub-class of all reasonable holomorphic foliations on $\CC \PP^2$. However, as long as we impose a mild condition on the set of singularities, it turns out that {\it every} singular holomorphic foliation on $\CC \PP^2$ is induced by a polynomial
differential equation of the form ($\ast$) in the affine chart $(x,y)\in \CC^2$. Fortunately, the mild condition is precisely what is needed in several complex variables: the singular set of the foliation must be an analytic subvariety
of codimension at least $2$, which is just a finite set in the case of the
projective plane.

Consider a closed orbit $\gamma$ of a smooth vector field in the real plane.
To describe the behavior of trajectories near $\gamma$, one has the simple and
useful concept of the {\it Poincar\'e first return map}: Choose a small 1-disk $\Sigma$
transversal to $\gamma$ at some point $p$, choose a point $q\in \Sigma$ near $p$,
and look at the first point of intersection with $\Sigma$ of the trajectory
passing through $q$. In this way, one obtains the germ of a smooth diffeomorphism
of $\Sigma$ fixing $p$. The iterative dynamics of this self-map of $\Sigma$
reflects the global behavior of the trajectories near $\gamma$.

For a singular holomorphic foliation on $\CC \PP^2$ induced by an equation
($\ast$), a similar notion, called the {\it monodromy mapping}, had already been used by Petrovski\u\i\  and Landis, and extensively utilized by Il'yashenko. A closed orbit should now be replaced by a non-trivial loop $\gamma$ on the leaf passing through a given point $p$, small transversal $\Sigma$ is a 2-disk, and the result of traveling over $\gamma$ on the leaf passing through a point on $\Sigma$ near $p$ gives the germ of a biholomorphism of $\Sigma$ fixing $p$, called
the {\it monodromy mapping associated with $\gamma$}. Note that all points in the
orbit of a given point on $\Sigma$ under this biholomorphism lie in the same leaf.
In this way, to each non-trivial loop in the fundamental group of the leaf we associate a self-map of $\Sigma$ reflecting the behavior of nearby leaves as
one goes around the loop. It is easily checked that composition of loops
corresponds to superposition of monodromy mappings, so the fundamental group of the leaf maps homomorphically onto a subgroup of the group of germs of biholomorphisms
of $\Sigma$ fixing $p$; the latter subgroup will be called the {\it monodromy
group} of the leaf. Thus a global problem can be reduced to a large extent
to the study of germs of biholomorphisms of $\CC$ fixing (say) the origin.
The group of all these germs is denoted by Bih$_0(\CC)$.

Now here is the crucial observation: For ``almost every'' equation of the form
($\ast$), the line at infinity of $\CC \PP^2$ with finitely many points deleted
is a leaf of the extended foliation. On the other hand, no leaf can be bounded
in $\CC^2$ (this is in fact more subtle than the Maximum Principle), so
every leaf has a point of accumulation on the line at infinity. Therefore, the monodromy group of the {\it leaf at infinity}, which is finitely-generated since the leaf is homeomorphic to a finitely-punctured Riemann sphere, gives us much
information about the global behavior of all leaves.

What Il'yashenko observed was the fact that almost all the dynamical properties of leaves have a discrete interpretation in terms of finitely-generated
subgroups of Bih$_0(\CC)$, the role of which is played by the monodromy
group of the leaf at infinity. Therefore, he proceeded to study these subgroups
and deduced theorems about the behavior of leaves of singular holomorphic foliations. The density and ergodicity theorems for these foliations are direct consequences
of the corresponding results for subgroups of Bih$_0(\CC)$. The density theorem
asserts that for ``almost every'' equation ($\ast$), every leaf other than the
leaf at infinity is dense in $\CC \PP^2$. The ergodicity theorem says that
``almost every'' foliation induced by an equation of the form ($\ast$) is {\it ergodic}, which
means that every measurable saturated subset of $\CC \PP^2$ has zero or full measure.

Yet another advantage of the reduction to the discrete case is the possibility of
studying consequences of equivalence between two such foliations. Two singular
holomorphic foliations on $\CC \PP^2$ are said to be {\it equivalent} if there
exists a homeomorphism of $\CC \PP^2$ which sends each leaf of the first foliation to a
leaf of the second one. Such an equivalence implies the equivalence between the
monodromy groups of the corresponding leaves at infinity. The latter equivalence
has the following meaning: $G,G'\subset {\rm Bih}_0(\CC)$ are
{\it equivalent} if there exists a homeomorphism $h$ of some neighborhood of
$0\in \CC$, with $h(0)=0$, such that $h\circ f\circ h^{-1}\in G'$ if and only if $f\in G$.
The correspondence $f\mapsto h\circ f\circ h^{-1}$ is clearly a group isomorphism.
Under typical conditions, one can show that the equivalence between two subgroups of Bih$_0(\CC)$ must in fact be holomorphic. This gives some invariants of the equivalence classes, and reveals a rigidity phenomenon. Applying these results to the monodromy groups
of the leaves at infinity, one can show the existence of {\it moduli of stability} and the phenomenon of {\it absolute rigidity} for these foliations.

All the above results are proved under typical conditions on the differential equations.
After all, the existence of the leaf at infinity (an algebraic leaf) plays a
substantial role in all these arguments. Naturally, one would like to know how to study
foliations which do not satisfy these conditions, in particular, those which do not admit any algebraic leaf. From this point of view, the dynamics of these foliations is far from being understood. The major achievements along this line have been made by the Brazilian and French schools. In their interesting paper \cite{C-L-S1}, C. Camacho, A. Lins Neto and P. Sad study non-trivial minimal sets of these foliations and show several properties of the leaves within a non-trivial minimal set. A {\it minimal set} of a foliation on $\CC \PP^2$ is a compact, saturated, non-empty subset of $\CC \PP^2$ which is minimal with respect to these three properties. A {\it non-trivial minimal set} is one which is not a singular point. It follows that the existence of a non-trivial minimal set is equivalent to the existence of a leaf which does not accumulate on any singular point. The fact that non-trivial minimal sets do not exist when
the foliation admits an algebraic leaf makes the problem much more challenging.
The major open problem in this context is the most primitive one: ``Does there exist a singular holomorphic foliation on $\CC \PP^2$ having a non-trivial minimal set?''

An intimately connected question concerns limit sets of the
leaves of these foliations. The classical theorem of Poincar\'e$-$Bendixson
classifies all possible limit sets for foliations on
the $2$-sphere: Given a smooth vector field on the $2$-sphere with a finite number of
singular points, the $\omega$- (or $\alpha$-) limit set of any point is either a singular point or a closed orbit or a chain of singular points and trajectories starting from one of these singular points and ending at another one. This enables us to understand the asymptotic behavior
of {\it all} trajectories. Naturally, one is interested in proving a complex version
of the Poincar\'e$-$Bendixson Theorem for foliations on $\CC \PP^2$. Now the concept of the limit set is defined as follows: Take a leaf ${\cal L}$ and let $K_1\subset
K_2\subset \cdots \subset K_n \subset \cdots $ be a sequence of compact subsets
of ${\cal L}$, with $\bigcup_{n\geq 1}K_n={\cal L}$. Then the {\it limit set} of ${\cal L}$
is by definition the intersection $\bigcap_{n\geq 1}\ \overline{{\cal L}\sm K_n}$. Incidentally, the limit set
of every non-singular leaf is non-vacuous, since it can be shown that no non-singular leaf
can be compact. Despite some results in this direction (see e.g. \cite{C-L-S2}), the problem of classifying possible limit sets is almost untouched.

Finally, it should be mentioned how the Hilbert's question for a real equation
($\ast$) can be translated into the complex language. It is easy to prove that
a limit cycle of a real equation ($\ast$) is a {\it non-trivial} loop on the
corresponding leaf of the complexified equation. So any piece of information about the
fundamental groups of leaves could be a step toward understanding the limit cycles. A non-trivial loop on a leaf of a complex equation ($\ast$) is said to
be a {\it (complex) limit cycle} if the germ of its associated monodromy is not
the identity map. Il'yashenko has shown that ``almost every'' equation ($\ast$) has a countably infinite number of homologically independent (complex) limit cycles; nevertheless this result does not have a direct bearing on the Hilbert's question. A complex version of the finiteness problem may be the following: Can the fundamental group of a (typical) leaf of such foliations be infinitely-generated? If not, does there exist an upper bound, depending only on the degree of $P$ and $Q$ in ($\ast$), for the number of generators of the fundamental groups? These questions, as far as I know, are not yet answered.
\newpage
\vspace*{4 mm}
\noindent
{\large {\bf Chapter 1}} \vspace{4mm}\\
{\Large{\bf Singular Holomorphic Foliations by Curves}} \\ \\ \\ \\ \\ \\ \\ \\ \\ \\ \\ \\
\thispagestyle{plain}
\noindent
This chapter introduces the concept of a singular holomorphic
foliation by (complex) curves on a complex manifold, which will be quite essential in subsequent chapters.
The condition imposed on the set of singularities is rather strong so that the local nature of the foliation can be
described by methods of several complex variables. In particular,
when the underlying manifold is $\CC \PP^2$ (2-dimensional complex projective space),
the very special geometry of the space allows one to do some algebraic geometry to show that all such foliations are induced by polynomial 1-forms $\omega $ on $\CC\PP^2$, the
leaves being the solutions of $\omega =0$. The main tools here are
extension theorems of several complex variables and the rigidity properties
of holomorphic line bundles on projective spaces.\\ \\
\vspace{4 mm}
{\large \bf Holomorphic Foliations on Complex Manifolds}\\
{\bf 1.1 Definition}\ \ Let $M$ be a complex manifold of dimension
$n$ and $0<m<n$. A {\it holomorphic (non-singular) foliation ${\cal F}$ of codimension}
$m$ on $M$ is an analytic atlas ${\cal A}= \{(U_i, \varphi_i)\}_{i\in I}$ for $M$
which is maximal with respect to the following properties:
\begin{enumerate}
\item[(i)]
For each $i\in I$, $\varphi _i$ is a biholomorphism $U_i\rightarrow A_i \times B_i$, where $A_i$ and $B_i$ are open polydisks in $\CC^{n-m}$ and $\CC^m$,
respectively.
\item[(ii)]
If $(U_i, \varphi _i)$ and $(U_j,\varphi _j)$ are in ${\cal A}$
with $U_i\cap U_j \neq \emptyset $, then $\varphi_{ij}:=\ \varphi_i \circ \varphi_j^{-1} \ :\ \varphi_j (U_i\cap U_j)\rightarrow \varphi_i (U_i \cap U_j)$
has the form $\varphi_{ij}(z,w)=(\psi_{ij}(z,w), \eta _{ij}(w))$, where
$(z,w)\in \CC^{n-m}\times \CC^m$, and $\psi_{ij}$ and $\eta _{ij}$
are holomorphic mappings into $\CC^{n-m}$ and $\CC^m$, respectively.
\end{enumerate}

Condition (ii) may also be expressed in the following way: Using coordinates
$(z_1, \ldots , z_n)$ for $\CC^n$, the mapping $(z_1,\ldots ,z_n) \stackrel{\varphi _{ij}}{\longmapsto }(\varphi_{ij}^1,\ldots , \varphi _{ij}^n)$
is required to satisfy $\partial \varphi _{ij}^k / \partial z_l =0$ for
$n-m+1\leq k\leq n$ and $1\leq l\leq n-m$.

Each $(U_i, \varphi_i)\in$ ${\cal A}$ is called a {\it foliation chart} (or a {\it flow box}). Given any foliation chart
$(U_i, \varphi_i)$, the sets $\varphi _i^{-1}(A_i \times \{w \}) , w \in B_i$, are called {\it plaques}
of ${\cal F}$ in $U_i$. Evidently these form a partition of $U_i$ into connected pieces of complex
submanifolds of dimension $ n-m $. Each $p\in M$ lies in at least one plaque. Two
points $p$ and $q$ are called {\it equivalent} if there exists a sequence
$P_1,\ldots , P_k$ of plaques such that $p\in P_1 , q\in P_k$, and $P_i \cap P_{i+1} \neq \emptyset , 1\leq i \leq k-1$. The {\it leaf} of ${\cal F}$ through {\it p},
denoted by ${\cal L}_p$, is the equivalence class of $p$ under this relation. Each leaf has a natural structure of a
connected $ (n-m)$-dimensional complex manifold which is injectively immersed in a holomorphic way in {\it M}. Two leaves are disjoint or else identical.\\ \\
{\bf 1.2  Remarks}

(a) A holomorphic foliation of codimension $(n-1)$ on an {\it n}-dimensional
complex manifold $M$ is also called a {\it holomorphic foliation by curves}.
Its leaves are immersed Riemann surfaces. As we will discuss later,
the field of complex lines tangent to the leaves determines a holomorphic line bundle on $M$.

(b) It can be easily checked that the above definition is equivalent to the
following, which is more natural from the geometric viewpoint: A holomorphic foliation $\cal F$ of codimension $m$ on an {\it n}-dimensional complex
manifold $M$ is a partition of $M$ into disjoint connected subsets $ \{ {\cal L}_\alpha \}$
(called the {\it leaves} of ${\cal F}$) such that for each $p \in M$ there exists a chart $(U,\varphi )$ around {\it p} with the property that $\varphi : U\rightarrow A \times B$
$(A \subset \CC^{n-m} \ {\rm and} \ B \subset \CC^m$ are open polydisks) maps
the connected components of ${\cal L}_\alpha \cap U$ to the {\it level sets}
$A \times \{w\}, w\in B.$\\ \\
\vspace {4 mm}
{\large \bf Singular Holomorphic Foliations by Curves}\\
Now we study those foliations which are allowed to have some ``mild'' singularities.
Recall that a subset $E$ of a complex manifold $M$ is called an {\it analytic subvariety} if
each $p\in M$ has a neighborhood $U$ on which there are holomorphic
functions $f_j: U\rightarrow \CC$, $1\leq j\leq k$, such that $E\cap U=\{x\in U: f_j(x)=0 , 1\leq j \leq k \}$.
Evidently every analytic subvariety of $M$ is closed, hence $M \sm E$ is itself a complex manifold of the same dimension as $M$.\\ \\
{\bf 1.3 Definition} \ \ Let $M$ be a complex manifold. A {\it singular holomorphic foliation by curves}
${\cal F}$ on $M$ is a holomorphic foliation by curves on $M \sm E$, where $E$ is an analytic
subvariety of $M$ of codimension $>1$. A point $p\in E$ is called a {\it removable
singularity} of ${\cal F}$ if there exists a chart $(U, \varphi)$ around $p$, compatible
with the atlas ${\cal A}$ of ${\cal F}$ restricted to $M \sm E$, in the
sense that $\varphi \circ \varphi _i^{-1}$ and $\varphi _i \circ \varphi ^{-1}$
have the form 1.1(ii) for all $(U_i,\varphi _i)\in {\cal A}$. The set of
all non-removable singularities of ${\cal F}$ in $E$ is called the {\it singular set} of
${\cal F}$, and is denoted by sing$({\cal F})$.\\

Naturally, as in the case of real 1-dimensional singular foliations on real surfaces, the most
important examples of singular holomorphic foliations are furnished by vector fields.\\ \\
{\bf 1.4 A Basic Example} \ \ Let $X=\sum _{i=1}^n f_i \ \partial /\partial z_i$
be a holomorphic vector field on a domain $U\subset \CC^n$. We further assume that the
Jacobian $(\partial f_i/\partial z_j)_{1\leq i, j\leq n}$ has rank$>$1 throughout the domain $U$.
Then $X$ vanishes on the analytic variety $\{z: f_1(z)=\cdots =f_n (z)=0 \}$
which has codimension $>1$ (possibly the empty set). By definition, a {\it solution}
of the differential equation $\stackrel{ \cdot }{z} =X(z)$ with {\it initial condition}
$p\in U$ is a holomorphic mapping $\eta : \DD(0,r)\rightarrow U$ such that
$\eta (0)=p$ and for every $T\in \DD(0,r),\ d\eta(T)/dT=X(\eta (T))$. The
image $\eta (\DD(0,r))$ is called a {\it local integral curve} passing through
$p$. By the Existence and Uniqueness Theorem of solutions of holomorphic differential equations \cite{I-Y},
each $p$ has such a local integral curve passing through it, and two local integral curves through $p$
coincide on some neighborhood of $p$. It follows that if $X(p)=0$, then its integral curve will be the point $p$
itself. If $X(p)\neq 0$, any local integral curve through $p$ is a disk holomorphically embedded in $U$.


\realfig{figthesis1}{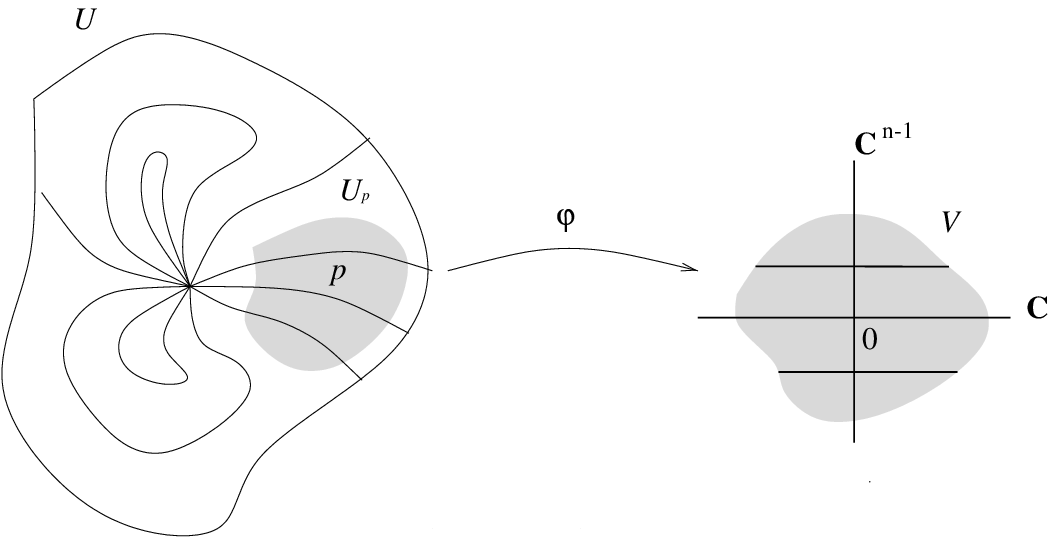}{\sl Straightening a holomorphic vector field near a non-singular point}{11cm}


Now suppose that $X(p)\neq 0$. By the Straightening Theorem for holomorphic vector fields \cite{I-Y},
there exist neighborhoods $U_p \subset U$ of $p$ and $V \subset \CC^n$ of 0 and a biholomorphism
$\varphi : U_p \rightarrow V$ such that $\varphi (p)=0$  and $\varphi _\ast (X|_{U_p})= \partial /\partial z_1$
(Fig. 1). Thus the connected components of the
intersection of integral curves and $U_p$ are mapped by $\varphi $ to ``horizontal'' lines
$\{z_2={\rm const}., \ldots,z_n={\rm const}. \}$. In other words, $X$ induces a singular holomorphic foliation by curves on $U$,
called ${\cal F}_X$, with sing$({\cal F}_X)=\{X=0 \}$, whose plaques are local integral curves of $X$.

Conversely, every singular holomorphic foliation by curves is locally induced by a holomorphic
vector field (cf. Proposition 1.14).\\

There is one property of these foliations which is quite elementary but will be frequently used in subsequent
arguments. Observe that there are much stronger results based on the theorem of analytic
dependence of solutions on initial conditions \cite{I-Y} as well as the concept of the holonomy mapping (see 2.1), but the following is sufficient for our purposes. \\ \\
{\bf 1.5 Proposition} \ \ {\it Let ${\cal F}$ be a singular holomorphic foliation by curves on M. Suppose that $p\in {\overline {\cal L}}_q$. Then
${\overline {\cal L}}_p \subset {\overline {\cal L}}_q.$}\\ \\
{\bf Proof}. There is nothing to prove if $p \in $ sing$({\cal F})$, so let {\it p} be non-singular.
Let $p'\in {\cal L}_p$ and join $p$ to $p'$ by a continuous path
$\gamma :[0,1]\rightarrow {\cal L}_p$ (thus avoiding sing$({\cal F})$) such that
$\gamma (0)=p$ and $\gamma (1)=p'$. Choose a partition $0=t_0<t_1<\cdots<t_n=1$
of [0,1] and foliation charts $(U_i,\varphi_i), \ 0\leq i \leq n-1$, such that
$\gamma [t_i,t_{i+1}]\subset U_i$. It follows from the local picture of plaques in
$\varphi _i(U_i)$ that $ \gamma (t_1)\in {\overline {\cal L}}_q$, so by repeating this argument,
$p'\in {\overline {\cal L}}_q$. Since $p'$ was arbitrary, we have ${\overline {\cal L}}_p\subset {\overline {\cal L}}_q$. $\hfill \Box$\\ \\
{\bf 1.6 Example} \ \ Example 1.4 has a counterpart in the context of differential forms defined on
domains in $\CC^n$. Our main example which will be seen to be quite general is the
following. Let $\omega =P(x,y) dy - Q(x,y)dx$ be a holomorphic 1-form on $\CC^2$, where $P$
and $Q$ are relatively prime polynomials. By definition, the singular foliation
induced by $\omega , \ {\cal F} _\omega :\{\omega =0\}$, is one which is induced on $\CC^2$
by the vector field $X(x,y)=P(x,y) \partial /\partial x +Q(x,y) \partial /\partial y $ as in Example 1.4.
In the language of 1-forms its leaves are obtained as follows: Take any
$p\in \CC^2$ at which $P$ and $Q$ are not simultaneously zero, and let $T \stackrel {\eta}{\mapsto} (x(T), y(T))$
be a holomorphic mapping on some disk $\DD(0,r)$ such that $\eta (0)=p$ and
\begin{equation}
P(x(T), y(T))y' (T)-Q(x(T),y(T))x' (T)=0
\end{equation}
for all $T\in \DD(0,r)$. The plaque through $p$ is the image under $\eta $
of some possibly smaller neighborhood of 0. Note that sing$({\cal F}_{\omega})$ is a finite set by the Bezout's Theorem. Observe that
\begin{equation}
{\cal F}_{\omega}={\cal F}_{f \omega}
\end{equation}
\noindent
for all holomorphic functions $f:\CC^2\rightarrow \CC$ which are
nowhere zero. The importance of this property lies in the fact that it
allows one to extend the foliation from $\CC^2$ to the complex projective plane
$\CC \PP^2$.\\

{\it Convention}. For simplicity, the term ``Singular Holomorphic Foliation by Curves'' will be abbreviated as SHFC from now on.\\ \\
\setcounter{equation}{0}
\vspace{4 mm}
{\large \bf Polynomial SHFC's on $\CC \PP^2$}\\
Now we make a digression to study polynomial SHFC's on $\CC \PP^2$ which are obtained by extending an ${\cal F}_{\omega}$
induced by a polynomial 1-form $\omega $ on $\CC^2$.\\ \\
{\bf 1.7 Geometry of $\CC \PP^2$} \ \ Consider $\CC^3 \sm \{(0,0,0)\}$ with the action of
$\CC^{\ast}$ defined by $\lambda . (x_0,x_1,x_2)$ $=(\lambda x_0, \lambda x_1 , \lambda x_2)$. The orbit of
$(x_0, x_1,x_2)$ is denoted by $[x_0,x_1,x_2]$. The quotient of $\CC^3 \sm \{(0,0,0)\}$
modulo this action (with the quotient topology) is called the {\it complex projective plane}
$\CC \PP^2$, and the natural projection is denoted by $\pi $. $\CC \PP^2$
can be made into a compact complex 2-manifold in the following way: Cover $\CC \PP^2$ by three open sets
\begin{equation}
U_i:=\{[x_0,x_1, x_2] : x_i\neq 0 \}, \ \ \ \ \ i=0,1,2
\end{equation}
and define homeomorphisms $\phi _i : \CC^2\rightarrow U_i$ by
\begin{equation}
\begin{array}{rl}
\phi_0(x,y) & =[1,x,y],\\
\phi_1(u,v) & =[u,1,v],\\
\phi_2(r,s) & =[r,s,1].
\end{array}
\end{equation}
The change of coordinates $\phi_{ij} = \phi _j^{-1}\circ \phi_i$ are given by
\begin{equation}
\begin{array}{l}
\phi_{01}(x,y)= \phi ^{-1}_1\circ \phi _0 (x,y)=\displaystyle{(\frac{1}{x} , \frac{y}{x})},\vspace{0.1in}\\
\phi _{12}(u,v)= \phi^{-1}_2\circ \phi _1 (u,v) =\displaystyle{ (\frac{u}{v} , \frac{1}{v})},\vspace{0.1in}\\
\phi _{20}(r,s)= \phi _0^{-1}\circ \phi _2 (r,s)=\displaystyle{ (\frac{s}{r} , \frac{1}{r})}.
\end{array}
\end{equation}
\noindent
These being holomorphic, the atlas $\{ (U_i , \phi _i^{-1} ), i=0,1,2 \}$
determines a unique complex structure on $\CC \PP^2$ for which the $\phi _i$
are biholomorphisms. Intuitively, $\CC \PP^2$ is a one-line compactification of $\CC^2$ for the
following reason: Each $(U_i, \phi^{-1}_i)$ is called an {\it affine chart}
of $\CC \PP^2$. Each $L_i:=\CC \PP^2 \sm U_i$ has a natural structure of the
Riemann sphere ${\overline \CC}$, since for example $L_0=\{[0,x,y]:(x,y)\in \CC^2 \}$
and it can be identified with $\{[x,y] \in \CC \PP^1:(x,y)\in \CC^2 \} \simeq {\overline \CC}$
under the restriction to $\CC^2$ of the action of $\CC^\ast $. Each
$L_i$ is called the {\it line at infinity} with respect to the affine chart $(U_i, \phi^{-1}_i)$. Every pair of the lines $L_i$ are intersecting in one point (Fig. 2).


\realfig{figthesis2}{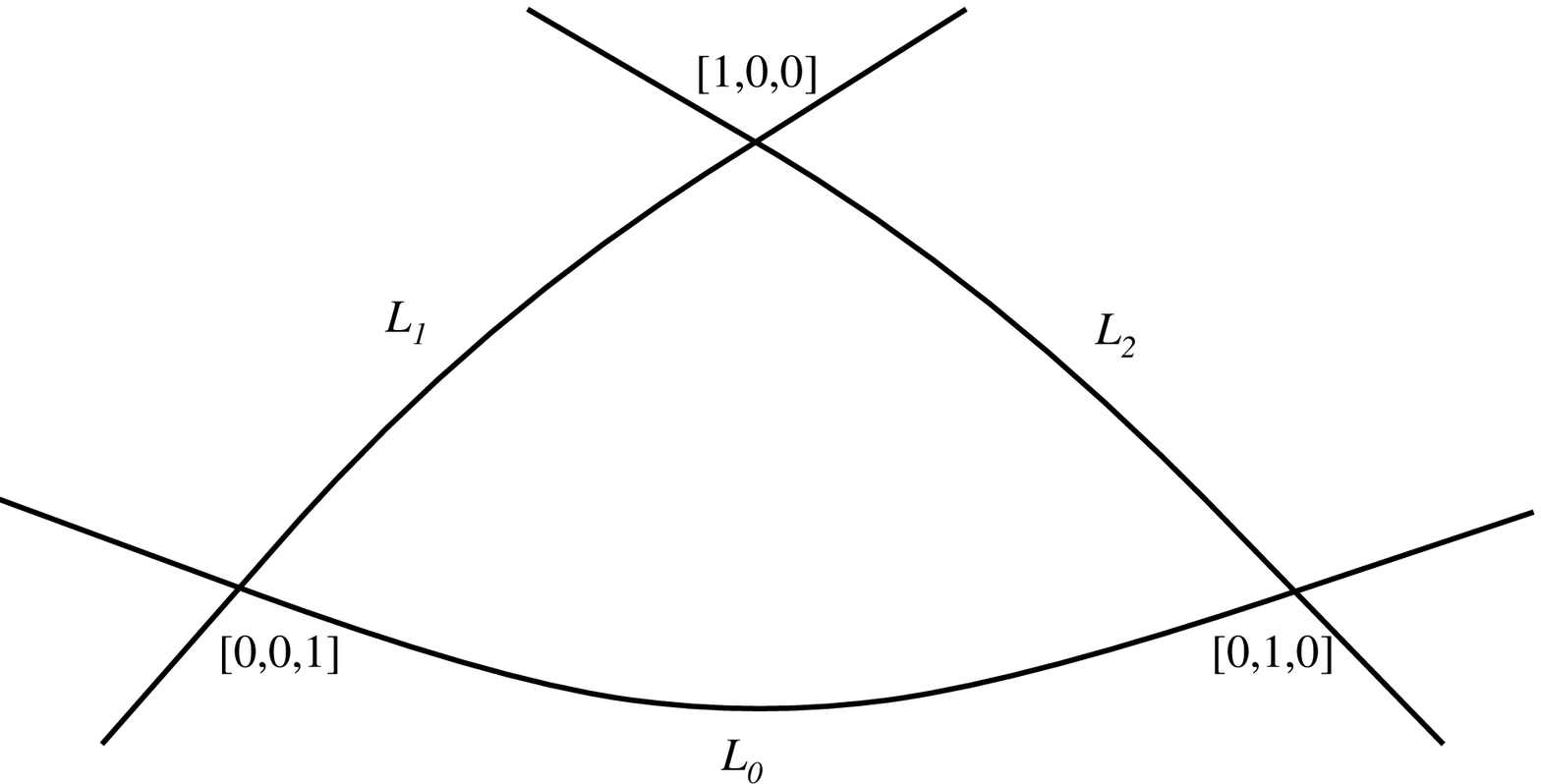}{\sl Geometry of $\CC \PP^2$}{10cm}


It is easy to check that given any {\it projective line} $L$ in $\CC \PP^2$, i.e., the projection under $\pi $ of any plane $ax_0+bx_1+cx_2=0$ in $\CC^3$, one
can choose a biholomorphism $\phi : \CC^2 \rightarrow \CC \PP^2 \sm L$. In this way, $L$ may be viewed as the line at infinity with respect to some affine chart.\\ \\
\setcounter{equation}{0}
{\bf 1.8 Algebraic Curves in $\CC \PP^2$} \ \ Suppose that $P=P(x,y)=\sum a_{ij}\ x^i y^j$
is a polynomial of degree $k$ on $(x,y)\in \CC^2$. Write $P$ in two other affine charts as (see 1.7 (3))
$$P \circ \phi_{10}(u,v)=P(\frac{1}{u} , \frac{v}{u})=u^{-k} \sum a_{ij} \ u^{k-(i+j)} v^j,$$
$$P \circ \phi _{20}(r,s)=P(\frac{s}{r} , \frac{1}{r})=r^{-k} \sum a_{ij} \ r^{k-(i+j)} s^i.$$
Set $P'(u,v)=\sum a_{ij} \ u^{k-(i+j)} v^j$ and $P'' (r,s)=\sum a_{ij} \ r^{k-(i+j)}s^i$. Then
the {\it algebraic curve} $S_P$ in $\CC \PP^2$ is defined as
$$\phi_0 \{(x,y):P(x,y)=0 \} \cup \phi_1 \{(u,v): P' (u,v)=0\} \cup \phi_2 \{(r,s): P'' (r,s)=0 \}.$$
Another way of viewing this curve is by introducing the homogeneous polynomial
$H_P$ of degree $k$ in $\CC^3$ as\\
$$H_P(x_0,x_1,x_2):= x_0^k P(\frac{x_1}{x_0} , \frac{x_2}{x_0})= \sum a_{ij}\ x_1^i x_2^j x_0^{k-(i+j)}.$$
It is then easily verified that $S_P=\pi \{(x_0,x_1,x_2): \ H_P (x_0,x_1,x_2)=0 \}.$ A
projective line is an algebraic curve $S_P$ for which deg $P=1$.\\ \\
{\bf 1.9 Extending Polynomial 1-froms on $\CC \PP^2$} \ \ Consider the polynomial
1-from $\omega =Pdy-Qdx$ on $\CC^2$ and its corresponding SHFC
${\cal F}_{\omega}$, as in 1.6. Using coordinate map $\phi_0$ (see 1.7(2))
one can transport ${\cal F}_\omega $ to $U_0$. To complete this picture to all of $\CC \PP^2$
we have to define the foliation on $L_0$. This can be done as follows. First
transport ${\cal F} _{\omega}$ to the affine chart $(u,v)$. To this end, write
$$\begin{array}{rl}
\tilde{\omega} (u,v)& :=(\phi_{10}^\ast \omega )(u,v)\vspace{0.1in}\\
 &=\displaystyle{P(\frac{1}{u},\frac{v}{u})d(\frac{v}{u})-Q(\frac{1}{u},\frac{v}{u})d(\frac{1}{u})}\vspace{0.1in}\\
 &=\displaystyle{u^{-1} P(\frac{1}{u},\frac{v}{u})dv-u^{-2} (vP(\frac{1}{u},\frac{v}{u})-Q(\frac{1}{u},\frac{v}{u}))du}.
\end{array}$$
Set $R(x,y):=yP(x,y)-xQ(x,y)$. Then
$$\tilde{\omega} (u,v)=u^{-1} P(\frac{1}{u},\frac{v}{u})dv-u^{-1} R(\frac{1}{u},\frac{v}{u})du.$$
Let $k$ be the least positive integer such that $\omega ':=u^{k+1}\tilde{\omega} $
is a polynomial 1-form on $(u,v)\in \CC^2$. Two foliations ${\cal F}_{\omega'}$ and ${\cal F}_{\tilde{\omega}} $
are then identical on $\{(u,v)\in \CC^2 :u\neq 0\}$ by $1.6 (2);$ however,
${\cal F}_{\omega'}$ which is defined on all of $(u,v)\in \CC^2$ is a well-defined extension of ${\cal F}_{\tilde{\omega}} $.
Now transport ${\cal F}_{\omega'}$ by $\phi_1$ to $U_1$. It is easily checked that the foliation induced by $({\cal F}_\omega ,\phi _0)$
coincides with that of $({\cal F}_{\omega'},\phi _1)$ on $U_0\cap U_1$.

In a similar way, ${\cal F}_\omega$ can be transported to the affine chart $(r,s)$
by $\phi_{20}$ to obtain a foliation ${\cal F}_{\omega''}$ induced by a polynomial 1-form $\omega ''$
on $(r,s)\in \CC^2$. Then ${\cal F}_{\omega ''}$ is transported to $U_2$ by $\phi _2$.

We still denote the extended foliation on $\CC \PP^2$ by ${\cal F}_\omega$ and frequently identify
${\cal F}_\omega , {\cal F}_{\omega'}, {\cal F}_{\omega''}$ with their transported companions on
$\CC \PP^2$. Thus, without saying explicitly, the affine charts $(x,y),(u,v),$ and $(r,s)$ are
considered as {\it subsets} of $\CC \PP^2$ itself by identifying them with $U_0,U_1$ and $U_2$, respectively.

It follows from the above construction that in each affine chart, ${\cal F}_\omega $ is
given by the integral curves of the following vector fields:\\
\begin{equation}
\begin{array}{ll}
{\rm In} \ (x,y) \in U_0, & X_0 = P(x,y) \displaystyle{\frac{\partial}{ \partial x} + Q(x,y) \frac{\partial}{ \partial y}},\vspace{0.1in} \\
{\rm In} \ (u,v) \in U_1, & X_1 = \displaystyle{u^{k} P(\frac{1}{u} , \frac{v}{u}) \frac{\partial}{\partial u}+ u^{k} R(\frac{1}{u} , \frac{v}{u}) \frac{\partial}{ \partial v}},\vspace{0.1in} \\
{\rm In} \ (r,s) \in U_2, & X_2 = \displaystyle{-r^{l} Q(\frac{s}{r},\frac{1}{r}) \frac{\partial}{\partial r} +r^{l} R(\frac{s}{r},\frac{1}{r}) \frac{\partial}{\partial s}},\\
\end{array}
\end{equation}
\setcounter{equation}{0}

\noindent
where $k$ and $l$ are least positive integers making the above vector fields into polynomial ones.
This shows that ${\cal F}_\omega $ is a SHFC on $\CC \PP^2$ with sing$({\cal F}_\omega )= S_P \cap S_Q \cap S_R$.\\ \\
{\bf 1.10 Holomorphic Vector Fields on $\CC \PP^2$} \ \ So far it should seem to be that
the SHFC's ${\cal F}_\omega $ constructed on $\CC \PP^2$ are realized as integral curves of holomorphic
{\it line fields} on $\CC \PP^2$ rather than {\it vector fields}. This distinction
should have been observed when we multiplied $\tilde{\omega}$ by a power of $u$ to cancel the pole at $u=0$
(cf. 1.9). Indeed, special geometry of
the projective plane makes the set of holomorphic vector fields on $\CC \PP^2$  very small, so that
in order to obtain sufficiently rich examples we have to allow SHFC's
which arise as integral curves of holomorphic line fields.

The following proposition describes all holomorphic vector fields on $\CC \PP^2$. The same characterization is true for $\CC \PP^n$, as can be shown by a coordinate-free argument \cite{C-K-P}, but here we present a very simple proof for $\CC \PP^2$.\\ \\
{\bf 1.11 Proposition} \ \ {\it Every holomorphic vector field on $\CC \PP^2$ is induced by a linear vector field on $\CC^3$.}\\ \\
{\bf Proof}.\ \ Let $X$ be a holomorphic vector field on $\CC \PP^2$, which has
the following expressions in three affine charts $U_0,U_1,$ and $U_2$:
\begin{center}
${\rm In} \ (x,y) \in U_0 \ : \ X_0 = \displaystyle{f_0\  \frac{\partial}{\partial x}+ g_0\  \frac{\partial}{\partial y}},$\vspace{0.1in}\\
${\rm In} \ (u,v) \in U_1 \ : \ X_1 = \displaystyle{f_1\  \frac{\partial}{\partial u}+ g_1\  \frac{\partial}{\partial v}},$\vspace{0.1in}\\
${\rm In} \ (r,s) \in U_2 \ : \ X_2 = \displaystyle{f_2\ \frac{\partial}{\partial r} + g_2\ \frac{\partial}{\partial s}}.$\\
\end{center}
Since $(\phi_{01})_\ast X_0=X_1$ and $(\phi _{02})_\ast X_0=X_2$, we obtain
\begin{equation}
f_1(u,v)=-u^2 f_0(\frac{1}{u} , \frac{v}{u})
\end{equation}
\begin{equation}
g_1(u,v)=-uv\ f_0(\frac{1}{u} , \frac{v}{u} )+u\ g_0 (\frac{1}{u} , \frac{v}{u}),
\end{equation}
and
\begin{equation}
f_2(r,s)=-r^2 g_0 (\frac{s}{r} , \frac{1}{r})
\end{equation}
\begin{equation}
g_2(r,s)=r f_0(\frac{s}{r} , \frac{1}{r}) -rs\ g_0 (\frac{s}{r} , \frac{1}{r}).
\end{equation}
\noindent
Consider the power series expansions $f_0(x,y)=\sum _{i,j\geq 0} a_{ij}\ x^i y^j$ and
$g_0(x,y)=\sum _{i,j \geq 0} b_{ij}\ x^i y^j$. Since $f_1,g_1,f_2$,$g_2$ are holomorphic
on $\CC^2$, we have the following:
\begin{center}
$a_{ij}=0 \ \ \ \ \ \ \ \ \ {\rm if} \ \ \ \ \ \ i+j \geq 3, \ \ \ \ \ \ \ \ \ \ \ \ \ \ \ \ \ \ ({\rm by} (1))$\vspace{0.1in}\\
$b_{ij}=0 \ \ \ \ \ \ \ \ \ {\rm if} \ \ \ \ \ \ i+j \geq 3, \ \ \ \ \ \ \ \ \ \ \ \ \ \ \ \ \ \ ({\rm by} (2))$\vspace{0.1in}\\
$a_{02}=b_{20}=0 \ , \ a_{20} = b_{11} \ , \ b_{02}= a_{11}, \ \ \ \ \ \ \ \ \ ({\rm by} (2))$\\
\end{center}
and (3) and (4) give no new relations. Now it can be easily checked that $X$ is induced by the linear vector field
$$\begin{array}{ll}
\tilde{X} & =(x_0-a_{20} x_1 - a_{11} x_2) \displaystyle{\frac{\partial}{\partial x_0}}\vspace{0.1in}\\
 & +(a_{00} x_0 +(a_{10}+1)x_1+a_{01}x_2) \displaystyle{\frac{\partial}{\partial x_1}}\vspace{0.1in}\\
 & +(b_{00} x_0 +b_{10}x_1+(b_{01}+1)x_2) \displaystyle{\frac{\partial}{\partial x_2}}\\
\end{array}$$
on $\CC^3$. Conversely, every linear vector field on $\CC^3$ gives a holomorphic vector field on $\CC \PP^2$, and we are done. $\hfill \Box$ \\

As a result, every holomorphic vector field on $\CC \PP^2$ is seen in the affine chart
$(x,y)$ as $X_0=f_0\ \partial /\partial x+g_0\ \partial /\partial y$, where
\begin{equation}
\begin{array}{l}
f_0(x,y)=a+bx+cy+x(dx+ey)\\
g_0(x,y)=a'+b'x+c'y+y(dx+ey),
\end{array}
\end{equation}
\setcounter{equation}{0}
\noindent
for some complex constants $a, a'$, etc.\\

\noindent
{\large \bf Rigidity of SHFC's on $\CC \PP^2$ }\vspace{4 mm} \\
Our aim here is to give in detail the proof of the remarkable fact that every SHFC on $\CC \PP^2$
is of the form ${\cal F}_\omega$ for some polynomial 1-form $\omega$ on $\CC^2$ (equivalently, a polynomial
vector field on $\CC^2$). The same proof works for every $\CC \PP^n$, $n\geq 2$, with only minor modifications. We refer the reader to \cite{G-O} for a more general set up.
Another proof for this fact can be achieved by methods of algebraic geometry (see \cite{I4}).

The proof goes along the following lines: First we associate to each SHFC on a complex manifold $M$
a (holomorphic) line bundle $B'\hookrightarrow TM$ over $M\sm E$ ($E$ is an analytic subvariety of $M$, as in Definition 1.3). Then we show that $B'$ may be
extended to a tangent line bundle $B$ over $M$ (Theorem 1.15). This establishes a natural relationship between
SHFC's on $M$ and holomorphic ``bundle maps'' $\beta :B\rightarrow TM$. Rigidity of line bundles in the case $M=\CC \PP^2$ will then be applied to show that each bundle map $\beta :B\rightarrow T\CC \PP^2$ is induced
by a polynomial 1-form on $\CC^2$ (Corollary 1.21). The foundational material on holomorphic line bundles
on complex manifolds used here can be found in \cite{G-H} or \cite{Kod}.\\ \\
{\bf 1.12 SHFC's and Line Bundles} \ \ Let ${\cal F}$ be a SHFC on a complex manifold $M$. Let $\{(U_i,\varphi_i)\}_{i\in I}$ be the
collection of foliation charts on $M':=M\sm E$. By 1.1(ii) the transition functions $\varphi_{ij}=\varphi_i\circ\varphi_j^{-1}=(\varphi_{ij}^1,\ldots,\varphi_{ij}^n)$ satisfy $\partial\varphi_{ij}^k/\partial z_1=0$ for $2\leq k\leq n$. Applying the chain rule to $\varphi^1_{ij}=\varphi^1_{ik}\circ \varphi^1_{kj}$, we obtain
\begin{equation}
\frac {\partial \varphi_{ij}^1} {\partial z_1}(p)=\frac {\partial \varphi_{ik}^1} {\partial z_1}(\varphi_{kj}(p)) \frac {\partial \varphi_{kj}^1} {\partial z_1}(p)
\end{equation}
for every $p\in \varphi_i(U_i\cap U_j\cap U_k)$. Define $\xi_{ij}: U_i\cap U_j\rightarrow \CC^{\ast}$ by
\begin{equation}
\xi_{ij}(p):=\frac {\partial\varphi_{ij}^1} {\partial z_1}(\varphi_j(p)).
\end{equation}
Evidently one has the cocycle relation $\xi_{ij}=\xi_{ik} \cdot \xi_{kj}$ on $U_i\cap U_j\cap U_k$, thus obtaining a holomorphic line bundle $B'$ on $M'$.
This line bundle has a natural holomorphic injection into $TM$ in such a way that the image of the fiber over $p$ under this injection
coincides with the tangent line to ${\cal L}_p$ at $p$. To see this, observe that
$$B'= \bigcup_{i\in I} (U_i\times\CC)/\sim,$$
where $(p,t)\in U_i\times\CC$ is identified with $(q,t')\in U_j\times\CC$ if and only if $p=q$ and $t=\xi_{ij}(p)t'$. Let $B_p' \ , p\in M',$ be the fiber of
$B'$ over $p$ and define $\beta_i':U_i\times\CC\rightarrow TM|_{U_i}$ by
\begin{equation}
\beta_i'(p,t):=t((\varphi_i^{-1})_\ast\frac{\partial}{\partial z_1})(p).
\end{equation}
\noindent
It follows then from (2) that if $p\in U_i\cap U_j$ and $(p,t)\sim (p,t')$, then $\beta_i'(p,t)=\beta_j'(p,t').$ Therefore, the $\beta_i$ give rise to a holomorphic
bundle map $\beta':B'\rightarrow TM$ which is injective and $\beta'(B_p')$ is the tangent line to ${\cal L}_p$ at $p$, a subspace of $T_pM$.

Now there is a natural question: Can $B'$ be extended to a line bundle over all of $M$? The crucial point for the answer, which is affirmative, is the condition
that the codimension of $E$ is $>1.$ Recall the following classical theorem of F. Hartogs (see \cite{W} for a proof).\\
\setcounter{equation}{0}

\noindent
{\bf 1.13 Theorem} \ \ {\it Let $U\subset \CC^n$ be a domain and $E\subset U$ be an analytic subvariety of $U$ of codimension $>1.$ Then every holomorphic (resp. meromorphic) function
on $U\sm E$ can be extended to a holomorphic (resp. meromorphic) function on U}. $\hfill \Box$ \\

Using this theorem one can extend line bundles induced by SHFC's. To this end, we have to prove the following simple but remarkable proposition (cf. \cite{G-O}). The letter $E$ will always denote an analytic subvariety of the ambient
space which satisfies Definition 1.3. \\ \\
{\bf 1.14 Proposition} \ \ {\it Let $U\subset \CC^n$ be a domain and ${\cal F}$ be a SHFC on U. Then for each $p\in U$ there exists a holomorphic vector field X on some neighborhood
$U_p$ of p such that X is non-vanishing on $U_p\sm E$ and is tangent to the leaves of ${\cal F}$. This X is unique up to multiplication by a holomorphic function
which is non-zero in a neighborhood of p. Moreover, $q\in U_p\cap E$ is a removable singularity of ${\cal F}$ if and only if $X(q)\neq 0$.}\\ \\
{\bf Proof.} There is nothing to prove if $p\not \in E$, so let $p\in E$ and let $U_p$ be a connected neighborhood of $p$ in $U$. Then $U_p':=U_p\sm E$
is open and connected. Each $q\in U_p'$ has a small connected neighborhood $U_q\subset U_p'$ on which there is a holomorphic vector field representing ${\cal F}$ on $U_q$. Without loss of generality
we may assume that the first component of these vector fields is not identically zero over $U_p'$. If $(Y_1,\ldots,Y_n)$ is the vector field representing ${\cal F}$ on $U_q$, then $(1,Y_2/Y_1,\ldots,Y_n/Y_1)$ is a
meromorphic vector field on $U_q$ representing ${\cal F}$ on $U_q\sm \{ {\rm zeros\ of\ } Y_1\}$. Repeating this argument for each $q\in U_p'$, one concludes that there are meromorphic functions $F_2,\ldots,F_n$ defined on
$U_p'$ such that $(1,F_2,\ldots,F_n)$ represents ${\cal F}$ on $U_p'\sm \{ {\rm poles\  of\ the\ } F_i$ in $U_p'\}$. By Theorem 1.13 each function $F_i$ may be extended to a meromorphic function on $U_p$ (still denoted by $F_i$)
since the codimension of $U_p\cap E$ is $>1.$ The germ of $F_i$ can be uniquely written as $F_i=f_i/g_i$ by choosing $U_p$ small enough, where $f_i$ and $g_i$ are relatively prime holomorphic functions
on $U_p$. Then $X:=(g,gF_2,\ldots,gF_n)$ is a holomorphic vector field on $U_p$ representing ${\cal F}$ on $U_p'\sm \{ {\rm zeros\ of\ } g\}$, where $g$ is the least common multiple of the $g_i$. Note that  codim $\{ z: X(z)=0\} >1.$

Now let $q\in U_p'$, and $(U,\varphi )$ be a foliation chart around $q$. The vector field $\varphi_\ast X$ has the form $h\ \partial /\partial z_1$ since it is tangent to the horizontal lines $\{ z_2={\rm const.},\ldots,z_n={\rm const.}\} $
away from zeros of $g\circ \varphi^{-1}$. Since the zero set of $h$ either is empty or has codimension 1, while the zero set of $\varphi _\ast X$ has codimension $>$1, it follows that $h$ is nowhere zero and $X(q)\neq 0$. Thus
$X(q)$ is tangent to ${\cal L}_q$ at $q$.

For the uniqueness part, let $\tilde{X}$ be another such vector field. Then $\tilde{X}=\xi X$ on $U_p'$, where $\xi :U_p'\rightarrow\CC^\ast $ is holomorphic. Extend $\xi $ over $U_p$ by Theorem 1.13. This new $\xi $ is nowhere vanishing, since
its zero set, if non-empty, would have codimension 1.

Finally, let $q\in U_p\cap E$ be a removable singularity of $\cal F$. Choose a compatible chart $(U,\varphi)$ around $q$ (cf. 1.3) and let $\tilde{X}= \varphi^{-1}_\ast (\partial/\partial z_1)$. Then $\tilde{X}$ describes $\cal F$ on $U$, so by the
above uniqueness we have $\tilde{X}=\xi X$, with $\xi$ being a holomorphic non-vanishing function on some neighborhood of $q$. Thus $X(q)\neq 0$. Conversely, suppose that $X(q)\neq 0$, and let $(U,\varphi)$ be a local chart around $q$ straightening
the integral curves of $X$, i.e., $\varphi_\ast (X|_U)=\partial/\partial z_1$. Then it is evident that $(U,\varphi)$ is compatible with every foliation chart of $\cal F$ in $U_p'$. $\hfill \Box$\\

Now let ${\cal F}$ be a SHFC on $M$ and $\beta':B'\rightarrow TM$ be the bundle map constructed in 1.12. According to Proposition 1.14, each $p\in M$ has a neighborhood $U_i$ and a holomorphic vector field $X_i$ defined on $U_i$ representing
${\cal F}$ on $U_i\sm E$. By the uniqueness part of 1.14, whenever $U_i\cap U_j \neq \emptyset $, we have $X_j=\xi _{ij} X_i$ on $U_i\cap U_j$, where $\xi _{ij}:U_i\cap U_j\rightarrow \CC^\ast$ is holomorphic.
Let $B$ be the line bundle over $M$ defined by the cocycle $\{ \xi _{ij}\}$. As in the construction of $\beta'$ in $1.12(3)$, define $\beta_i:U_i\times \CC\rightarrow TM|_{U_i}$ by
\begin{equation}
\beta_i(p,t):=tX_i(p).
\end{equation}
\noindent
Definition of $B$ shows that $\beta_i$'s patch together to yield a well-defined
bundle map $\beta :B\rightarrow TM$ for which $\beta (B_p)$ is the tangent line
to ${\cal L}_p$ at $p$ if $p\not \in E$. Note that $B$ is an extension of $B'$ since
$\beta^{-1}\circ \beta':B'\rightarrow B|_{M'}$ is an isomorphism of line bundles.
If $\eta :B\rightarrow TM$ is another bundle map which represents ${\cal F}$
away from $E$, then $\eta :U_i\times \CC\rightarrow TM|_{U_i}$ satisfies $\eta
(p,t)=t\lambda_i(p)X_i(p)$ for $p\in U_i\sm E$, where
$\lambda_i:U_i\sm E\rightarrow
\CC$ is non-vanishing. Extend $\lambda_i$ to $U_i$ by
Theorem 1.13. Note that  the action of $\eta$ on $B_p$ is
well-defined, so that $\lambda_i(p)=\lambda_j(p)$ if
$p\in U_i\cap U_j$. Defining $\lambda:M\rightarrow
\CC$ by $\lambda|_{U_i}:=\lambda_i$, we see that
$\eta =\lambda \cdot \beta $. Furthermore, $\lambda$ can
only vanish on $E$, but its zero set,
if non-empty, would have to have codimension 1. Therefore, $\lambda$ does not
vanish at all. Finally, the last part of Proposition 1.14 and relation (1) show
that $p\in E$ is a removable singularity of ${\cal F}$ if and only if $\beta (B_p)\neq 0$.

Summarizing the above argument, we have proved the following\\
\setcounter{equation}{0}

\noindent
{\bf 1.15 Theorem} \ \ {\it The line bundle $B'$ over $M\sm E$ associated
to a SHFC \ ${\cal F}$ on M can be extended to a line bundle B over M. There exists
a holomorphic bundle map $\beta :B\rightarrow TM$ for which $\beta (B_p)$ is the tangent line to
${\cal L}_p$ at p if $p\not \in E$. $\beta $ is unique up to multiplication by
a nowhere vanishing holomorphic function on M. A point $p\in E$ is a
removable singularity of ${\cal F}$ if and only if $\beta (B_p)\neq 0$}. $\hfill \Box $\\ \\
{\bf 1.16 Remarks}

(a) It follows from the last part of Theorem 1.15 that sing(${\cal F}$) is precisely
$\{ p\in M: \beta (B_p)=0\}$; in particular, after removing all the removable singularities
of ${\cal F}$ in $E$, sing(${\cal F}$) turns out to be an analytic subvariety of $M$ of codimension $>1$.

(b) If $M$ is a {\it compact} complex manifold (in particular, if $M=\CC \PP^2$),
then every two bundle maps $\beta ,\beta':B\rightarrow TM$ representing ${\cal F}$
differ by a non-zero multiplicative constant.

(c) Suppose that ${\cal F}$ is a SHFC on $M$ and $\beta :B\rightarrow TM$ and
$\tilde{\beta} :\tilde{B}\rightarrow TM$ are two bundle maps, where $B$ and $\tilde{B}$
are two line bundles over $M$ representing ${\cal F}$. By the construction of
our bundle maps, we may suppose that $B$ (resp. $\tilde{B}$) is formed by $(\{U_i\},\{ \xi_{ij}\} )$
(resp. $(\{V_k\},\{ \eta_{kl}\} )$) and there are holomorphic vector fields $X_i$
on $U_i$ (resp. $Y_k$ on $V_k$) satisfying $X_j=\xi_{ij}X_i$ on $U_i\cap U_j$
(resp. $Y_l=\eta_{kl}Y_k$ on $V_k\cap V_l$) and $\beta_i(p,t):=tX_i(p)$ (resp.
$\tilde{\beta}_k(p,t):=tY_k(p)$). Since $B$ and $\tilde{B}$ both represent ${\cal F}$,
for every $i,k$ with $U_i\cap V_k\neq \emptyset$, there is a nowhere vanishing holomorphic function $\lambda_{ik}$ on $(U_i\cap V_k)\sm {\rm sing}({\cal F})$ such that
$X_i(p)=\lambda_{ik}(p)Y_k(p)$ for all $p\in (U_i\cap V_k)\sm {\rm sing}({\cal F})$.
By Theorem 1.13 $\lambda_{ik}$ can be extended to $U_i\cap V_k$. Note that the
extended function cannot vanish at all, since its only possible zero set is
$U_i\cap V_k\cap {\rm sing}({\cal F})$ which has codimension $>$1. Now define
$\psi :B\rightarrow \tilde{B}$ by mapping the class of $(p,t)\in U_i\times \CC$
to the class of $(p, t\lambda_{ik}(p))\in V_k\times \CC$. It is quite easy to
see that $\psi$ defines an isomorphism of line bundles with $\beta =\tilde{\beta}\circ \psi$.

Conversely, let $\beta :B\rightarrow TM$ represents ${\cal F}$. If $\tilde{B}$
is any line bundle over $M$ with bundle map $\tilde{\beta } : \tilde {B} \rightarrow TM$
isomorphic to $B$ by $\psi : B\rightarrow \tilde {B}$ such that $\beta =\tilde {\beta } \circ \psi$,
then $\tilde{\beta}:\tilde{B}\rightarrow TM$ also represents ${\cal F}$.\\ \\
{\bf 1.17 Line Bundles and $H^1(M,{\cal O}^\ast )$}   Every holomorphic line bundle $B$ on a
complex manifold $M$ is uniquely determined by an open covering ${\cal U}=\{ U_i\} _{i\in I}$
of $M$ and a family $\{ \xi_{ij}\}_{i,j\in I} $ of non-vanishing holomorphic functions on each
$U_i\cap U_j\neq \emptyset $ satisfying the cocycle relation $\xi_{ij}=\xi_{ik}\cdot \xi_{kj}$
on $U_i\cap U_j\cap U_k$. Thus $\{ \xi_{ij}\} $ may be regarded as an element
of $Z^1({\cal U},{\cal O}^\ast )$, the group of \v{C}ech 1-cocycles with coefficients
in the sheaf of non-vanishing holomorphic functions on $M$, with respect to the
covering ${\cal U}$. Let $B'$ be another line bundle on $M$ defined
by cocycle $\{ \eta_{ij}\} $. It is not difficult to see that $B'$ is isomorphic to
$B$ if and only if there exist non-vanishing holomorphic functions $f_i:U_i\rightarrow \CC^\ast $
such that $\eta_{ij}=(f_i/f_j)\xi_{ij}$ on $U_i\cap U_j$.
Interpreting $\{ \xi_{ij}\} $ and $\{ \eta_{ij}\} $ as elements of $Z^1({\cal U},{\cal O}^\ast )$,
the last condition may be written as $\{ \xi_{ij}/\eta_{ij}\} \in B^1({\cal U},{\cal O}^\ast)$,
the group of \v{C}ech 1-coboundaries. We conclude that {\it two holomorphic 
line bundles on M are isomorphic if and only if they represent the same element of the \v{C}ech cohomology group} $H^1(M,{\cal O}^\ast ):=Z^1(M,{\cal O}^\ast )/B^1(M,{\cal O}^\ast )$.\\ \\
{\bf 1.18 Line Bundles on $\CC \PP^2$}   Most of the results presented here are true
for $\CC \PP^n,\ n\geq 2$. However, we only treat the case $n=2$ for simplicity
of exposition.

Consider the short exact sequence of sheaves
$$0 \rightarrow {\ZZ}\stackrel {2\pi i}{\longrightarrow} {\cal O}\stackrel {\rm exp}{\longrightarrow} {\cal O}^\ast \rightarrow 0$$
on $M$. From this sequence we obtain the long exact sequence of cohomology groups
\begin{equation}
\cdots \rightarrow H^1(M,{\cal O})\rightarrow H^1(M,{\cal O}^\ast )\stackrel{c_1}{\longrightarrow} H^2(M,{\ZZ})\rightarrow H^2(M,{\cal O})\rightarrow \cdots
\end{equation}
For every line bundle $B\in H^1(M,{\cal O}^\ast ),\ c_1(B)\in H^2(M,{\rm {\ZZ}})$
is called the {\it first Chern class} of $B$. Explicitly, if $B=[\{ \xi_{ij}\} ]$
is defined for a covering ${\cal U}=\{ U_i\} $ for which every $U_i$ is connected,
every $U_i\cap U_j$ is simply-connected, and every $U_i\cap U_j\cap U_k$ is connected,
then $c_1(B)=[\{ c_{ijk}\} ]$, where $c_{ijk}:=1/{(2\pi \sqrt{-1})}\{ \log \xi_{jk}-\log \xi_{ik}+\log \xi_{ij}\} $
and the branches of logarithms are arbitrarily chosen \cite{Kod}.

There are two basic facts about the sequence (1) in the case $M=\CC \PP^2$.
First, for every $i\geq 1$ we have $H^i(\CC \PP^2,{\cal O})=0$ as a consequence of the Hodge Decomposition Theorem \cite{G-H}. Second, $H^2(\CC \PP^2,{\ZZ})\simeq {\ZZ}$ \cite{G-H}.
Therefore, $c_1$ is an isomorphism in (1), $H^1(\CC \PP^2,{\cal O}^\ast )$ is the infinite cyclic group $\ZZ$, and every line bundle on $\CC \PP^2$ is determined up to isomorphism by its first Chern class (cf. 1.17).

Now consider $\CC \PP^2$ with affine charts $\{ (U_i,\phi_i^{-1})\} _{i=0,1,2}$ , as in
1.7. For every integer $n$, define a line bundle $B(n)$ on $\CC \PP^2$ whose
cocycle $\{ \xi _{ij}\}_{i,j=0,1,2}$ is given by
$$\xi_{ij}:U_i\cap U_j\rightarrow \CC^\ast ,\ \ \ i,j=0,1,2$$
\begin{equation}
\xi_{ij}[x_0,x_1,x_2]:={ \left (\frac {x_j}{x_i} \right )}^n.
\end{equation}
\setcounter{equation}{0}
It is not difficult to check that $B(-1)$, called the {\it canonical line bundle} over $\CC \PP^2$, and its dual bundle $B(1)$ are both generators for the infinite cyclic group $H^1(\CC \PP^2,{\cal O}^\ast )$. By a standard convention we define $c_1(B(1))=1$. It follows that $c_1(B(n))=n$ for all integers $n$.

Adding up the above remarks, it follows that {\it every holomorphic line bundle $B$ on $\CC \PP^2$ is isomorphic to $B(n)$, where $n=c_1(B)$}.

Now suppose that ${\cal F}$ is a SHFC on $\CC \PP^2$, and let $\beta :B\rightarrow T\CC \PP^2$
be its associated bundle map (see Theorem 1.15). Choose $n=c_1(B)$ so that
$B$ is isomorphic to $B(n)$ by some $\psi :B\rightarrow B(n)$. By 1.16(c), if $\tilde{\beta}$ is defined
by $\tilde{\beta}:= \beta \circ \psi^{-1}$, then $\tilde{\beta}:B(n)\rightarrow T\CC \PP^2$ is a
bundle map which also represents ${\cal F}$.\\ \\
{\bf 1.19 Explicit Form of SHFC's on $\CC \PP^2$} \ In view of the remarks in the last section, we are now
going to determine the explicit form of every holomorphic bundle map $\beta :B(-n+1)\rightarrow T\CC \PP^2$ for every integer $n$. (We choose $-n+1$ instead of $n$ just to make later formulations easier.)

Once again consider $\CC \PP^2$ equipped with three affine charts $(U_i,\phi_i^{-1}), i=0,1,2$.
Restricting $\beta $ to $U_i\times \CC$, there exists a holomorphic vector
field $X_i$ on $U_i$ representing the action of $\beta_i :U_i\times \CC\rightarrow TU_i \simeq U_i \times \CC^2$. In other words, $\beta_i(p,t)=(p, t X_i(p))$. Let $X_0:=f(x,y)\partial /\partial x+g(x,y)\partial /\partial y$ and $X_1:=\tilde{f}(u,v)
\partial /\partial u+\tilde{g}(u,v)\partial /\partial v$, where $f,g,\tilde{f},\tilde{g}$
are holomorphic on $\CC^2$. Set $\tilde{U}_0:=U_0\sm \{(x,y):x=0\} $ and
$\tilde{U}_1:=U_1\sm \{(u,v):u=0\} $. We have the following commutative diagram:\\
$$\begin{CD}
\tilde{U}_0 \times \CC @> \xi >> \tilde{U}_1 \times \CC\\
@VV{\beta_0}V  @VV{\beta_1}V \\
T\tilde{U}_0 \simeq \tilde{U}_0 \times \CC^2 @> (\phi_{01})_\ast>> T\tilde{U}_1\simeq \tilde{U}_1 \times \CC^2
\end{CD}$$
where $\xi (x,y,t):=(u,v,t / \xi_{01}(x,y))$. Note that by 1.18(2), $\xi_{01}(x,y)=\xi_{01}[1,x,y]=x^{-n+1}=u^{n-1}$ so that $\xi (x,y,t)=(u,v,u^{-n+1}t)$. It follows that $(\phi_{01})_{\ast}X_0=u^{-n+1} X_1$, or 
\begin{equation}
\begin{array}{l}
\tilde{f}(u,v)=\displaystyle{-u^{n+1}f(\frac{1}{u} , \frac{v}{u})},\vspace{0.1in}\\
\tilde{g}(u,v)=\displaystyle{-u^{n}[vf(\frac{1}{u} , \frac{v}{u})-g(\frac{1}{u} , \frac{v}{u})]}.
\end{array}
\end{equation}
Since $\tilde{f}$ and $\tilde{g}$ are holomorphic on $\CC^2$, (1) shows
that if $n\leq -1$, then $f$ and $g$ both vanish. So let us assume that $n\geq 0$.

Suppose that $f=\sum _{k=0} ^{\infty}f_k$ and $ g=\sum _{k=0}^{\infty}g_k$, where $f_k$ and $g_k$ are the homogeneous parts of degree $k$ of the power series expansions of $f$ and $g$. It follows then from (1) that $f_k\equiv g_k\equiv 0$ for $k\geq n+2$, and
\begin{equation}
vf_{n+1}(1,v)-g_{n+1}(1,v) \equiv 0.
\end{equation}
Coming back to the affine chart $(x,y)\in U_0$, we obtain from (2) that
\begin{equation}
yf_{n+1}(x,y)-xg_{n+1}(x,y)\equiv 0.
\end{equation}
It is easy to see that there are no other restrictions on these homogeneous polynomials. Thus we have proved the following \\ \\
{\bf 1.20 Theorem}  {\it Let $\beta:B(-n+1)\rightarrow T\CC \PP^2$ be a holomorphic bundle map. If $n\leq -1$, then $\beta \equiv 0$. If $n\geq 0$, then in the affine chart $(x,y)\in U_0$, $\beta$ is given by a polynomial vector field of the form
$$\sum_{k=0}^{n+1} f_k\frac{\partial}{\partial x}+\sum_{k=0}^{n+1}g_k\frac{\partial}{\partial y},$$
where $f_k$ and $g_k$ are homogeneous polynomials of degree k, and $yf_{n+1}-xg_{n+1}\equiv 0$}. $\hfill \Box $\\ 

Now let ${\cal F}$ be a SHFC on $\CC \PP^2$ and let $\beta : B \rightarrow T\CC \PP ^2$ represent ${\cal F}$. Put $n=-c_1(B)+1$. Then, by the above theorem ${\cal F}$ is induced by a polynomial $1$-form 
$$\omega = [ (\sum_{k=0}^n f_k) + xh ] \, dy- [ (\sum_{k=0}^n g_k) +yh] \, dx,$$
in the affine chart $U_0$, where $h=f_{n+1}/x = g_{n+1}/y$ is a homogeneous polynomial of degree $n$ or $h \equiv 0$. Let us assume for a moment that the latter happens, i.e., $f_{n+1} \equiv g_{n+1} \equiv 0$. Then, rewriting 1.19(1) gives us
\setcounter{equation}{0}
\begin{equation}
\begin{array}{l}
\displaystyle{\tilde{f}(u,v)=-u^{n+1} \sum_{k=0}^n u^{-k} f_k(1,v)},\vspace{0.08in}\\
\displaystyle{\tilde{g}(u,v)=-u^{n} \sum_{k=0}^n u^{-k} [v f_k(1,v)-g_k(1,v)]}.
\end{array}
\end{equation}
This shows that $v f_n(1,v)-g_n(1,v) \not \equiv 0$ since otherwise $\tilde{f}$ and $\tilde{g}$ would have a common factor $u$, meaning that $\beta$ would vanish on the entire line $u=0$ (contradicting the fact that $\beta$ must vanish at finitely many points). We summarize the above observations in the following \\ \\
\setcounter{equation}{0}
{\bf 1.21 Corollary}  {\it Let ${\cal F}$ be a SHFC on $\CC \PP^2$ and $\beta:B
\rightarrow T\CC \PP^2$ be any holomorphic bundle map representing ${\cal F}$. Then $c_1(B)\leq 1$. If $n=-c_1(B)+1$, then ${\cal F}$ is given by a unique (up to a multiplicative constant) polynomial 1-form  
\begin{equation}
\omega =(f+xh) \, dy - (g+yh) \, dx
\end{equation}
in the affine chart $(x,y)\in U_0$. Here} 
\begin{enumerate}
\item[$\bullet$]
{\it $f=\sum_{k=0}^n f_k$ and $g=\sum_{k=0}^n g_k$, with $f_k$ and $g_k$ being homogeneous polynomials of degree $k$,}
\item[$\bullet$]
{\it either $h$ is a non-zero homogeneous polynomial of degree n, or if $h\equiv 0$, then $yf_n-xg_n \not \equiv 0$,}
\item[$\bullet$]
{\it two polynomials $f+xh$ and $g+yh$ have no common factor.} $\hfill \Box$\\
\end{enumerate}
{\large \bf Geometric Degree of a SHFC on $\CC \PP^2$}\vspace{4mm}\\
Here we show that the ``cohomological degree'' $n=-c_1(B)+1$ of a SHFC ${\cal F}$ given by Corollary 1.21 coincides with a geometric invariant which we will call the ``geometric degree'' of ${\cal F}$. Roughly speaking, the geometric degree of ${\cal F}$ is the number of points a generic projective line is tangent to the leaves of ${\cal F}$. This notion allows us to give a stratification of the space of all SHFC's on $\CC \PP^2$.\\ \\
{\bf 1.22 Definition of Geometric Degree} \ \ Let ${\cal F}$ be a SHFC on $\CC \PP^2$ and $L$ be any projective line such that $L \sm {\rm sing} ({\cal F})$
is not a leaf of ${\cal F}$. A point $p\in L$ is called a {\it tangency point}
of ${\cal F}$ and $L$ if either $p\in {\rm sing}({\cal F})$, or $p\not \in$ sing$({\cal F})$ and the tangent lines to $L$
and ${\cal L}_p$ at $p$ coincide.

Let ${\cal F}$ be of the form ${\cal F}_\omega :\{\omega =Pdy-Qdx=0 \}$
in the affine chart $(x,y)\in U_0$, and $L\cap U_0$ be parametrized by
$\ell(T)=(x_0+aT,y_0+bT)$, where $p=(x_0,y_0)=\ell(0)$. Then it is clear
that $p$ is a tangency point of ${\cal F}$ and $L$ if and only if $T=0$ is a root of the polynomial $T\mapsto b P(\ell(T))-a Q(\ell(T))$. The {\it order} of tangency of ${\cal F}$ and $L$ at $p$ is defined to be the multiplicity
of $T=0$ as a root of this polynomial. Define
$$m({\cal F} , L):=\sum _p ( {\rm order \ of \ tangency \ of } \ {\cal F} \ {\rm and } \ L \ {\rm at} \ p ),$$
where the (finite) sum is taken over all the tangency points.

It is quite elementary to show that $m({\cal F}, L)$ does not depend on
$L$ (as long as $L \sm {\rm sing} ({\cal F})$ is not a leaf), so that we can call it the {\it geometric degree} of ${\cal F}$.\\ \\
{\bf 1.23 Theorem} \ \ {\it Let ${\cal F}$ be a SHFC on $\CC \PP^2$ and $n=-c_1(B)+1$, where $B$ is the line bundle associated with ${\cal F}$. Let $L$ be any projective line such that $L \sm {\rm sing} ({\cal F})$ is not a leaf. Then $m({\cal F}, L)=n$. In particular, a SHFC ${\cal F}$ has geometric degree $n$ if and only if ${\cal F}$ is induced by a 1-form $\omega$ as in 1.21(1) in which $h$ is a non-zero homogeneous polynomial of degree $n$, or $h\equiv 0$ and $yf_n-xg_n \not \equiv 0$}. $\hfill \Box$\\ \\
{\bf Proof.} Let $n=-c_1(B)+1$ so that ${\cal F}$ is induced by a $1$-form $\omega$ as in 1.21(1). Take a projective line $L$ such that $L \sm {\rm sing} ({\cal F})$ is not a leaf. Since the normal form 1.21(1) is invariant under projective transformations, we may assume that $L$ is the $x$-axis. By 1.21(1), $(x,0)\in L\cap U_0$ is a tangency point if and only if $g(x,0)=\sum_{k=0}^n g_k(x,0)=0$. As for the point at infinity for $L$, consider the affine chart $(u,v)=(1/x, y/x)\in U_1$ in which $L$ is given by the line $\{ v=0 \}$. By 1.19(1), the foliation is described by the polynomial 1-form $\omega'=\tilde{f}\, dv -\tilde{g}\, du$, where
\setcounter{equation}{0}
\begin{equation}
\begin{array}{l}
\displaystyle{\tilde{f}=-\sum_{k=0}^n u^{n+1-k}f_k(1,v)-h(1,v),}
\vspace{0.08in}\\
\displaystyle{\tilde{g}=-\sum_{k=0}^n u^{n-k}[vf_k(1,v)-g_k(1,v)].}
\end{array}
\end{equation}
This shows that $L$ has a tangency point at infinity if and only if $u=0$ is a root of $\tilde{g}(u,0)=0$, or $g_n(1,0)=0$ by 1.23(1). To prove the theorem, we distinguish two cases:

(i) The polynomial $x \mapsto g_n(x,0)$ is not identically zero. In this case, $\sum_{k=0}^n g_k(x,0)$ is a polynomial of degree $n$ in $x$, so there are exactly $n$ finite tangency points on $L$ counting multiplicities. Note that the point at infinity for $L$ is not a tangency point since $g_n(1,0)\neq 0$. So in this case, $m({\cal F}, L)=n$.

(ii) There is a largest $0\leq j<n$ such that $x \mapsto g_j(x,0)$ is not identically zero (otherwise $\sum_{k=0}^n g_k(x,0)$ would be everywhere zero, so $L \sm {\rm sing} ({\cal F})$ would be a leaf). This means that $g(x,0)=0$ has exactly $j$ roots counting multiplicities. In this case, the point at infinity for $L$ is a tangency point of order $n-j$. In fact, 1.23(1) shows that $\tilde{g}(u,0)=-\sum_{k=0}^n u^{n-k}g_k(1,0)=-\sum_{k=0}^j u^{n-k}g_k(1,0)$, which has a root of multiplicity $n-j$ at $u=0$. Thus, there are $n$ tangency points on $L$ altogether, so again $m({\cal F}, L)=n$. $\hfill \Box$\\ \\
As an example, 1.11(5) shows that {\it ${\cal F}$ is induced by a holomorphic
vector field on $\CC \PP^2$ if and only if the geometric degree of ${\cal F}$ is $\leq 1$.}\\

The set of all SHFC's on $\CC \PP^2$ of geometric degree $n$ is denoted by ${\cal D}_n$. Each ${\cal D}_n$ is topologized in the natural way: a neighborhood of
${\cal F} \in {\cal D}_n$ consists of all foliations of geometric degree $n$ whose defining
polynomials have coefficients close to that of ${\cal F}$, up to multiplication
by a non-zero constant. To be more accurate, consider the complex linear space of all
polynomial 1-forms $\omega $ as in 1.21(1). By 1.16(b), $\omega $
and $\omega'$ define the same foliation if and only if there exists a non-zero constant
$\lambda $ such that $\omega '=\lambda \omega $. Therefore ${\cal D}_n$ may be considered
as an open subset of the complex projective space $\CC \PP^N$, where $N$
is the dimension of the above linear space minus one, that is $N=2 \sum ^{n+1}_{k=1} k+(n+1)-1=n^2+4n+2$.
It is not difficult to check that ${\cal D}_n$ is connected and dense in this
projective space.\\ \\
{\bf 1.24 Corollary} \ \ {\it The set ${\cal D}_n (n\geq 0)$ of all SHFC's of geometric degree $n$ on $\CC \PP^2$ can be identified with an open, connected and dense subset of the complex projective space $\CC \PP^N$, where
$N=n^2+4n+2$. So we can equip ${\cal D}_n$ with the induced topology and a natural Lebesgue measure class}. $\hfill \Box $ \\ \\
The definition of ${\cal D}_n$ allows us to decompose the space ${\cal S}$ of all SHFC's on $\CC \PP^2$ into a disjoint union $\bigcup {\cal D}_n$ and topologize it in a natural way. A subset $\cal U$ of $\cal S$ is open if and only if ${\cal U} \cap {\cal D}_n$ is open for every $n$. Hence in this topology every ${\cal D}_n$ is a connected component of $\cal S$. Similarly, $\cal S$ inherits a natural Lebesgue measure class: A set ${\cal U}\subset {\cal S}$ has measure zero if and only if ${\cal U}\cap {\cal D}_n$ has measure zero in ${\cal D}_n$.\\ \\  
\vspace {4 mm}
{\large \bf Line at Infinity as a Leaf}\\
\noindent
Let us find conditions on a 1-form $\omega $ which guarantee that
the line at infinity $L_0= \CC \PP^2 \sm U_0 $ with singular points of ${\cal F}_{\omega}$ deleted is a leaf. Consider a SHFC ${\cal F}\in {\cal D}_n$ induced by a polynomial 1-form $\omega$ as 1.21(1): $\omega =(f+xh)dy-(g+yh)dx$, where $ f= \sum ^n_{k=0} f_k \ , \ g=\sum ^n_{k=0} g_k$ 
and $h$ is either a non-zero homogeneous polynomial of degree $n$, or $h\equiv 0$ but $yf_n-xg_n \not \equiv 0$.\\ \\
{\bf 1.25 Theorem } \ {\it The line at infinity $L_0$ with singular points of
${\cal F} \in {\cal D}_n$ deleted is a leaf of ${\cal F}$ if and only if $h\equiv 0$.}\\ \\
{\bf Proof.} This follows easily from the proof of Theorem 1.23. In fact, $L_0 \sm {\rm sing}({\cal F})$ is a leaf if and only if the line $\{ u=0\}$ is a solution of $\omega'=\tilde{f} dv -\tilde{g}du=0$ in 1.23. By 1.23(1), this happens if and only if $h(1,v) \equiv 0$. Since $h$ is a homogeneous polynomial, the latter condition is equivalent to $h\equiv 0$. $\hfill \Box$\\ \\
{\bf 1.26 Remark}\ Here is an alternative notation for the polynomial 1-forms which will be used in many subsequent discussions. Let 
$$\omega=(f+xh)dy-(g+yh)dx=Pdy-Qdx,$$
as in 1.21(1). 
Define 
$$R(x,y)=yP(x,y)-xQ(x,y)=yf(x,y) -xg(x,y)$$
as in 1.9(1), and note that deg $R \leq n+1$. Then ${\cal F}|_{U_1}$ is given by $\{ \omega '=0 \}$, where
\setcounter{equation}{0}
$$\omega '(u,v)=u^{k}P(\frac{1}{u} , \frac{v}{u})dv-u^{k} R(\frac{1}{u} , \frac{v}{u})du,$$
and, as in 1.9(1), $k$ is the least positive integer which makes $\omega '$ a
polynomial 1-form. (Note that this representation should be identical to 1.23(1) up to a multiplicative constant.) If $h \not \equiv 0$, then deg $P=n+1$ and so $k=n+1$. On the other hand, if $h \equiv 0$, then $yf_n-xg_n \not \equiv 0$ by Theorem 1.23 which means deg $R= n+1$. So again we have $k=n+1$.\\ 

For a fixed ${\cal F}$ we denote $L_0\sm {\rm sing}({\cal F})$ by ${\cal L}_\infty$, if it is indeed a leaf of ${\cal F}$. Frequently, we refer to ${\cal L}_\infty$ as the {\it leaf at infinity}.

As can be seen from the above theorem, for a foliation ${\cal F} \in {\cal D}_n$ the line at infinity $L_0 \sm {\rm sing}({\cal F})$ is unlikely to be a leaf since this condition is equivalent to vanishing of a homogeneous polynomial. This is a result of the way we topologized the space $\cal S$ of all SHFC's on $\CC \PP^2$ using the topologies on the ${\cal D}_n$. The decomposition $\bigcup {\cal D}_n$ is quite natural from the geometric point of view; however, it leads to a rather peculiar condition on the polynomials describing the associated $1$-form (Theorem 1.23). The situation can be changed in a delicate way by choosing a different decomposition ${\cal S}=\bigcup {\cal A}_n$ which is more natural from the point of view of differential equations in $\CC^2$ but has no longer an intrinsic geometric meaning. Elements of ${\cal A}_n$ are simply determined by the maximum degree of their defining polynomials.\\ \\  
{\bf 1.27 Definition} \ \ Fix the affine chart $(x,y)\in U_0$ and let $n\geq 0$. A SHFC ${\cal F}$ is said to belong to the class ${\cal A}_n$ if it is induced by a polynomial 1-form $\omega = Pdy-Qdx$ with max $\{ {\rm deg} P, {\rm deg} Q \}=n$ and $P, \, Q$ relatively prime. The number $n$ is called the {\it affine degree} of ${\cal F}$ (with respect to $U_0$).\\ \\
Note that ${\cal A}_n$ is well-defined since if ${\cal F}_\omega ={\cal F}_ {\omega '}$, then $\omega '=\lambda \omega$ for some non-zero constant $\lambda $.
It is important to realize that unlike the condition ${\cal F} \in {\cal D}_n$ (normal form 1.21(1)), whether or not ${\cal F} \in {\cal A}_n$ strongly depends on the choice of a particular affine coordinate system, and that is why we call $n$ the ``affine degree.''  

Consider the complex linear space of all polynomial 1-forms $ \omega =Pdy-Qdx $ with
max $\{{\rm deg} P, {\rm deg} Q \} \leq n$, which has dimension $(n+1)(n+2).$ Then, as in the case of ${\cal D}_n$, one has\\ \\
{\bf 1.28 Corollary} \ \ {\it The set ${\cal A}_n (n \geq 0)$ of all SHFC's of affine degree $n$ on $\CC \PP^2$ can be identified with 
an open, connected and dense subset of the complex projective space $\CC \PP^N$, where $N=n^2+3n+1$. So we can equip ${\cal A}_n$ with the induced topology and a natural Lebesgue measure class. $\hfill \Box $}\\

Using the decomposition of $\cal S$ into the disjoint union of the ${\cal A}_n$, we can define a new topology and measure class on $\cal S$ in the same way we did using the ${\cal D}_n$ (see the remarks after Corollary 1.24). In this new topology, each class ${\cal A}_n$ turns into a connected component of $\cal S$. The two topologies and measure classes are significantly different.  As a first indication of this difference, let ${\cal F}: \{\omega =Pdy-Qdx=0 \} \in {\cal A}_n$ and decompose $P=\sum ^n_{k=0} P_k$ and $Q=\sum ^n_{k=0} Q_k$ into the sum of homogeneous polynomials $P_k$ and $Q_k$ of degree $k$. Then it easily follows from Theorems 1.23 and 1.25 that\\ \\
{\bf 1.29 Corollary} \ \ {\it The line at infinity $L_0 \sm$ {\rm sing} $({\cal F})$ is a leaf of ${\cal F} : \{ Pdy-Qdx=0 \} \in {\cal A}_n$ if and only if $yP_n-xQ_n \not \equiv 0$.} $\hfill \Box$ \\ \\
One concludes that in ${\cal A}_n $ it is very likely to have $L_0 \sm$ {\rm sing} $({\cal F})$ as a leaf, contrary to what we observed before in ${\cal D}_n$. \\ \\
{\bf 1.30 Corollary}\ \ {\it Fix the affine chart $(x,y)\in U_0$ and a SHFC ${\cal F} \in {\cal A}_n$.}
\begin{enumerate}
\item[$\bullet$]
{\it If the line at infinity $L_0 \sm \mbox{sing}({\cal F})$ is a leaf, then ${\cal F} \in {\cal A}_n \cap {\cal D}_n$ so that }
\begin{center}
{\it affine degree of ${\cal F}=$ geometric degree of ${\cal F}.$}
\end{center} 
\item[$\bullet$] 
{\it Otherwise, ${\cal F} \in {\cal A}_n \cap {\cal D}_{n-1}$ so that}
\begin{center}
{\it affine degree of ${\cal F}=($geometric degree of ${\cal F}) + 1.$} \vspace{0.1in}
\end{center} 
\end{enumerate}

The main reason for the contrast between ${\cal A}_n$ and ${\cal D}_n$ is the fact that ${\rm dim} {\cal A}_n <{ \rm dim} {\cal D}_n < {\rm dim } {\cal A}_{n+1}$. In fact, it is not hard to see that ${\cal D}_n \subset {\cal A}_n \cup {\cal A}_{n+1}$ and ${\cal A}_n \subset {\cal D}_{n-1} \cup {\cal D}_n$. Fig.\ 3 illustrates the first property while Fig. 4 is a schematic diagram of the set-theoretic relations between these classes.\\ \\
\realfig{figad}{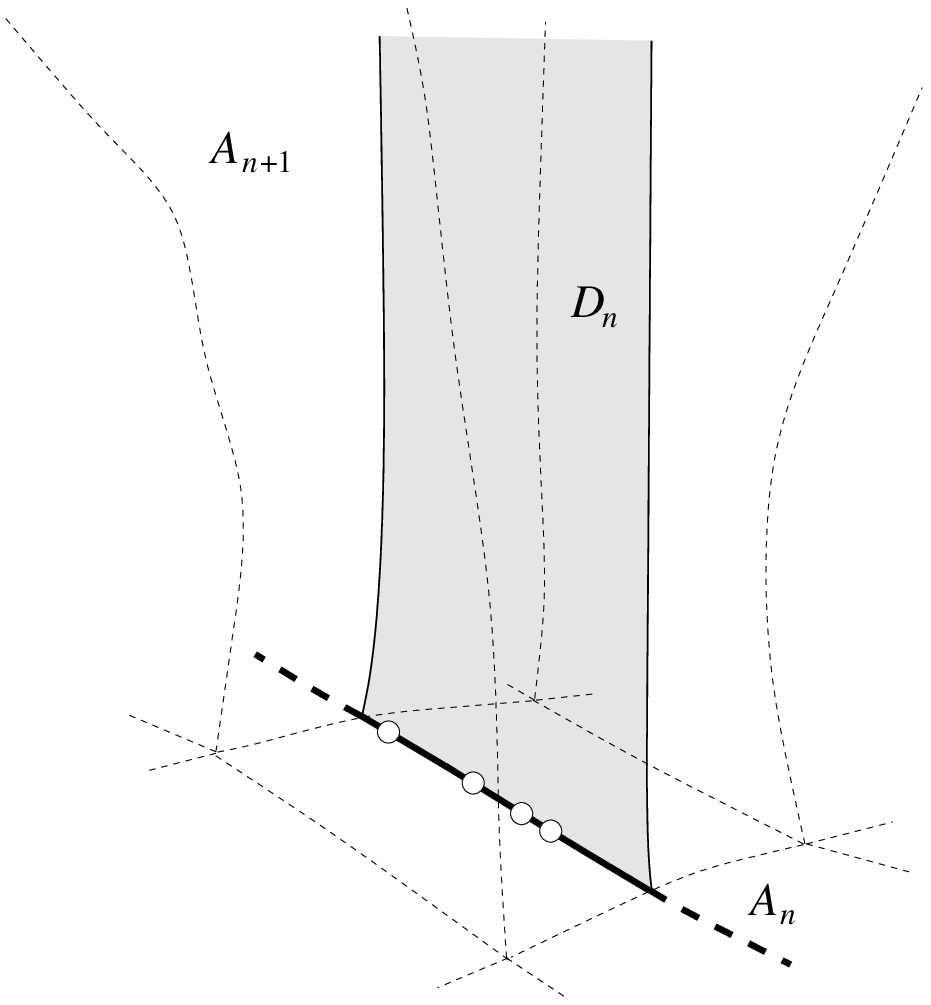}{${\rm dim} {\cal A}_n <{ \rm dim} {\cal D}_n < {\rm dim } {\cal A}_{n+1}$}{7cm}
\realfig{figthesis3}{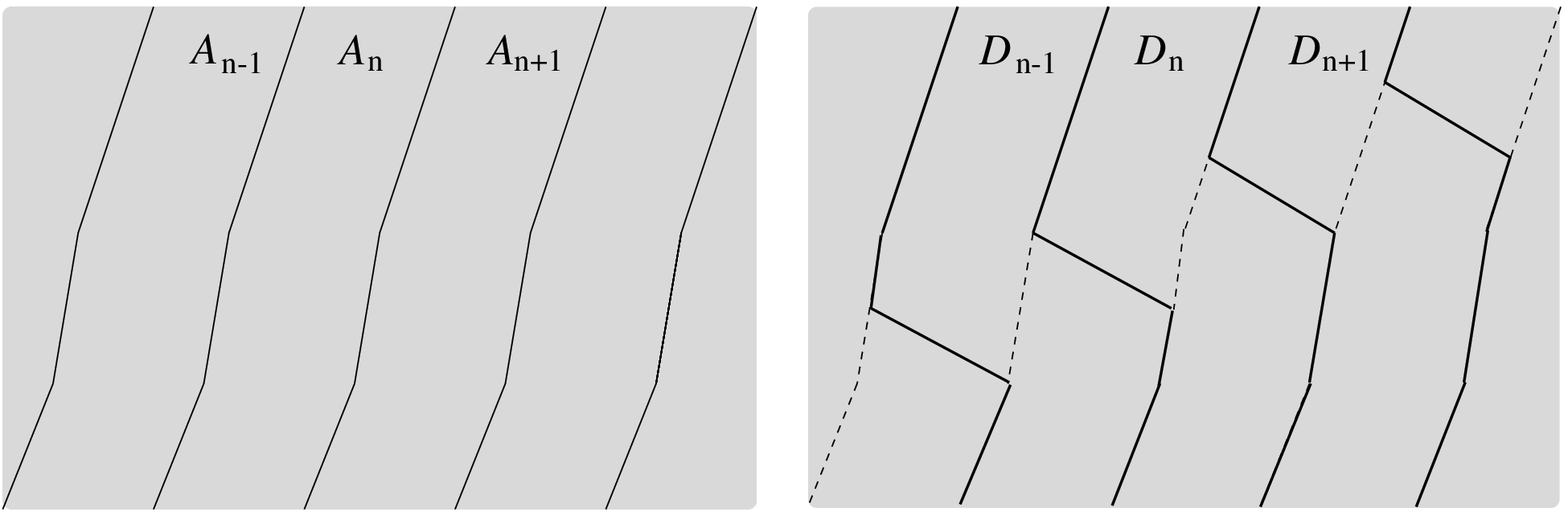}{\sl Set-theoretic relations between $\{ {\cal D}_n \}$ and $\{ {\cal A}_n \}$}{10cm}
\noindent
{\bf 1.31 Examples} The two SHFC's
$${\cal F}_1 : \{ xdy-ydx=0 \} \hspace{2cm} {\cal F}_2 : \{ ydy-xdx=0 \}$$ 
both belong to ${\cal A}_1$ so they both have affine degree $1$. However, ${\cal F}_1$ belongs to ${\cal D}_0$ hence has geometric degree $0$, while ${\cal F}_2$ belongs to ${\cal D}_1$ and so it has geometric degree $1$. Note that the line at infinity is {\it not} a leaf of ${\cal F}_1$ but it {\it is} a leaf of ${\cal F}_2$. 

As another example, let us illustrate how the topologies coming from the two decompositions ${\cal S}=\bigcup {\cal A}_n=\bigcup {\cal D}_n$ are different. Consider the two SHFC's 
$${\cal F} : \{ x^2 dy-y^2 dx=0 \},$$
$${\cal F}_{\varepsilon} : \{ (x^2+\varepsilon xy^2) dy-(y^2+\varepsilon y^3) dx=0 \}.$$
As $\varepsilon \rightarrow 0$, ${\cal F}_{\varepsilon} \rightarrow {\cal F}$ in the topology induced by $\bigcup {\cal D}_n$ but ${\cal F}_{\varepsilon}$ does not converge in the topology induced by $\bigcup {\cal A}_n$.\\ \\
{\bf 1.32 A Complex One-Dimensional Analogy}\ The following simple example may help understand the difference between the two decompositions $\{ {\cal D}_n \}$ and $\{ {\cal A}_n \}$: Let $\bf S$ be the space of all linear conjugacy classes of complex polynomial maps in $\CC$ of degree at most $4$ which are tangent to the identity map at the origin. This space can be naturally decomposed by the order of tangency near the fixed point: For $0\leq n \leq 2$, consider the sets ${\bf D}_n$ of conjugacy classes of normalized polynomials as follows:
$$\begin{array}{rl}
{\bf D}_0= & \langle z\mapsto z+z^4 \rangle \simeq \mbox{point} \\
{\bf D}_1= & \langle z\mapsto z+z^3+az^4 \rangle \simeq \CC \\
{\bf D}_2= & \langle z\mapsto z+z^2+az^3+bz^4 \rangle \simeq \CC^2.
\end{array}$$
Clearly ${\bf S}=\bigcup_{n=0}^2 {\bf D}_n$ and this decomposition induces a topology $\tau _{\bf D}$ and a measure class $\mu_{\bf D}$ on $\bf S$. On the other hand, one can consider the following conjugacy classes determined by the degree of polynomials (i.e., by their behavior near infinity):
$$\begin{array}{rl}
{\bf A}_0= & \langle z\mapsto z+z^2 \rangle \simeq \mbox{point} \\
{\bf A}_1= & \langle z\mapsto z+az^2+z^3 \rangle \simeq \CC \\
{\bf A}_2= & \langle z\mapsto z+az^2+bz^3+z^4 \rangle \simeq \CC^2.
\end{array}$$
This gives rise to a second decomposition ${\bf S}=\bigcup_{n=0}^2 {\bf A}_n$ hence a corresponding topology $\tau _{\bf A}$ and measure class $\mu_{\bf A}$ on $\bf S$. One has the relations
$$\begin{array}{ll}
{\bf A}_0 \subset {\bf D}_2 & {\bf D}_0 \subset {\bf A}_2\\
{\bf A}_1 \subset {\bf D}_1 \cup {\bf D}_2 & {\bf D}_1 \subset {\bf A}_1 \cup {\bf A}_2
\end{array}$$
The topologies $\tau _{\bf D}$ and $\tau _{\bf A}$ and measure classes $\mu_{\bf D}$ and $\mu_{\bf A}$ are very different. For example, ${\bf D}_1 \subset {\bf S}$ is an open set in $\tau_{\bf D}$, but it is {\it not} open in $\tau _{\bf A}$ since ${\bf D}_1 \cap {\bf A}_1$ is a single point. Similarly, ${\bf D}_1 \subset {\bf A}_1 \cup {\bf A}_2$ has measure zero with respect to $\mu_{\bf A}$ but this is certainly not true with respect to $\mu_{\bf D}$.\\ \\
{\bf 1.33 Typical Properties}  Certain geometric or dynamical properties often hold for ``most'' and not all SHFC's in ${\cal A}_n$ or ${\cal D}_n$. In these cases, we can use the Lebesgue measure class to make sense of this fact. A property $\cal P$ is said to be {\it typical} for elements of ${\cal A}_n$, or we say that a {\it typical} SHFC in ${\cal A}_n$ satisfies $\cal P$, if $\{ {\cal F}\in {\cal A}_n: {\cal F}\ {\rm does\ not\ satisfy}\ {\cal P} \}$ has Lebesgue measure zero in ${\cal A}_n$. We can define a typical property in ${\cal D}_n$ in a similar way.\\ 

For example, it follows from Theorem 1.25 and Corollary 1.29 that {\it having the line at infinity as a leaf is not typical in ${\cal D}_n$ but it is typical in ${\cal A}_n$}. \\ \\
{\bf 1.34 Definition} \  Let ${\cal F}\in {\cal A}_n$. We say that ${\cal F}$ has {\it Petrovski\u\i-Landis property} if $L_0\sm {\rm sing}({\cal F})$ is a
leaf of ${\cal F}$ and $L_0\cap {\rm sing}({\cal F})$ consists of $(n+1)$ distinct points in the affine chart $(u,v)\in U_1$. The class of all such ${\cal F}$ is denoted by ${\cal A}_n'$.\\

Assuming ${\cal F}:\{ Pdy-Qdx=0\} \in {\cal A}_n'$, one has $yP_n-xQ_n\not \equiv 0$
by Corollary 1.29. On the other hand, using the notation of 1.26, if $R=yP-xQ$, we have deg$R=n+1$ and
$L_0\cap {\rm sing}({\cal F})=\{ (0,v):u^{n+1}\displaystyle{R(\frac{1}{u} , \frac{v}{u})}|_{u=0}=0\} $ in the affine
chart $(u,v)\in U_1$. It follows that $u^{n+1}\displaystyle{R(\frac{1}{u} , \frac{v}{u})}|_{u=0}$ must have $n+1$ distinct
roots in $v$. The above two conditions on $P$ and $Q$ show that\\ \\
{\bf 1.35 Corollary}\ {\it A typical SHFC in ${\cal A}_n$ has the Petrovski\u\i-Landis property }.$\hfill \Box$
\newpage
\vspace*{4 mm}
\noindent
{\large {\bf Chapter 2}}\vspace{4mm}\\
{\Large{\bf The Monodromy Group of a Leaf}}\\ \\ \\ \\ \\ \\ \\ \\ \\ \\ \\ \\
\thispagestyle{plain}
\noindent
Given a SHFC ${\cal F}$ on $\CC \PP^2$ one can study individual leaves as
Riemann surfaces. However, to study the so-called {\it transverse dynamics}
of the foliation one needs a tool to describe rate of convergence or divergence of
nearby leaves. The concept of holonomy, and in particular the monodromy of
a leaf, first introduced by C. Ehresmann, is the most natural and essential tool
for describing the transverse dynamics near the leaf. The point is that the
transverse dynamics of a leaf depends directly on its fundamental
group: the smaller $\pi_1 ({\cal L})$ is , the simpler the behavior of the
leaves near ${\cal L}$ will be.\\ \\
\vspace{4 mm}
{\large \bf Holonomy Mapping and the Monodromy Group}\\
Let ${\cal F}$ be a SHFC on $\CC \PP^2$ and ${\cal L}$ be a non-singular leaf
of ${\cal F}$. Fix $p,q \in {\cal L}$, and consider small sections $\Sigma , \Sigma ' \simeq \DD$
transversal to ${\cal L}$ at $p,q,$ respectively. Let $\gamma :[0,1]\rightarrow {\cal L}$
be a continuous path with $\gamma(0)=p, \gamma (1)=q$. For each $z\in \Sigma $
near $p$ one can ``travel'' on ${\cal L}_z$ ``over'' $\gamma [0,1]$ to reach
$\Sigma '$ at some point $z '$. To be precise, let $\{(U_i , \varphi _i)\}_{0\leq i \leq n}$
be foliation charts of ${\cal F}$ and $0=t_0 <t_1 <\cdots <t_n <t_{n+1}=1$ be
a partition of $[0,1]$ such that if $U_i\cap U_j \neq \emptyset$ then $U_i \cup U_j$
is contained in a foliation chart, and $\gamma [t_i,t_{i+1}] \subset U_i$
for $0\leq i \leq n$. For each $1\leq i \leq n $ choose a section $\Sigma_i\simeq \DD$
transversal to ${\cal L}$ at $\gamma (t_i)$, and let $\Sigma _0=\Sigma $ and
$\Sigma _{n+1}=\Sigma '$ (see Fig. 5).
Then for each $z\in \Sigma_i$ sufficiently close to $\gamma (t_i)$ the plaque of
$U_i$ passing through $z$ meets $\Sigma _{i+1}$ in a unique point $f_i(z)$, and
$z\mapsto f_i(z)$ is holomorphic, with $f_i(\gamma (t_i))=\gamma (t_{i+1})$.
It follows that the composition $f_\gamma :=f_n \circ \cdots \circ f_0$ is defined for
$z\in \Sigma $ near $p$, with $f_\gamma (p)=q$.\\ \\
{\bf 2.1 Definition } \ The mapping $f_\gamma $ is called the {\it holonomy}
associated with $\gamma $.\\

There are several remarks about this mapping which can be checked directly
from the definition (cf. \cite{C-L}). \\ \\


\realfig{figthesis4}{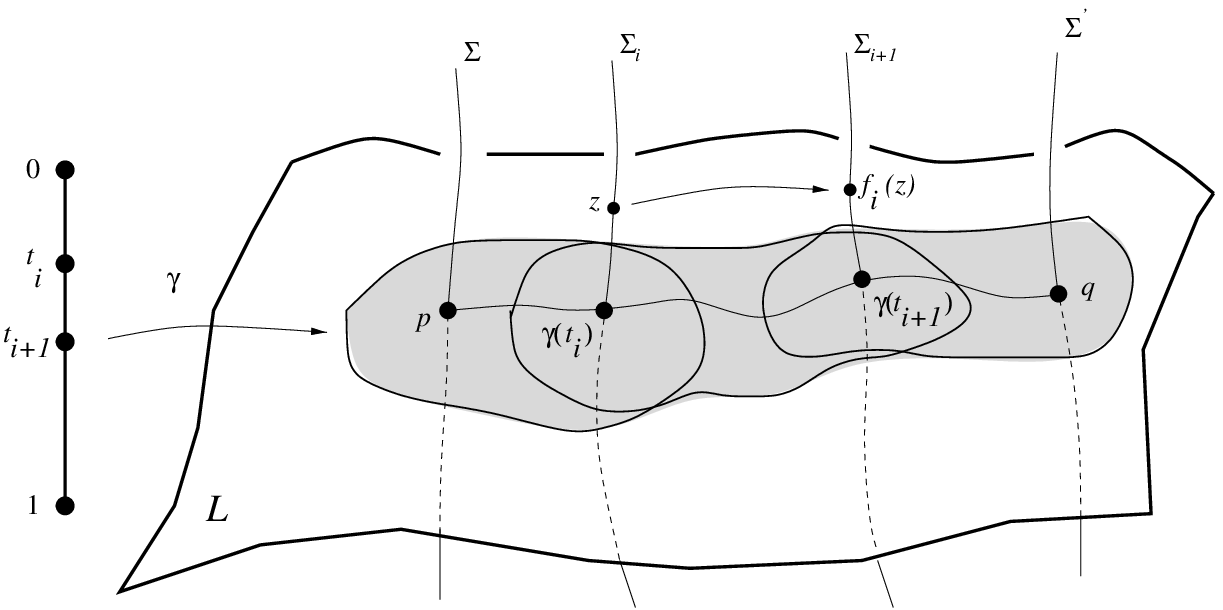}{\sl Definition of holonomy}{11cm}

\noindent
{\bf 2.2 Remarks}

(a) $f_\gamma $ is independent of the chosen transversals $\Sigma _i, 1\leq i \leq n$, and
the foliation charts $U_i$. Hence $\Sigma , \Sigma '$, and $\gamma $ determine its germ at $p$
uniquely.

(b) The germ at $p$ of $f_\gamma $ depends only on the homotopy class of $\gamma $
rel$\{ 0,1 \}$, that is, if $\eta $ is another path joining $p$ and $q$ in
${\cal L}$ which is homotopic to $\gamma $ with $\eta(0)=\gamma(0)$ and $\eta(1)=\gamma(1)$, then the germ of $f_\eta$ at $p$
coincides with that of $f_\gamma $.

(c) If $\gamma ^{-1}(t):=\gamma (1-t)$, then $f_{\gamma^{-1}}={(f_\gamma)}^{-1}$. In particular,
$f_\gamma $ determines the germ of a {\it conformal } mapping $\Sigma \rightarrow \Sigma '$.

(d) Let $\Sigma _1$ and $\Sigma _1'$ be other sections transversal to ${\cal L}$ at $p$ and $q$,
respectively. Let $h:\Sigma \rightarrow \Sigma_1$ and $\tilde {h}:\Sigma ' \rightarrow \Sigma _1'$ be
projections along the plaques of ${\cal F}$ in a neighborhood of $p$ and $q$, respectively.
Then the holonomy $g_\gamma $ with respect to $\Sigma _1,\Sigma _1'$ satisfies
$g_\gamma =\tilde {h}\circ f_\gamma \circ h^{-1}$.\\

Now in the particular case where $p=q$ we have an interesting generalization of the notion of the Poincar\'e
first return map for real vector fields.\\ \\
{\bf 2.3 Definition} \ Let ${\cal L}$ be a non-singular leaf of a SHFC ${\cal F}$ on $\CC \PP^2 , p\in {\cal L}$,
and $\Sigma $ be a section transversal to ${\cal L}$ at $p$. For each $[\gamma ]\in \pi_1 ({\cal L}, p)$, the
holonomy mapping $f_\gamma :\Sigma \rightarrow \Sigma $ is called the {\it monodromy mapping} of
${\cal L}$ associated with $\gamma $ (Fig. 6).\\

Note that by 2.2(b), the germ of $f_\gamma$ at $p$ depends only on the homotopy class $[\gamma ]$. It is quite
easy to see that $[\gamma ]{\stackrel{\Delta}{\longmapsto}} f_\gamma $ is a group
homomorphism from $\pi_1({\cal L},p)$ into the group of germs of biholomorphisms of $\Sigma $
fixing $p$: $\Delta (\gamma \circ \eta )=\Delta (\gamma)\circ \Delta (\eta )$.\\


\realfig{figthesis5}{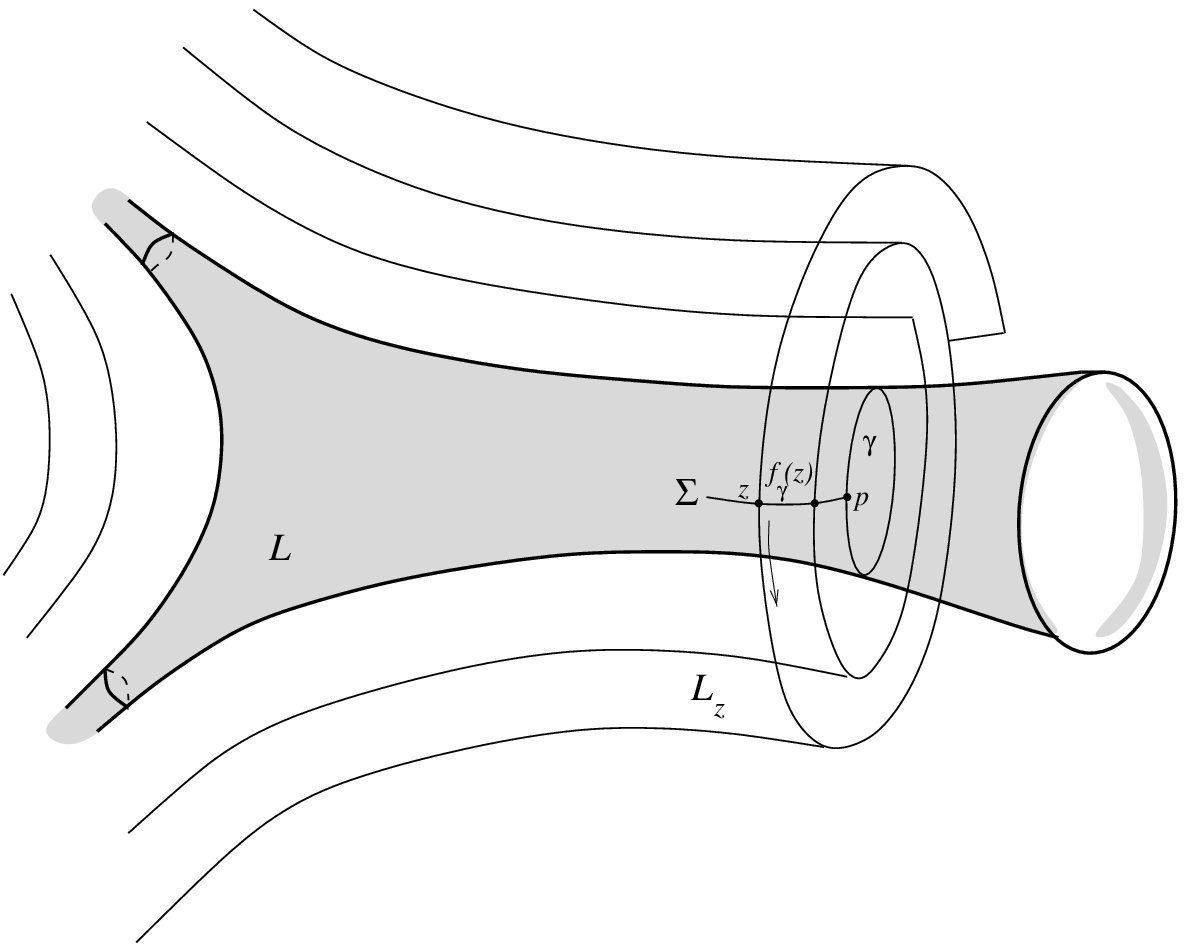}{\sl Monodromy map associated with $\gamma$}{11cm}


\noindent
{\bf 2.4 Remarks}

(a) It should be noted that each $[\gamma ] \in \pi _1 ({\cal L}, p)$
determines only the {\it germ} at $p$ of a biholomorphism $f_\gamma $ of
$\Sigma $, for changing $\gamma $ in its homotopy class results in
changing the domain of definition of $f_\gamma $.

(b) If the transversal $\Sigma $ is replaced by another one $\Sigma _1$,
then by 2.2(d) there exists a germ of biholomorphism $h:\Sigma \rightarrow \Sigma_1$
fixing $p$ such that the monodromy mapping $g_\gamma $ associated with $\gamma $
and $\Sigma _1$ satisfies $g_\gamma =h\circ f_\gamma \circ h^{-1}$. In other words,
the monodromy mapping only depends on $\gamma $ {\it up to conjugacy}. Since
conjugate germs of biholomorphisms fixing $p$ have the same iterative dynamics near $p$,
it is not important which transversal we choose at $p$.\\ \\
{\bf 2.5 Conventions}

(a) We can always fix some arbitrary $p\in {\cal L}$ as the base point of the fundamental group,
since $\pi _1({\cal L}, p) \simeq \pi_1 ({\cal L}, q)$ for every $p,q\in {\cal L}$.

(b) We will always fix some transversal $\Sigma $ at $p$. Moreover, we will choose
a coordinate on $\Sigma $ in which $p=0$. In such a way, every monodromy mapping
$f_\gamma $ may be regarded as an element of ${\rm Bih}_0(\CC)$, the group
of germs at 0 of biholomorphisms $\CC\rightarrow \CC$ fixing the origin.\\ \\
{\bf 2.6 Definition} \ The image under $\Delta $ of $\pi_1 ({\cal L})$ is
called the {\it monodromy group} of the leaf ${\cal L}$, and is denoted by
$G({\cal L})$. Since we always fix the transversals, it can be identified with a subgroup of ${\rm Bih}_0(\CC)$.\\

Given a leaf ${\cal L}$ whose $\pi_1$ is finitely-generated, it is natural to fix
some loops as the generators of $\pi_1({\cal L})$. So we arrive at the
following definition:\\ \\
{\bf 2.7 Definition} \ A non-singular leaf ${\cal L}$ of a SHFC ${\cal F}$
on $\CC \PP^2$ is called a {\it marked leaf} if $\pi_1({\cal L})$ is finitely-generated and a set of loops $\{ \gamma _1, \ldots , \gamma _k \}$ is given
as its generators. Similarly, a finitely-generated subgroup $G \subset {\rm Bih}_0(\CC)$
is called a {\it marked subgroup} if a set of local biholomorphisms $\{f_1,\ldots , f_k\}$
is given as its generators.\\

Clearly, the monodromy group of a marked leaf is a marked subgroup if one
chooses $\{f_{\gamma _1} ,\cdots , f_{\gamma _k} \}$ as its generators.\\ \\
\vspace{4 mm}
{\large \bf Monodromy Pseudo-Group of a Leaf}\\
There is an obvious ambiguity in the domain of definition of an element
of the monodromy group of a leaf. Even if our leaf is marked, so that the
generators of $G({\cal L})$ are fixed, it is not evident what the domain of
definition of an arbitrary $f\in G ({\cal L})$ which does not belong to the generator set is.
On the other hand, it may happen that $f$ can be analytically continued to some
point $z \in \Sigma $, but the result of this continuation differs from the correct
value of the monodromy. For example \cite{I4}, if ${\cal F}$ is induced by a
Hamiltonian form $\omega =dH$, then the monodromy group of ${\cal L}_\infty $
(see 2.14) is Abelian. Hence for any $\gamma \in [\pi _1 ({\cal L}_\infty), \pi_1 ({\cal L}_\infty )]$
the germ of $f_\gamma $ at 0 is the identity, so it can be analytically continued
over $\Sigma $. But for $z \in \Sigma $ sufficiently far away from 0 the
result of continuation of ${\cal L}_z$ over $\gamma $ may differ from $z$,
showing that $f_\gamma (z)$ is not defined for all $z\in \Sigma $.

Since the transverse dynamics of the leaf ${\cal L}$ is reflected in the orbit
of points in $\Sigma $ under the action of $G ({\cal L})$, it is quite natural to be careful about the domains of definitions.\\ \\
{\bf 2.8 Definition} \ Let $G \subset {\rm Bih}_0(\CC)$ be a marked
subgroup with generators $f_1,\ldots , f_k$, all defined on some
domain $\Omega $ around 0. The {\it pseudo-group} $PG$ consists of all
pairs $(f, \Omega _f)$, where $f\in G$ and $\Omega _f$ is a domain on which
$f$ is conformal, with the group operation $(f,\Omega _f)\circ (g, \Omega _g):= (f\circ g, \Omega _{f\circ g})$.
The domain $\Omega _f$ is defined as follows: Let
\setcounter{equation}{0}
\begin{equation}
f=\prod ^N_{i=1} f_{j_i}^{\epsilon_i}\ \ \ \ \ , \ j_i \in \{1,\ldots, k\} \ , \ \epsilon _i \in  \{-1,1\}
\end{equation}
be any representation of $f$ in terms of the generators. Any germ $\prod^n_{i=1} f^{\epsilon_i}_{j_i} , n\leq N$,
is called an {\it intermediate representation} of $f$. The domain $\Omega _{\Pi f}$
associated to the representation (1) is defined as the maximal starlike domain
centered at 0 contained in $\Omega $ on which all the intermediate representations can be
conformally continued, with
$$\prod ^n_{i=1} f^{\epsilon_i}_{j_i} (\Omega _{\Pi f}) \subset \Omega,\ \  {\rm for} \ n\leq N.$$
Finally, $\Omega _f$ is defined to be the union of $\Omega _{\Pi f}$'s for all
possible representations of $f$ of the form (1). It is clear that $f$ is a conformal mapping on $\Omega _f$.\\

By the above definition to each marked leaf ${\cal L}$ there corresponds a {\it monodromy pseudo-group}
$PG({\cal L})$. The construction above allows us to define the {\it orbit} of
$z\in \Sigma $ as $\{f(z) :(f, \Omega _f) \in PG({\cal L})\ {\rm and} \ z\in \Omega _f\}$.
Note that {\it the orbit of z under $PG({\cal L})$ always lies in ${\cal L}_z \cap \Sigma$.}\\ \\
\vspace{4 mm}
{\large \bf Multiplier of a Monodromy Mapping}\\
As we will see later, the dynamics of an $f\in {\rm Bih}_0(\CC)$ is
dominated essentially by its derivative $f'(0)$ at the fixed point 0,
usually called the {\it multiplier} of $f$ at 0. Therefore it would not be
surprising that the behavior of leaves near a given leaf ${\cal L}$ is
determined to a large extent by the multipliers at 0 of the monodromy mappings $f_\gamma $
for $\gamma \in \pi _1 ({\cal L})$. Our aim here is to give a formula for
$f'_\gamma (0)$ in terms of a path integration. We will use this formula
in the next section, where we will compute the multipliers of the monodromy mappings of the
leaf at infinity of an ${\cal F} \in {\cal A} _n'$. It must be mentioned that a
similar formula is proved by G. Reeb \cite{Re} for codimension 1 real foliations. Of
course that argument does not work here since our foliations have real codimension 2. However, the
fact that these foliations are induced by polynomial 1-forms allows us to find a short
interesting proof.

Let ${\cal F}:\{\omega =Pdy-Qdx=0 \} \in {\cal A}_n$, and suppose that ${\cal L}_{\infty}$, the line at infinity
with sing$({\cal F})$ deleted, is a leaf of ${\cal F}$ (see 1.29). Fix
some non-singular leaf ${\cal L} \neq {\cal L}_\infty $. It
follows that ${\cal L}$ is completely contained in the affine chart $(x,y)\in U_0 $.
Without loss of generality we may assume that $E:= {\cal L} \cap S_P$ is discrete in the leaf topology,
hence countable (see 1.8 for definition of $S_P$). In fact, if $E$ is not discrete, it has a limit point $z_0$ in
${\cal L}$. Parametrizing ${\cal L}$ around $z_0$ by $T\mapsto (x(T), y(T))$
with $z_0 =(x(0), y(0))$, we conclude that $P(x(T),y(T))\equiv 0$.
By analytic continuation, $P(x,y)=0$ for every $(x,y) \in {\cal L}$. Therefore we can
pursue the argument with ${\cal L} \cap S_Q$, which is finite since $P$ and $Q$ are
assumed to be relatively prime.

So assume that $E$ is discrete. By integrability of $\omega $, there exists a
meromorphic 1-form $\alpha =X dy -Y dx$, holomorphic on ${\cal L} \sm E$, such that
$d\omega =\omega \wedge \alpha$.\\ \\
{\bf 2.9 Theorem} \ {\it Given $\gamma \in \pi _1 ({\cal L} , p_0)$ one has}
\setcounter{equation}{0}
\begin{equation}
f'_\gamma (0) =\exp (- \int _\gamma \alpha ).
\end{equation}
{\bf Proof.} \ Consider ${\cal L}$ as the graph of a (multi-valued) function $y=\varphi (x)$
over some region of the $x$-axis. Let $E={\cal L} \cap S_P$ and $\tilde {E}:= \{x \in \CC :{\rm There\  exists\ } y\in \CC$
such that $(x,y) \in E \}$. Without loss of generality we can assume that the
base point $p_0 =(x_0 , y_0)$ is not in $E$. Moreover, we may replace $\gamma $ by
a path in its homotopy class that avoids $E$, if necessary. Take a vertical section $\Sigma $ parallel to
the $y$-axis, transversal to ${\cal L}$ at $p_0 $ , and let $y$ be the coordinate on $\Sigma $
(Fig. 7). For $y \in \Sigma $ near $y_0 $, each leaf ${\cal L}_y$ is the
graph of the solution $\Phi (x,y)$ of $dy/dx =Q(x,y) /P(x,y)$, i.e.,
$$\frac {\partial \Phi}{\partial x} (x,y)= \frac {Q (x, \Phi (x,y))}{P(x, \Phi (x,y))},\ \ \ \ \ \Phi (x_0 , y)=y.$$
Define $\xi (x):= \partial \Phi /\partial y (x, y_0).$ Note that\\
$$\begin{array}{ll}
{\displaystyle \frac {d \xi}{dx} (x)} & = {\displaystyle \frac{ \partial^2 \Phi}{\partial x \partial y} (x,y_0) } \\
 & \\
 & \displaystyle{=\frac {\partial}{\partial y} \left [ \frac {Q(x, \Phi (x,y))}{P(x, \Phi (x,y))} \right ] (x,y_0)}\\
 & \\
 & \displaystyle{= \left [ \frac {Q_y(x,\varphi (x)) P(x, \varphi (x)) -P_y (x, \varphi (x)) Q(x, \varphi (x))}{P^2 (x, \varphi (x))} \right ] \xi (x)}\\
 & \\
 & \displaystyle{=:T(x)\xi (x)},
\end {array}$$
\noindent
with $\xi (x_0)=1$. Thus $\xi (x)=\exp (\int ^x_{x_0} T(\tau ) d\tau)$, where the path of integration
avoids $\tilde {E}$. 


\realfig{figthesis6}{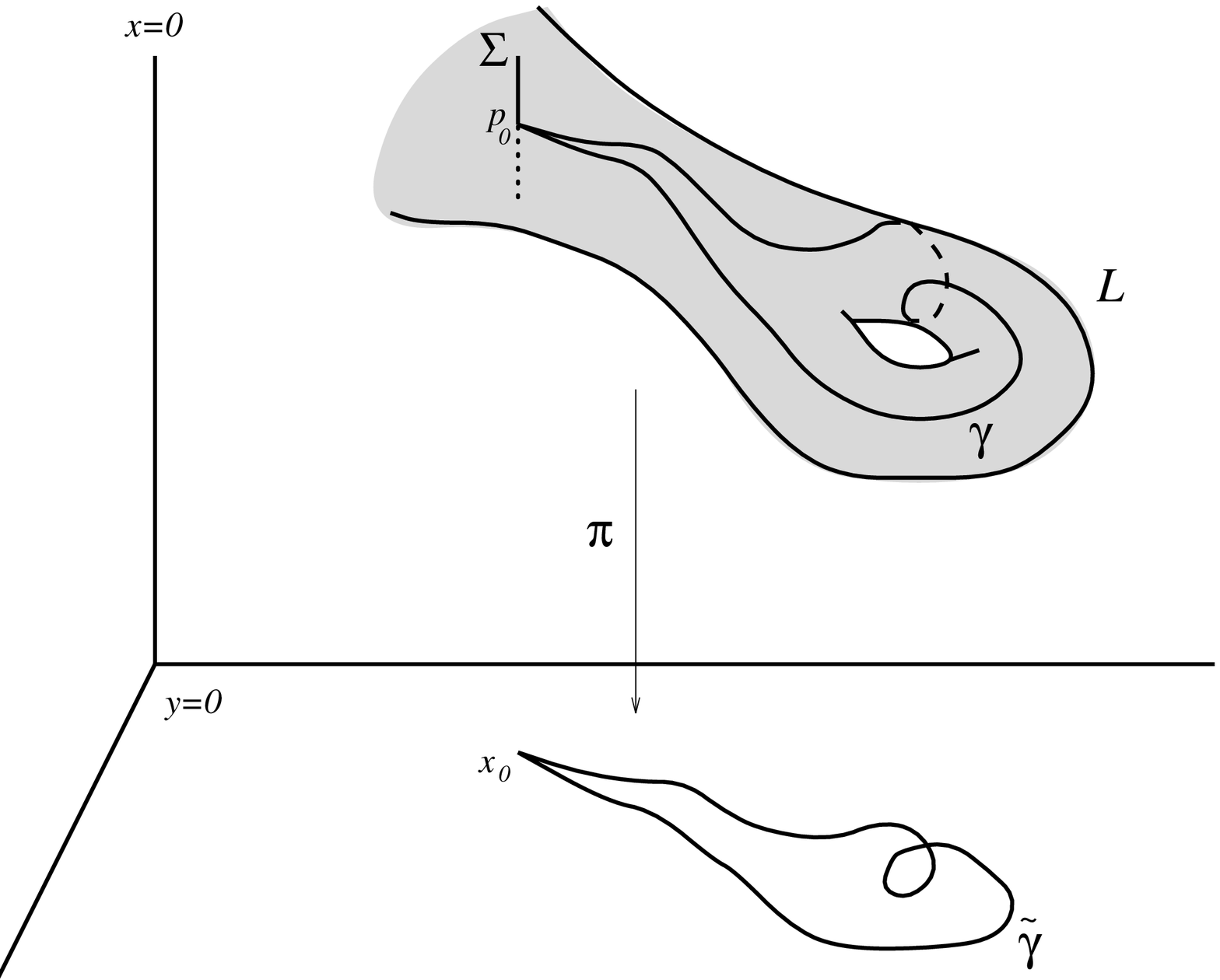}{}{9cm}


But $ f'_\gamma (0)$ is
the result of analytic continuation of $\xi $ along $\tilde {\gamma}$, the projection on
the $x$-axis of $\gamma $, so that $f'_\gamma (0)=\exp (\int _{\tilde {\gamma }}T(x)dx)$.
On the other hand, the condition $d\omega =\omega \wedge \alpha $ implies that
\begin{equation}
YP-XQ=P_x+Q_y,
\end{equation}
so that on an open neighborhood $W\subset {\cal L}$ of $\gamma $ one has
$$\begin{array}{ll}
\alpha |_W & = (Xdy-Ydx)|_W \vspace{4 mm}\\
& \displaystyle{ = \left ( \frac{XQ}{P}-Y \right ) \rule[-2.5mm]{.1mm}{8mm} _{\ W} \ dx} \vspace{4 mm}\\
& \displaystyle{ = -\left ( \frac {P_x+Q_y}{P}\right ) \rule[-2.5mm]{.1mm}{8mm}_{\ y=\varphi (x)} \ dx . \ \ \ \ \ \ \ \ \ \ {\rm (by (2))}}
\end{array}$$
Finally, note that
$$\begin{array}{ll}
\displaystyle{\int _{\tilde {\gamma }} T(x)dx+\int_\gamma \alpha} & \displaystyle{ = \int _{\tilde {\gamma }} \left [ T(x)-\left ( \frac {P_x+Q_y}{P} \right )\rule[-2.5mm]{.1mm}{8mm} _{\ y=\varphi (x)} \right ] dx} \vspace{5 mm}\\
& \displaystyle{= \int _{\tilde {\gamma }} \frac {Q_y P-P_yQ-P(P_x+Q_y)}{P^2}\rule[-2.5mm]{.1mm}{8mm} _{\ y=\varphi (x)} \ dx}\vspace{5 mm}  \\
& \displaystyle{= -\int _{\tilde {\gamma }} \frac {P_yQ+P_xP}{P^2} \rule[-2.5mm]{.1mm}{8mm} _ {\ y=\varphi (x)} \ dx } \vspace{5 mm} \\
& \displaystyle{= -\int _{\tilde {\gamma }} \frac {\partial /\partial x (P(x,\varphi (x)))}{P} \ dx} \vspace{5 mm}\\
& \displaystyle{= -\int _{\tilde {\gamma }} \frac {\partial } {\partial x } \log P(x, \varphi (x)) \ dx} \vspace{2 mm}\\
& \displaystyle{= 2\pi in}
\end{array}$$
\noindent
for some integer $n$ by the Argument Principle, proving the result. $\hfill \Box $\\ \\
\vspace{4 mm}
{\large \bf Monodromy Group of the Leaf at Infinity}\\
Let ${\cal F}\in {\cal A}_n'$ be a SHFC on $\CC \PP^2$ having Petrovski\u\i-Landis
property, i.e., $L_0\sm {\rm sing}({\cal F})$ is a leaf and
$L_0\cap {\rm sing}({\cal F})$ consists of exactly $n+1$ points $\{p_1,\ldots,
p_{n+1}\} $. Having an algebraic leaf homeomorphic to the punctured Riemann
sphere imposes some severe restrictions on the global behavior of the leaves.
Roughly speaking, for ${\cal F}\in {\cal A}_n'$ the global behavior of the
leaves is essentially determined by their local behavior in some neighborhood
of the leaf at infinity ${\cal L}_\infty $. The point is that in such a case
each non-singular leaf must accumulate on $L_0$ (Corollary 2.13 below).
Next, we can simply study the monodromy group of ${\cal L}_\infty$ to learn about the
behavior of nearby leaves.\\ \\
{\bf 2.10 Theorem} \ {\it Let $X=\sum_{j=1}^n f_j\ \partial/\partial z_j$ be a
holomorphic vector field on $\CC^n$. Then every non-singular solution of
the differential equation $dz/dT=X(z)$ is unbounded.}\\ \\
{\bf Proof}. Fix $z_0\in \CC^n$ with $X(z_0)\neq 0$, and suppose by
way of contradiction that the integral curve ${\cal L}_{z_0}$ passing through
$z_0$ is bounded. Let $T\mapsto \eta (T)$ be a local parametrization of ${\cal L}
_{z_0}$, with $\eta (0)=z_0$. Let $R >0$ be the largest radius such that $\eta (T)$
can be analytically continued over $\DD(0,R)$. If $R=+\infty $, then $\eta (T)$
will be constant by the Liouville's Theorem, contrary to the assumption $X(z_0)
\neq 0$. So $R<+\infty$. Recall from the Existence and Uniqueness Theorem
of solutions of holomorphic differential equations that for each $T_0\in
\CC$ and each $p\in \CC^n$, if $X$ is holomorphic on $\{ z\in \CC^n:
|z-p|<b\} $, then there exists a local parametrization $T\mapsto \eta _p(T)$ of
${\cal L}_p$, with $\eta _p(T_0)=p$, defined on $\DD(T_0, b/(M+kb))$, where $M={\rm sup}\{|X(z)|:|z-p|<b\}$ and $k={\rm sup}\{|dX(z)/dz|:|z-p|<b\}$ (see e.g. \cite{Co-L}). Now in our case, since ${\cal L}_{z_0}$ is bounded, $|X(\cdot)|$ and $|dX(\cdot)/dz|$
have both finite supremums on ${\cal L}_{z_0}$, and hence $b/(M+kb)$ may be chosen uniformly for
all $p\in {\cal L}_{z_0}$. This being so, for each $T_0 \in \DD (0,R)$ with
$R-|T_0| < b/(M+kb)$, there is a parametrization of ${\cal L}_{z_0}$ around $p=\eta (T_0)$,
so that $\eta (T)$ can be analytically continued over some larger disk, a contradiction. $\hfill \Box $\\ \\
{\bf 2.11 Remark} \ The same argument in the real case gives another proof of the
fact that an integral curve of a differential equation on ${\Bbb R}^n$ is either unbounded or it is
bounded and parametrized by the whole real line.\\ \\
{\bf 2.12 Corollary} \ {\it Let $ L\subset \CC \PP^2$ be any projective line and ${\cal L}$
be any non-singular leaf of a SHFC  \ ${\cal F}$. Then $\overline {{\cal L}} \cap L\neq \emptyset$}.\\ \\
{\bf Proof.} \ Choose an affine chart $(x,y)$ for $\CC \PP^2 \sm L \simeq \CC^2$.
In this coordinate system, $L=L_0$. Note that the leaves of ${\cal F}|_{\CC^2}$
are the integral curves of a (polynomial) vector field. Now the result follows from Theorem 2.10. $\hfill \Box $\\ \\
{\bf 2.13 Corollary} \ {\it Any non-singular leaf of a SHFC in ${\cal A}_n'$ has an accumulation
point on the line at infinity}. $\hfill \Box$\\

If the accumulation point is not singular, then the whole line at infinity is
contained in the  closure of the leaf by Proposition 1.5.

Knowing  that each leaf accumulates on $L_0$, we now proceed to study the
monodromy group of the leaf at infinity.\\
\setcounter{equation}{0}

\noindent
{\bf 2.14 Definition of $G_\infty$} \ Let ${\cal F}:\{ \omega =Pdy-Qdx=0\} \in {\cal A}_n'$
and $L_0\cap {\rm sing}({\cal F})=\{p_1,\ldots ,p_{n+1}\}$. Recall from 1.26
that ${\cal F}$ is described in the affine chart $(u,v)\in U_1$ by
\begin{equation}
\omega '(u,v)=u^{n+1}P(\frac{1}{u} , \frac{v}{u})dv-u^{n+1}R(\frac{1}{u} , \frac{v}{u})du=0,
\end{equation}
where $R(x,y)=yP(x,y)-xQ(x,y)$. The line at infinity $L_0$ is the closure of
$\{ (0,v):v\in \CC\}$, and $p_j:=(0,a_j), \ 1\leq j\leq n+1$, where $a_j$'s
are distinct roots in $v$ of the polynomial $u^{n+1}\displaystyle{R(\frac{1}{u} , \frac{v}{u})}|_{u=0}$.

The leaf at infinity ${\cal L}_\infty =L_0\sm \{ p_1,\ldots ,p_{n+1}\}$
can be made into a marked leaf by choosing fixed loops $\{ \gamma_1,\ldots ,
\gamma_n\}$ as generators of $\pi_1({\cal L}_\infty)$. Fixing a base point
$a\not \in \{p_1,\ldots,p_{n+1}\}$, each $\gamma_j$ goes around $p_j$ once in
the positive direction and does not encircle $p_j$ for $i\neq j$ (Fig. 8).

The monodromy mappings $f_{\gamma_j}$ for $1\leq j\leq n$ generate
the {\it monodromy group of the leaf at infinity}, denoted by $G_\infty$.\\
\setcounter{equation}{0}

\noindent
{\bf 2.15}  \ Let us compute the multiplier of each monodromy mapping $f_{\gamma_j}$
in terms of our data from $P$ and $Q$.

In the affine chart $(u,v)\in U_1$, the foliation is induced by the polynomial
1-form $\omega'$ of 2.14(1). Write
\begin{equation}
u^{n+1}P(\frac{1}{u} , \frac{v}{u})=:u\tilde{P}(u,v)\ \ \ \ \ {\rm and} \ \ \ \ \ u^{n+1}R(\frac{1}{u} , \frac{v}{u})=:\tilde{R}(u,v),
\end{equation}
where $\tilde{P}$ and $\tilde{R}$ are polynomials with $\tilde{R}(0,a_j)=0$ for
$1\leq j\leq n+1$. Note that ${\cal F}|_{U_1}$ may be viewed as the foliation
induced by the following vector field (cf. 1.9(1)):
$$X_1=u\tilde {P}(u,v)\frac{\partial}{\partial u}+\tilde {R}(u,v)\frac{\partial}{\partial v}.$$

\realfig{figthesis7}{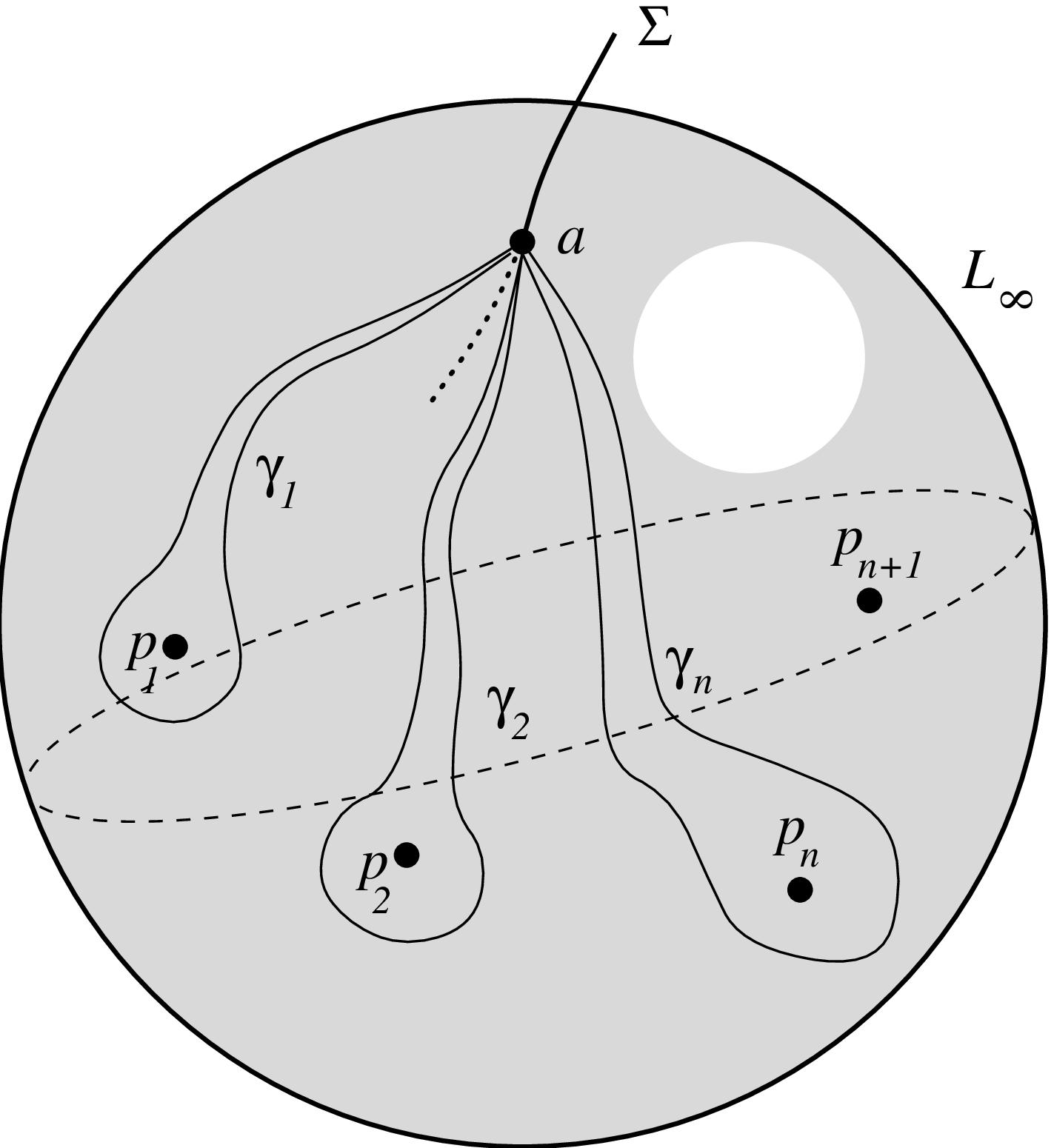}{}{6cm}

\noindent
Let us consider the Jacobian matrix $DX_1$ at the singular point $p_j$:
\begin{equation}
DX_1(p_j)=\left (
\begin{array}{cc}
\tilde {P}(0,a_j) & 0 \\
\tilde {R}_u(0,a_j) & \tilde {R}_v(0,a_j) \\
\end{array} \right )
\end{equation}
The quotient
\begin{equation}
\lambda _j:= \frac {\tilde {P}(0,a_j)}{\tilde {R}_v(0,a_j)}
\end{equation}
\noindent
of the eigenvalues of the matrix (2) is called the {\it characteristic number}
of the singularity $p_j$. Note that since the roots of $\tilde {R} (0,v)$ are
simple by assumption, we have $\tilde {R}_v (0,a_j)\neq 0$ and $\lambda _j$
of (3) is well-defined. On the other hand, the characteristic number is evidently
independent of the vector field representing ${\cal F}$ near $p_j$, since any two such
vector fields are the same up to multiplication by a nowhere vanishing holomorphic
function near $p_j$.

For simplicity we denote $f_{\gamma _j}$ by $f_j$ and $f_j '(0)$ by $\nu_j$.\\ \\
{\bf 2.16 Proposition} \ \ $\nu_j=e^{2\pi i\lambda _j}.$\\ \\
{\bf Proof.} \ By Theorem 2.9, we have $\nu _j=\exp (-\int _{\gamma_j} \alpha )$,
where $\alpha $ is any meromorphic 1-form which satisfies $d\omega '=\omega ' \wedge \alpha$. Choose, for example, $\alpha =-(u\tilde {P}_u+\tilde {P}+\tilde {R}_v)/\tilde {R} \ dv.$
Then we have
$$\begin{array}{ll}\vspace{4 mm}
\displaystyle{- \int _{\gamma _j} \alpha} & \displaystyle{=\int _{\gamma _j} \frac {u \tilde{P}_u+\tilde{P}+\tilde {R}_v}{\tilde {R}}\rule[-2.5mm]{.1mm}{8mm} _{\ u=0}\ dv}\\ \vspace{4 mm}
& \displaystyle{ =\int _{\gamma _j} \frac {\tilde {P}(0,v) +\tilde{R} _v(0,v)}{\tilde {R} (0,v)}\ dv }\\ \vspace{2 mm}
& \displaystyle{ =2\pi i \ {\rm Res} \left [ \frac {\tilde {P} (0,v)+\tilde {R} _v(0,v)}{\tilde {R} (0,v)} ; a_j \right ]}\\
& \displaystyle{ =2\pi i (\lambda _j+1),}
\end{array}$$
\noindent
so that $\exp (- \int _{\gamma _j} \alpha )=e^{2\pi i \lambda_j}$. $\hfill \Box$ \\ \\
{\bf 2.17 Remark} \ Since $f_1 \circ \cdots \circ f_{n+1} =$id, one has $\nu _1 \cdots \nu _{n+1} =1$
so that $\sum ^{n+1}_{j=1} \lambda _j$ is an integer by the above proposition.
However, this integer turns out to be 1 by the following argument. By 2.15(3), $\lambda _j$
is the residue at $a_j$ of the meromorphic function $\tilde {P} (0,v)/\tilde {R} (0,v)$
on $L_0 \simeq \overline {\CC}$. If $\tilde {R}(0,v)=c \prod _{j=1}^{n+1} (v-a_j),$
then $c$ is the coefficient of $y^{n+1}$ in $R(x,y)$ by 2.15(1), hence it is
the coefficient of $y^n$ in $P(x,y)$ since $R=yP-xQ$. So again by 2.15(1)
$\tilde {P}(0,v)$ is a polynomial in $v$ with leading term $cv^n$. It follows
that the residue at infinity of $\tilde {P}(0,v) /\tilde{R} (0,v)$ is $-1$. Hence
$\sum ^{n+1}_{j=1} \lambda _j-1=0$ by the Residue Theorem.\\ \\
\vspace{4 mm}
{\large \bf Equivalence of Foliations and Subgroups of Bih$_0(\CC)$ }\\
For singular smooth foliations by real curves on real manifolds there
are several notions of equivalence. One can think of topological or
$C^k$ {\it equivalences}, or topological or $C^k$ {\it conjugacies}.
In the case of equivalence, one is concerned only with the topology of the leaves, but in the case of conjugacy, the actual parametrization of the
leaves is also relevant.
Of course this is meaningful only when the foliation is described by a smooth
vector field on the ambient space.

In the complex analytic case one can still think of both equivalences and conjugacies
between SHFC's. But again the notion of conjugacy requires the existence of {\it
holomorphic vector fields on the whole space} representing our foliations.
Although one can develop a theory of conjugacies for holomorphic vector fields
defined on {\it open} subsets of $\CC^2$, it is rather awkward to deal with
conjugacies on {\it compact} complex manifolds, in particular $\CC \PP^2$
(cf. Proposition 1.11).

Consequently, the most natural notion of ``equivalence'' between SHFC's on $\CC \PP^2$
is furnished by the following\\ \\
{\bf 2.18 Definition} \ Two SHFC's ${\cal F}$ and ${\cal F'}$ on $\CC \PP^2$
are said to be {\it topologically} (resp. {\it holomorphically}) {\it equivalent}
if there exists a homeomorphism (resp. biholomorphism) $H:\CC \PP^2\rightarrow$
$\CC \PP^2$ which maps the leaves of ${\cal F}$ onto those of ${\cal F'}$.\\

The existence of equivalence between two SHFC's has the following implication
for the monodromy groups:\\ \\
{\bf 2.19 Proposition} \ {\it Let ${\cal F}$ and ${\cal F'}$ be topologically
(resp. holomorphically) equivalent SHFC's, linked by an equivalence
$H:\CC \PP^2\rightarrow \CC \PP^2$. Let p be a non-singular point for ${\cal F}$.
Then $G({\cal L}_p)$ is isomorphic to $G({\cal L}'_{H(p)})$. More precisely,
there exists a local homeomorphism (resp. biholomorphism) h on some neighborhood
of $0\in \CC$, with $h(0)=0$, such that for every $f\in G({\cal L}_p)$, $h\circ f=k(f)\circ h$,
where $k:G({\cal L}_p)\rightarrow G({\cal L'}_{H(p)})$ is a group isomorphism.}\\ \\
{\bf Proof.} Let $\gamma \in \pi_1({\cal L} _p,p)$ and $\Sigma $ be a section transversal
to ${\cal L}_p$ at $p. $ Set $q=H(p), \gamma ':=H_\ast \gamma \in \pi_1 ({\cal L} '_q, q)$, and
let $\Sigma '$ be a section transversal to ${\cal L} '_q$ at $q$. As in the definition of
the monodromy mapping, choose foliation charts $\{(U_i , \varphi _i )\}_{0\leq i \leq n} $
and transversals $\Sigma _i$ for ${\cal F}$, and the corresponding data $\{(U'_i, \varphi '_i )\}_{0\leq i \leq n} $
and $\Sigma '_i$ for ${\cal F'}$. Again we have a decomposition $f_\gamma =f_n \circ \cdots \circ f_0 $ and
$g_{\gamma '}=g_n\circ \cdots \circ g_0 $. Without loss of generality we may assume that
$U'_i=H(U_i)$. Since ${\cal F} |_{U_i}$ (resp. ${\cal F}'|_{U'_i})$ is a trivial foliation,
one has a projection along leaves $\pi _i:U_i\rightarrow \Sigma _i$ (resp. $\pi _i':U_i' \rightarrow \Sigma _i')$
which sends every $z$ to the unique intersection point of the plaque of $U_i$ (resp. $U_i'$) passing through $z$ with $\Sigma _i$ (resp. $\Sigma '_i)$.

Define $h_i:\Sigma _i\rightarrow \Sigma _i'$ by $h_i:=\pi '_i \circ H$. Since $H$
is a leaf-preserving homeomorphism (resp. biholomorphism) each $h_i$ is also a homeomorphism
(resp. biholomorphism) with inverse $h_i^{-1}=\pi_i \circ H^{-1}$ (see Fig. 9). Now
the definition of $f_i$ and $g_i$ shows that $h_{i+1} \circ f_i =g_i \circ  h_i$ for
$0\leq i \leq n$. Therefore two relations $f_\gamma =f_n\circ \cdots \circ f_0$
and $g_{\gamma'}=g_n\circ \cdots \circ g_0$ will show that $g_{\gamma'}\circ
h_0=h_0\circ f_\gamma$. To complete the proof, note that the mapping $f_\gamma
\mapsto k(f_\gamma):=h_0\circ f_\gamma \circ {h_0}^{-1}$ is an isomorphism between $G({\cal L}_p)$
and $G({\cal L}'_q)$. $\hfill \Box $\\

The above proposition suggests the following \\ \\
{\bf 2.20 Definition}  Two subgroups $G,G'\subset {\rm Bih}_0(\CC)$ are
said to be {\it topologically} (resp. {\it holomorphically}) {\it equivalent}
if there exists a homeomorphism (resp. biholomorphism) $h$ defined on some neighborhood of $0\in \CC$, with $h(0)=0$, such that $h\circ f\circ h^{-1}\in G'$ if and only if $f\in G$.


\realfig{figthesis8}{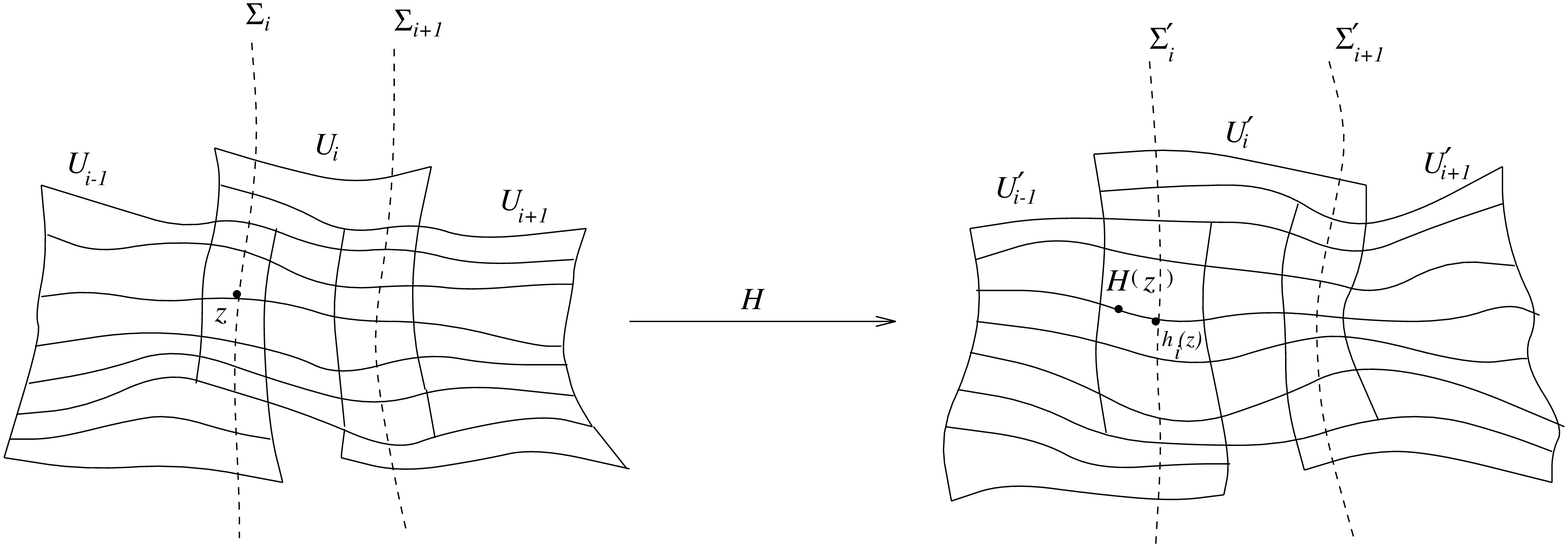}{}{14cm}


\setcounter{equation}{0}

It follows that the mapping $f\mapsto k(f):=h\circ f\circ h^{-1}$ is an isomorphism
$G \iso G'$, and the following diagram is commutative:
\begin{equation}
\begin{CD}
(\CC,0)@> f >> (\CC,0) \\
@VV{h}V      @VV{h}V   \\
(\CC,0)@> k(f) >> (\CC,0) \\
\end{CD}
\end{equation}
{\bf 2.21 Theorem } {\it If two SHFC's ${\cal F}$ and ${\cal F}'$ on $\CC \PP^2$
are topologically (resp. holomorphically) equivalent, then so are the monodromy
groups of the corresponding leaves}. $\hfill \Box $ \\

In view of this theorem, one naturally leads to the study of equivalence of subgroups
of ${\rm Bih}_0(\CC)$ in order to gain some informations about equivalence
of SHFC's. Below we give few important results concerning these equivalence
problems.\\ \\
{\bf 2.22 Theorem} \ {\it Let $G,G'\subset {\rm Bih}_0(\CC)$ be two topologically
equivalent marked subgroups linked by an orientation-preserving homeomorphism
h such that the diagram 2.20(1) is commutative. Suppose that G is non-Abelian,
and there exist $f_1,f_2\in G$ such that the multiplicative subgroup of $\CC ^{\ast}$ generated
by $f'_1(0), f'_2(0)$ is dense in $\CC$. Then h is actually a biholomorphism,
and G and $G'$ are holomorphically equivalent}. $\hfill \Box$ \\

This phenomenon is called {\it absolute rigidity} of subgroups of ${\rm Bih}_0(\CC)$.
Precisely, a subgroup $G\subset {\rm Bih}_0(\CC)$ is called {\it absolutely
rigid} if every subgroup which is topologically equivalent to $G$ is holomorphically equivalent to it. The proof of this theorem, which is based on
lifting the linearized germs to the universal covering of a punctured neighborhood of the
origin, as well as the approximation of linear germs by non-linear ones (see
Proposition 3.4), can be found in \cite{I4}.

Along this line, A. Shcherbakov \cite{Sh} has shown that \\ \\
{\bf 2.23 Theorem } \ {\it A non-solvable subgroup of\ {\rm Bih}$_0(\CC)$ is absolutely
rigid}. $\hfill \Box$ \\

For an almost complete topological and analytical classification of
germs in Bih$_0(\CC)$, see paper II of \cite{I5}.

The next result shows the existence of ``moduli of stability'' for SHFC's on
$\CC \PP^2$.\\ \\
{\bf 2.24 Theorem } \ {\it Let ${\cal F}, {\cal F}' \in {\cal A}'_n$ be two
SHFC's on $\CC \PP^2$ having no algebraic leaves other than the leaf at infinity.
Let $\{p_1, \ldots , p_{n+1} \}=L_0 \cap {\rm sing}({\cal F}),\ \{p'_1 , \ldots , p'_{n+1} \} =L_0 \cap {\rm sing} ({\cal F}')$,
and $\lambda _j$ and $\lambda '_j$ be the characteristic numbers of $p_j$ and $p_j'$, respectively.
Suppose that $\lambda _j$ and $\lambda _j'$ are non-zero, and ${\cal F}$ and
${\cal F}'$ are topologically equivalent by a homeomorphism $H:\CC \PP^2 \rightarrow \CC \PP^2$
with $H(p_j)=p'_j$. Then there exists an $
\Bbb R$-linear transformation
$A:\CC\rightarrow \CC$, with $A(\lambda _j)=\lambda '_j$ for
$1\leq j \leq n+1$}. $\hfill \Box $ \\

The proof of this result is based on investigating consequences of equivalence between
monodromy groups of the leaves at infinity, and uses the same techniques as the
proof of Theorem 2.22. It was proved for the first time by Yu. Il'yashenko \cite{I4} in the
case where $\lambda _j , \lambda _j'$ are not real numbers. Later, it was generalized by V. Naishul \cite{N}, who presented a much more difficult argument to handle the case where the characteristic numbers are non-zero and real.

Since it can be shown that the set of ${\cal F}\in {\cal A}_n'$ which do not
have any algebraic leaf other than ${\cal L}_\infty $ is open and dense in ${\cal A}_n$
(see 1.34 and proof of Proposition 3.22), it follows immediately that\\ \\
{\bf 2.25 Corollary} \ {\it No SHFC in ${\cal A}_n$ is structurally stable when
$n\geq 2 $}. $\hfill \Box$ \\

(Recall that ${\cal F}\in {\cal A}_n$ is structurally stable if there exists a neighborhood $\Omega \subset {\cal A}_n$ of ${\cal F}$ such
that every SHFC in $\Omega $ is topologically equivalent to ${\cal F}$.)

The following theorem, which proves a type of ``absolute rigidity'' for SHFC's, is a fundamental result first proved by Il'yashenko \cite{I4}.\\ \\
{\bf 2.26 Theorem }\ {\it A typical ${\cal F} \in {\cal A}_n$ is absolutely rigid. That
is, there exist neighborhoods $\Omega \subset {\cal A}_n$ of ${\cal F}$ and
U of the identity mapping on $\CC \PP^2$ in the uniform topology such that
every SHFC in $\Omega $ which is topologically equivalent to ${\cal F}$ by a
homeomorphism in U is holomorphically equivalent to ${\cal F}$}. $\hfill \Box $ \\

X. G\'omez-Mont \cite{G-O} has generalized the above theorem to SHFC's on projective
complex surfaces provided that there exists an algebraic leaf of sufficiently
rich homotopy group.
\newpage
\vspace*{4 mm}
\noindent
{\large {\bf Chapter 3}}\vspace{4 mm}\\
{\Large{\bf Density and Ergodicity Theorems}}\\ \\ \\ \\ \\ \\ \\ \\ \\ \\ \\ \\
\thispagestyle{plain}
\noindent
In the previous chapter we noted that the behavior of leaves
near the leaf at infinity ${\cal L}_\infty$ gives us information about
the global behavior of the leaves. The orbits of points under the action of the
monodromy group $G_\infty$ in turn give us a picture of the
behavior of leaves near ${\cal L}_\infty$. So a natural task is to consider
dynamics of germs in $G_\infty$, i.e., the iterations in a finitely-generated subgroup
of Bih$_0(\CC)$.

Here is a sketch of what will follow in this chapter. First we consider
elements of Bih$_0(\CC)$ without any attention to the relationship with
the monodromy groups of SHFC's. We study the linearization of hyperbolic germs
(Theorem 3.2), and approximation of a linear map by elements of a pseudo-group
of germs in Bih$_0(\CC)$ (Proposition 3.4) which leads us to a local density
theorem (Theorem 3.5). We then consider the notion of ergodicity in Bih$_0(\CC)$
and find conditions under which a finitely-generated subgroup of Bih$_0(\CC)$ is ergodic (Theorem 3.15).
Finally, these results will be applied to the monodromy group $G_\infty$ of a typical ${\cal F}\in {\cal A}_n$, leading to the density theorem of M. Khudai-Veronov (Theorem 3.25) and the ergodicity theorem of Yu. Il'yashenko and Ya. Sinai (Theorem 3.27).\\ \\
\vspace{4 mm}
{\large \bf Linearization of Elements in Bih$_0(\CC)$}\\
{\bf 3.1 Definition} A germ $f\in$ Bih$_0(\CC)$ is called {\it linearizable}
if there exists a holomorphic change of coordinate $\zeta =\zeta (z)$ near 0,
with $\zeta (0)=0$, such that
\setcounter{equation}{0}
\begin{equation}
\zeta (f(z))=f'(0) \cdot \zeta (z).
\end{equation}
In other words, we have the following commutative diagram:

$$\begin{CD}
(\CC,0) @> f >> (\CC,0) \\
@VV{\zeta}V  @VV{\zeta}V \\
(\CC,0) @> f'(0). >> (\CC,0) \\
\end{CD}$$
Therefore, a holomorphic change of coordinate {\it conjugates} a linearizable
germ with its tangent map at the fixed point 0.

If the germs $f$ and $g$ are conjugate, say if $\zeta \circ f \circ \zeta^{-1}=g$,
then $\zeta \circ f^n \circ \zeta^{-1}=g^n$ for every $n\geq 1$, so $f$ and $g$ have the same iterative dynamics near 0. Thus the possibility
of linearization can be very helpful in understanding the dynamics of iterations.

It is a remarkable fact that the possibility of linearization of a germ depends
crucially on the multiplier of the germ at the fixed point 0. In particular,
it was shown by G. Koenigs that the linearization of $f\in {\rm Bih}_0(\CC)$
is always possible if $f$ is {\it hyperbolic} in the sense that $|f'(0)|\neq 1$ \cite{M1}.\\ \\
{\bf 3.2 Theorem}\ {\it Every hyperbolic germ $f\in {\rm Bih}_0(\CC)$ is
linearizable. Moreover, the change of coordinate $\zeta$ is unique up to
multiplication by a non-zero constant.}\\ \\
{\bf Proof.} Let $\nu :=f'(0)$. Without losing generality one may assume that
$|\nu |<1$, for otherwise one can consider $f^{-1}$. Choose a constant $c<1$
so that $c^2<|\nu |<c$, and let $r>0$ be such that $|f(z)|\leq c|z|$ for
$z\in \DD(0,r)$. Thus for any $z_0\in \DD(0,r)$, the orbit $\{ z_n:=
f^n(z_0)\}_{n\geq 0}$ converges geometrically toward the origin, with $|z_n|
\leq rc^n$. But $|f(z)-\nu z|\leq k|z|^2$ for some constant $k$ and all $z\in
\DD(0,r)$, hence $|z_{n+1}-\nu z_n|\leq kr^2c^{2n}$. Set $\zeta_n(z):=f^n(z)/\nu^n$.
Then $|\zeta_{n+1}-\zeta_n|\leq (kr^2/|\nu |)(c^2/|\nu|)^n$. This estimate shows
that the sequence $\{ \zeta_n\}$ converges uniformly on $\DD(0,r)$ to a
holomorphic limit $\zeta$. The identity $\zeta \circ f=\nu \zeta$ is obvious.
Note that $\zeta'(0)=\lim_n\zeta'_n(0)=\lim_n\nu ^{-n}(f^n)'(0)=1$, so that
$\zeta$ is a local biholomorphism. If $\eta \circ f=\nu \eta$, then $\zeta \circ \eta^{-1}$
commutes with the linear map $z\mapsto \nu z$. Imposing this condition and comparing
the coefficients of the Taylor series expansions, it follows that $\zeta \circ \eta^{-1}=({\rm const.})z$,
and the required uniqueness follows. $\hfill \Box $\\ \\
{\bf 3.3 Remark} The problem of linearization of non-hyperbolic germs is extremely difficult and has a long history. In fact, if $\nu =e^{2\pi it}$ with $t\in (0,1)$ irrational, it turns out that the possibility of linearizing $f$ with $f'(0)=\nu$ depends on the asymptotic behavior of the denominators which appear in the rational approximations by the continued fraction expansion of $t$. Part of the linearization problem was solved by A. Brjuno in the mid $60$'s, but the complete solution has been achieved only in recent years with the work of J.C. Yoccoz and R. Perez-Marco \cite{P}.\\ \\ \\
\vspace{4 mm}
{\large \bf Approximation by Elements of a Pseudo-Group}\\
The next proposition shows how to approximate a linear map in a suitable
coordinate system by elements of a given pseudo-group $PG$ of germs in a finitely
generated subgroup $G\subset {\rm Bih}_0(\CC)$. The proof is nothing but
an elaboration of the following elementary fact: If $f\in {\rm Bih}_0(\CC)$
and $|\nu_1|<1$, then $f'(0)z=\lim_n \nu_1^{-n}f(\nu_1^nz)$.\\ \\
{\bf 3.4 Proposition} \ {\it Let G be a marked subgroup of ${\rm Bih}_0(\CC)$
with generators $ f_1,\ldots,f_k $, all defined on some domain $\Omega$
around $0$. Let $f_1$ be hyperbolic and $\zeta$ be a holomorphic coordinate change linearizing
$f_1$. Without loss of generality, assume that $|f'_1(0)|<1$, $\zeta$ is defined
on $\Omega$, and $\zeta(\Omega)=\DD(0,r)$ for some $r>0$. Let DG be the tangent group of G, i.e., the multiplicative subgroup of $\CC^\ast$
generated by $\nu_j:=f'_j(0)$, $1\leq j \leq k$. Then for every $\nu \in \overline{DG}$
there exists a sequence $F_n\in PG$ which converges to the linear map $\zeta \mapsto \nu \zeta$
uniformly on compact subsets of $\Omega \cap \nu^{-1} \Omega:=\{ z\in \Omega :
|\zeta(z)|<\min(r,r/|\nu|)\}.$}\\

\setcounter{equation}{0}
By an abuse of notation, we denote by $f(\zeta)$ the germ induced by $f$
in the coordinate $\zeta$, where $f\in {\rm Bih}_0(\CC)$ and $\zeta$ is
a holomorphic change of coordinate near 0.\\ \\
{\bf Proof.} \ It suffices to consider the case where $\nu \in DG$, that is $\nu =
f'(0)$ for some $f\in G$. The general case will then follow by
the uniformity of convergence and a standard diagonal argument. Define
\begin{equation}
F_n:=f_1^{-n}\circ f\circ f_1^n.
\end{equation}
First we claim that $F_n$ is defined on $\Omega$ for all sufficiently large $n$,
and $F_n(\zeta)\rightarrow \nu \zeta$ uniformly on compact subsets of $\DD(0,r)$ as $n\rightarrow \infty$.
In fact, if $f(\zeta)=\nu \zeta +\sum_{j=2}^\infty a_j \zeta^j$, then it is easily
seen from (1) that
\begin{equation}
F_n(\zeta)=\nu \zeta +\sum_{j=2}^\infty a_j \nu _1^{n(j-1)} \zeta^j.
\end{equation}
Since $f(\zeta)$ is holomorphic on $\DD(0,r')$ for some $0<r'<r$,
one has $\limsup_j \sqrt [j]{|a_j|}\leq 1/r'$, so for large $n,\ \limsup_j \sqrt [j]{|a_j||\nu_1|^{n(j-1)}}\leq 1/r $.
Thus the expression on the right hand side of (2) has an analytic continuation
over $\DD(0,r)$, i.e., $F_n$ is defined on $\Omega$.
Now the fact that $F_n(\zeta)\rightarrow \nu \zeta$ uniformly on compact subsets of $\DD(0,r)$ is an immediate consequence
of (2).

What remains to be shown is that for every compact set $K\subset \Omega \cap \nu^{-1} \Omega$
there is an $N=N(K)>0$ such that the domain $\Omega_{F_n}$ of $F_n$ as an element
of $PG$ (see 2.8) contains $K$ for all $n>N$.

Each intermediate representation of $F_n$ in (1) has the form
$$g_m:=f_1^m\ \ \ {\rm or}\ \ \ h_{mn}:=f_1^{-(n-m)}\circ f\circ f_1^n,\ 0\leq m\leq n.$$
We shall prove that for large $n$, $g_m(\zeta)$ and $h_{mn}(\zeta)$ have
conformal extensions over $\DD(0,r)$ and $K':=\zeta(K)$ is mapped into
$\DD(0,r)$ by them. In fact, $g_m(\zeta)=\nu_1^m \zeta$ and $h_{mn}(\zeta) \rightarrow \nu_1^m \nu \zeta$
uniformly in $\zeta$ and $m$ as $n\rightarrow \infty$, so that $g_m(\zeta)$
and $h_{mn}(\zeta)$ have conformal extensions over $\DD(0,r)$ for large
$n$. Now $g_m(\zeta (K))=\nu_1^m K'\subset \DD(0,r)$ for all $m\geq 0$.
Moreover, let $\delta :=\sup \{|\zeta (z)|: z\in K\} $, then $\delta <\min (r,r/|\nu|)$.
Choose $0<\epsilon <r-|\nu|\delta $, and find $N>0$ such that $|h_{mn}(\zeta)-\nu_1^m \nu \zeta|<\epsilon$
for all $0\leq m\leq n$ and $\zeta \in \DD(0,r)$ whenever $n>N$. Then if $\zeta \in K'$
we have $|h_{mn}(\zeta)|<\epsilon +|\nu_1|^m|\nu|\delta <r$ for $n>N$ and all
$0\leq m\leq n$. $\hfill \Box$\\ \\
{\bf 3.5 Theorem} \ {\it Let $G\subset {\rm Bih}_0(\CC)$ be a marked
subgroup, and suppose that the tangent group DG is dense in $\CC$. Then there exists an open neighborhood $\Omega$ of $0$ such that for every $z\in \Omega \sm \{0\}$ the orbit of $z$ under PG is dense in $\Omega$.}\\ \\
{\bf Proof.}  Since $\overline{DG}=\CC$, $G$ must contain at least one
hyperbolic germ, say $f_1$. Let $\Omega$ and $\zeta$ be as in the Proposition
3.4. By density of $DG$, to every $\nu \in \CC$ there corresponds a sequence
$\{F_n\}$ in $PG$ such that $F_n(\zeta)\rightarrow \nu \zeta$ uniformly on
compact subsets of $\Omega \cap \nu^{-1} \Omega$. Choose $z\in \Omega \sm \{0\}$
and let $w\in \Omega$ be arbitrary. Set $\nu :=\zeta(w)/\zeta(z)$. Then $F_n(\zeta (z))
\rightarrow (\zeta(w)/\zeta(z))\ \zeta(z)=\zeta(w)$ as $n\rightarrow \infty$, and
we are done. $\hfill \Box$\\ \\
\vspace{4 mm}
{\large \bf Ergodicity in Subgroups of Bih$_0(\CC)$}\\
\setcounter{equation}{0}
\noindent
Recall from ergodic theory that a measure-preserving transformation $T$ acting
on a probability space $X$ is called {\it ergodic} if every $T$-invariant
subset of $X$ has measure 0 or 1. In what follows, for two measurable sets $A,B \subset \CC$, the notation $A \doteq B$ means the symmetric difference $(A \sm B) \cup (B \sm A)$ has Lebesgue measure zero.\\ \\
{\bf 3.6 Definition} \ Two Lebesgue measurable subsets $A$ and $B$ of $\CC$
are said to be {\it equivalent at $0$} if there exists an open disk $U$ around
0 such that $A\cap U\doteq B\cap U$.
The {\it germ} of $A$ is the equivalence class of $A$ under
this relation, and is denoted by $[A]$. Given an $f\in {\rm Bih}_0(\CC)$,
a set $A$ is called $f$-$invariant$ if $[A]=[f(A)]$. Given a subgroup $G\subset
{\rm Bih}_0(\CC)$, a set $A$ is said to be $G$-$invariant$ if it is $f$-invariant
for every $f\in G$. The subgroup $G$ is called {\it ergodic} if for every $G$-invariant
set $A$, we have $[A]=[\CC]$ or $[\emptyset]$.\\ \\
{\bf 3.7 Example} \ Let $f(z)=\nu z,\ \nu \neq 0$, and $G$ be the subgroup of
Bih$_0(\CC)$ generated by $f$. Then $G$ is not ergodic. If $|\nu |\neq 1$,
this is obvious (consider a snail-like region spiraling toward the origin; compare Fig. 11).
If $|\nu |=1$, consider the set $A=\{re^{it}:t\in {\Bbb R} \ {\rm and} \ 1/(2n+1)\leq r\leq 1/(2n)\
{\rm for\ some}\ n\geq 1 \}$. However, we will see that if $G$ is generated by $f_1(z)=\nu_1 z$ and
$f_2(z)=\nu_2 z$, with the tangent group $DG=\langle \nu_1 , \nu_2 \rangle$ being dense in $\CC$, then $G$ is
ergodic (see Proposition 3.13).\\ \\
{\bf 3.8 Density Points } \ Let $A\subset \CC$ be Lebesgue measurable.
A point $z\in \CC$ is called a {\it density point} of $A$ if
\begin{equation}
\lim_{r \rightarrow 0} \frac {m (A \cap \DD(z,r))} {m (\DD(z,r))}=1,
\end{equation}
where $m$ denotes the Lebesgue measure on $\CC$.

It is a consequence of the theorem of existence of ``Lebesgue points'' for an
$f\in L^1(m)$ (see for example \cite{Ru}) that the limit on the left hand side of
(1) is just $\chi_A$, the characteristic function of $A$, up to a set of measure
zero. As a consequence, {\it almost every point of $A$ is a density point of $A$}.

Our main goal is to prove that under certain weak assumptions, a marked subgroup of Bih$_0(\CC)$ is ergodic. The proof is based on an idea due to E. Ghys.\\ \\
{\bf 3.9 Theorem}\ {\it Let f be a univalent function on $\DD(z,r)$. Then $f(\DD(z,r))$ contains $\DD(f(z),|f'(z)|r/4)$.}\\ \\
{\bf Proof.} This is the so-called ``Koebe 1/4--Theorem'' \cite{Ru}. $\hfill \Box$\\ \\
{\bf 3.10 Lemma}\ {\it Let $f(z)=\nu z$, with $0<|\nu| <1$, and $g\in {\rm Bih}_0(\CC)$.
Suppose that A is both f- and g-invariant. Then
there exists an open disk U around $0$ such that for all
$n\geq 0$, $A\cap U\doteq \nu^{-n}g\nu^n(A)\cap U$.}\\ \\
{\bf Proof.}\ For each integer $n$, $A$ is clearly invariant under $f^{-n}\circ g\circ f^n$.
To prove the assertion, it suffices to find a disk $U$ that works uniformly for
all sufficiently large $n$. Below we replace $\doteq $ by $=$ through
modification by a set of measure zero if necessary.

Choose a disk $V=\DD(0,r_1)$ such that
\setcounter{equation}{0}
\begin{equation}
A\cap V=(\nu A)\cap V=(\nu^{-1}A)\cap V=g(A)\cap V.
\end{equation}
Set $\mu :=1/g'(0)$. Since $\nu^{-n}g^{-1}\nu^n$ converges uniformly to $z\mapsto \mu z$,
there is a disk $U=\DD(0,r_2)\subset V$ and an $N>0$ such that $|\nu^{-n}g^{-1} (\nu^n z)|\leq
|\mu||z|+r_1/2$ for all $z\in U$ and all $n>N$. Moreover, we may choose
$r_2$ such that $|\mu|r_2<r_1/2$ and $N$ large enough such that $\nu^nU
\subset g(U)$ for all $n>N$.

Now suppose that $z\in A\cap U$. Then by (1), $w=\nu^n z \in A\cap U$.
Again by (1) we have $w=g(x)$, where $x\in A$. If $n>N$, $x$ is in fact in $A\cap U$
by the choice of $N$. If $y=\nu^{-n}x$, then $|y|=|\nu^{-n}g^{-1}(\nu^n z)|\leq |\mu||z|+r_1/2
\leq |\mu|r_2+r_1/2<r_1$, so that $y\in V$. Once again by (1) we obtain $y\in
A\cap V$. Hence, $z\in \nu^{-n}g\nu^n(A)\cap U$.

Similarly, one can show that $\nu^{-n}g\nu^n(A)\cap U\subset A\cap U$ for all
large $n$, and this completes the proof. $\hfill \Box$\\ \\
{\bf 3.11 Proposition}\ {\it Let $f_n\in {\rm Bih}_0(\CC)$ be defined on
some open neighborhood V of $0$ for all $n\geq 1$, and $f_n\rightarrow f\in {\rm Bih}_0(\CC)$
uniformly on compact subsets of V. Suppose that A is $f_n$-invariant for $n\geq 1$.
If there exists a disk $U\subset V$ around $0$ such that $A\cap U\doteq f_n(A)\cap U$
for all $n\geq 1$, then A is f-invariant.}\\ \\
{\bf Proof.}\ Without loss of generality we may assume that each point of $A$
is a density point. Choose a smaller disk $W\subset U\cap f(U)\cap f^{-1}(U)$ such that $f_n(W)
\subset U$ for all $n$. Choose $z_0\in A$ such that $f(z_0)\in W$. We show that
$f(z_0)$ is a density point for $A$. This proves that $f(A)\cap W\subset A
\cap W$. Next by the choice of $W$, $f_n^{-1}(A)\cap W=A\cap W$, so that the
same argument leads us to $A\cap W\subset f(A)\cap W$.

So let $z_0\in A$ with $f(z_0)\in W$. By the choice of $W$, we have $z_0\in U$.
Given a small $\epsilon >0$, there exists an $r=r(\epsilon )>0$ such that
$\DD(z_0,r)\subset U$, and for every $0<r'<r$,
\setcounter{equation}{0}
\begin{equation}
1-\epsilon <\frac{m(A\cap \DD(z_0,r'))}{m(\DD(z_0,r'))}\leq 1.
\end{equation}
Choose an arbitrary $\mu$ with $|f'(z_0)|-2\epsilon <\mu <|f'(z_0)|-\epsilon $.
By Theorem 3.9, $f_n(\DD(z_0,r))$ contains $\DD(f_n(z_0),\mu r/4)$ for
$n$ large enough (Fig. 10).

For simplicity, set $D:=\DD(z_0,r),\ D':=\DD(f(z_0),\mu r/4)$, and
$D_n':=\DD(f_n(z_0),\mu r/4)$. Since $D'_n\sm A\subset f_n(D)\sm
(A\cap f_n(D))=f_n(D)\sm (f_n(A\cap D))$, for $n$ large enough we have
$$\begin{array}{ll} \vspace{4 mm}
m(A\cap D'_n) & \geq m(D'_n)-m(f_n(D))+m(f_n(A\cap D))\\ \vspace{4 mm}
 & = \displaystyle {m(D'_n)-\int _D |f_n'|^2 dm + \int _{A\cap D} |f'_n|^2 dm }\\ \vspace{4 mm}
 & \geq m(D'_n)-(|f'(z_0)| + \epsilon )^2 m(D)+\mu ^2 m(A\cap D)\\ \vspace{4 mm}
 & \geq [1-16 (|f'(z_0)|+\epsilon )^2 \mu ^{-2} + 16 (1-\epsilon)] m(D'_n) \ \ \ \ \ \ \ \ \ {\rm (by (1))}\\ \vspace{4 mm}
 & =: \ell (\epsilon )m(D'_n).
\end{array}$$
As $n\rightarrow \infty$, $m(A\cap D'_n)\rightarrow m(A\cap D')$
by the Lebesgue's Dominated Convergence Theorem. Therefore
$$\frac {m(A\cap D')}{m(D')} \geq \ell (\epsilon ).$$
As $\epsilon \rightarrow 0, \ \ell (\epsilon ) \rightarrow 1$, and it follows that
$f(z_0)$ is a density point of $A$. $\hfill \Box$\\ \\ 
\realfig{figthesis9}{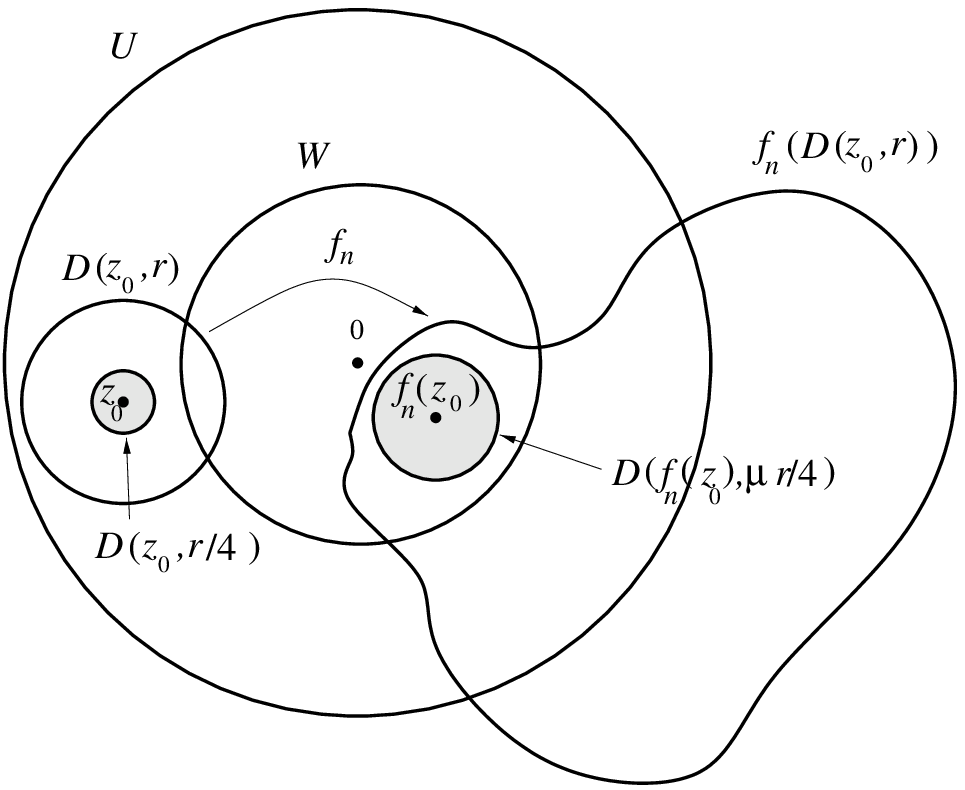}{}{7cm}
{\bf 3.12 Corollary} {\it Under the assumptions of Proposition 3.4, every G-invariant set A
is also DG-invariant.}\\ \\
{\bf Proof.} Let $\nu \in DG$ be equal to $g'(0)$ for some $g\in G$. In the
coordinate $\zeta $ one has $f_1(\zeta )=\nu_1 \zeta $, and by Proposition 3.4, $f_1^{-n} \circ g \circ f_1^n$
converges to $\zeta \rightarrow \nu \zeta $ uniformly on compact subsets of
$\Omega \cap \nu ^{-1} \Omega $. It follows then from Lemma 3.10 and
Proposition 3.11 that $A$ is invariant under $\zeta \mapsto \nu \zeta.\ \hfill \Box $ \\

The proof of the following proposition uses a standard technique in
ergodic theory (see Appendix 11 of \cite{A-A}).\\ \\
{\bf 3.13 Proposition} {\it Let $\nu_j\in \CC^\ast , 1 \leq j \leq k$, and G be
the subgroup of {\rm Bih}$_0(\CC)$ generated by $z \mapsto \nu _j z,\ 1 \leq j\leq k$. Then G is ergodic if and only if DG is dense in $\CC$.} \\ \\
{\bf Proof.} Suppose that $DG$ is dense in $\CC$. It follows that $G$ contains at
least one hyperbolic element $f(z)=\nu z$ with $|\nu |<1$. Let $A \subset \CC$ be any
$G$-invariant measurable set. Choose a disk $U$ around $0$ such that $A\cap U\doteq (\nu A)\cap U$. Take
the quotient of $U \sm \{0 \}$ under the action of the group
$\{f^n:n\in {\ZZ} \}$. This quotient is biholomorphic to a 2-torus $\TT^2$ (see Fig. 11).
Let $\tilde {G}$ be the induced group of translations of $\TT^2$, and $\tilde {A}$ be the induced
measurable subset of $\TT^2$. Note that $\tilde {A}$ is invariant under the action of $\tilde {G}$, and the orbit of each point in $\TT^2$
is dense under this action.

It suffices to show that $\tilde {A}$ or $\TT ^2 \sm \tilde {A}$ has measure zero. Expand the characteristic function of $\tilde {A}$ into the Fourier series
$$\chi _{\tilde {A}}(e^{2\pi ix},e^{2\pi iy}) = \sum _{m,n} a_{mn} e ^{2\pi i (mx+ny)} ,$$
where we identify $\TT^2$ with $\{ (e^{2 \pi i x}, e ^{2 \pi i y})\in S^1 \times S^1 \}$. Let
$\tilde {f} \in \tilde {G}$ be the translation $(e^{2 \pi i x}, e ^{2 \pi i y})\mapsto (e^{2 \pi i (x+\alpha ) }, e ^{2 \pi i (y+ \beta ) })$.
The $\tilde {G}$-invariance of $\tilde {A}$
shows that
$$\chi _{\tilde {A}} = \sum _{m,n} a_{mn} e^ {2 \pi i (m \alpha +n \beta )} e ^{2\pi i (mx+ny)}.$$
Therefore, for all $m,n \in {\ZZ},\ a_{mn} =a_{mn} e^{2 \pi i (m \alpha +n \beta )}$.
Since $\tilde {G}$ contains at least one irrational translation (otherwise, the orbit of each
point would be finite), we conclude that $a_{mn}=0$ for all $(m,n) \neq (0,0)$, therefore
$\chi _{\tilde {A}} =0$ or 1 almost everywhere.

Conversely, let $G$ be ergodic. Clearly $G$ contains a hyperbolic element
$f_0(z)=\nu z$, with $|\nu |<1$. Suppose by way of contradiction that $DG$ is not dense in $\CC$,
i.e. there exists an open disk $\DD(z,r)$ such that $\overline {DG} \cap \DD(z,r)=\emptyset $.
Set $A :=\bigcup_{f \in G} f(\DD(z,r/2))$. Then $A$ is an open, $G$-invariant
set such that $\overline {DG} \cap A =\emptyset $. The germ of $A$ at $0$, [$A$], is not
equal to [$\emptyset $] since every neighborhood of 0 contains $\nu ^n \DD(z,r/2)$ for
$n$ large enough. It follows that [$A$]=[$\CC$]. Invariance of $A$ and the
fact that $|\nu |<1$ will then show that $m(\CC \sm A)=0$. In particular,
$A$ is dense in $\CC$. Let $\{ \zeta _n\}$ be a sequence of elements of $A$
such that $\lim _n \zeta _n $=1. Let $\zeta _n =\sigma _n z_n$, where $\sigma _n \in DG$ and
$z_n \in \DD(z,r/2)$. Then it follows that the sequence $\{ \sigma _n^{-1}\} $
enters $\DD(z,r)$ for $n$ large enough. This contradicts the assumption
$\overline {DG} \cap \DD(z,r)=\emptyset $.\ $\hfill \Box $\\ \\
{\bf 3.14 Remark } The above proof shows how the notion of ergodicity for finitely-generated
subgroups of Bih$_0(\CC)$ containing a hyperbolic germ is related to the
usual notion of ergodicity for translations of tori, thus justifying Definition 3.6.

\realfig{figthesis10}{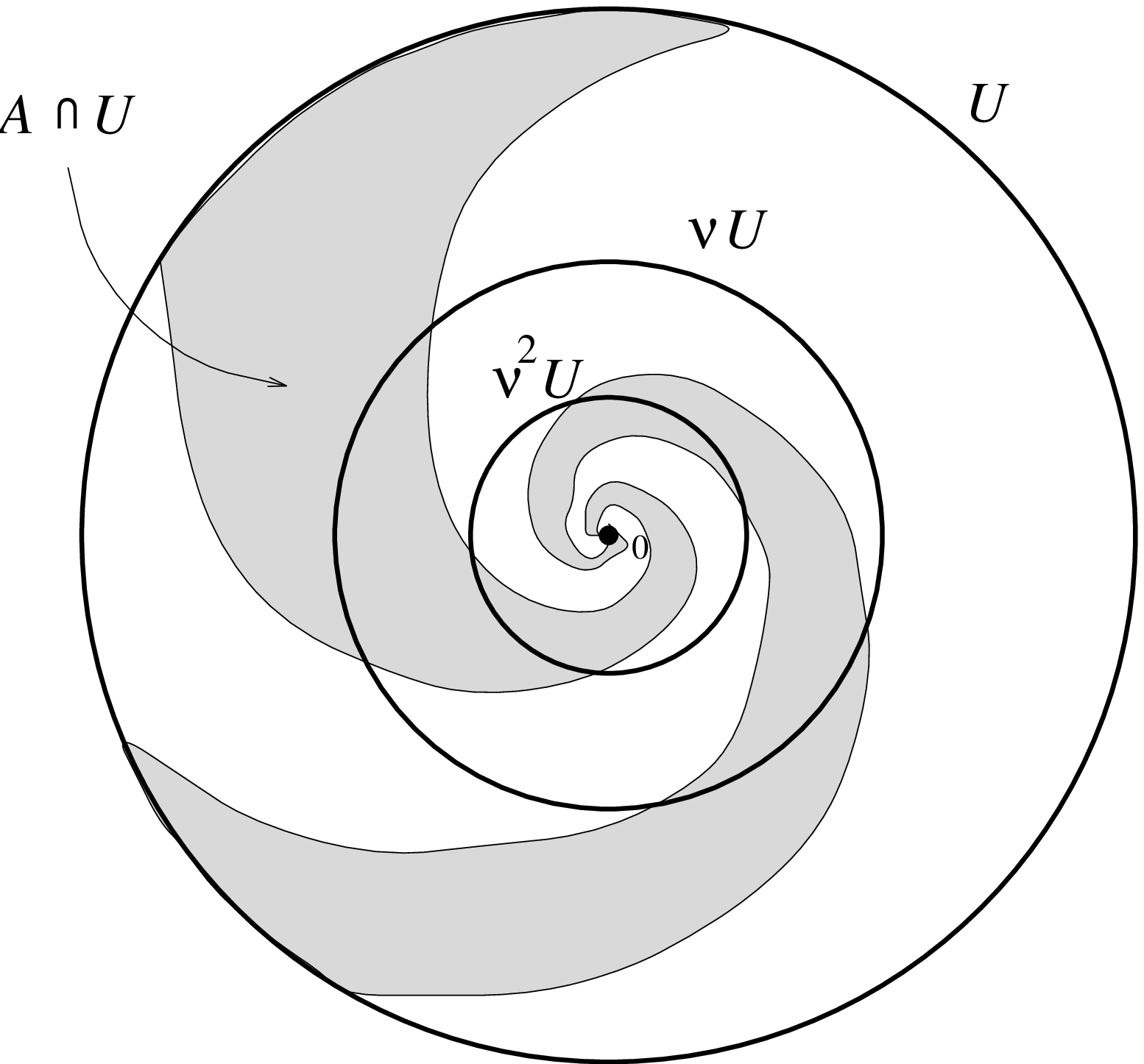}{}{6cm}

\noindent
{\bf 3.15 Theorem} {\it Let $G\subset {\rm Bih}_0 (\CC)$ be a marked subgroup and $DG$ be dense in $\CC$. Then $G$ is ergodic.} \\ \\
{\bf Proof.} Since $DG$ is dense in $\CC$, $G$ must contain at least one hyperbolic
element. By Corollary 3.12, every $G$-invariant set $A$ is also $DG$-invariant. By
Proposition 3.13, $[A]=[\CC]$ or $[\emptyset ]$. $\hfill \Box$\\ \\
\vspace{4 mm}
{\large \bf Density of Leaves of SHFC's on $\CC \PP^2$}\\
This section is devoted to a proof of the Khudai-Veronov's theorem on density of leaves of a typical SHFC on $\CC \PP^2$. We first prove a version of this
theorem which asserts that for a typical ${\cal F} \in {\cal A}_n$ all but a finite number of leaves are dense in $\CC \PP^2$. Next, by
applying a more elaborate argument, we show that among these exceptional leaves, only ${\cal L}_\infty $ is robust and in fact for a typical ${\cal F} \in {\cal A}_n$ {\it all}\ leaves except ${\cal L}_\infty$ are dense in $\CC \PP^2$.\\ \\
{\bf 3.16 Singularities of Hyperbolic Type } \ Let $X$ be a holomorphic
vector field defined on some domain $U$ containing $p\in \CC^2$, and let $p$
be an isolated singular point of $X$. Let $\sigma _1 \ {\rm and} \ \sigma _2$ be the eigenvalues of the
Jacobian matrix $DX(p)$. We say that $p$ is a {\it non-degenerate} singularity if
$\sigma _1 \sigma _2 \neq 0$; it is called a {\it hyperbolic} singularity if
it is non-degenerate and $\sigma _1/ \sigma _2 \not \in {\Bbb R}$. \\ \\
{\bf 3.17 Theorem } \ {\it Every hyperbolic singularity can be holomorphically linearized, i.e., there exist neighborhoods $V\subset U$ of p and W of $0\in \CC^2$ and a biholomorphism $\varphi: V\rightarrow W$ such that $(\varphi _\ast X)(x,y)=\sigma_1 x\ \partial /\partial x+ \sigma _2 y\ \partial /\partial y$}. $\hfill \Box $\\

In fact, in this two dimensional case, the linearization is possible when
either $\sigma _1/\sigma _2 \not \in {\Bbb R}$, or $\sigma _1/\sigma _2$ is positive
but not an integer or the inverse of an integer (the so-called ``non-resonant Poincar\'e case'' ).
For a proof, see \cite{Ar}.

By the above theorem, the local picture of a hyperbolic singularity can be
totally explored since one can easily integrate the linearized vector
field.

Recall that for a holomorphic vector field $X$ with an isolated singularity
at $p$, a {\it local separatrix } through $p$ is the image of a punctured disk
$\DD^\ast (0,r)$ under a holomorphic immersion $\eta $ such that $d\eta(T)/dT =X(\eta (T))$
for $T\in \DD^\ast (0,r)$, and $\lim _{T\rightarrow 0} \eta (T) =p$. In
other words, a local separatrix of $p$ is an integral curve of $X$ which
passes ``nicely'' through $p$. It is shown by C. Camacho and P. Sad that
every holomorphic vector field $X$ with an isolated singularity at $0\in \CC^2$
has a local separatrix through 0 \cite{C-S}. On the other hand, the same conclusion is false in higher dimensions, even in $\CC^3$ \cite{G-L}.\\ \\
{\bf 3.18 Corollary} \ {\it Let X be a holomorphic vector field defined near
an isolated hyperbolic singularity $p\in \CC^2$. Then X has exactly two local
separatrices through p, and every other integral curve passing
through the linearizing neighborhood accumulates on these two separatrices.}\\ \\
Fig. 12 is an attempt to visualize this situation.\\ \\
\realfig{figthesis11}{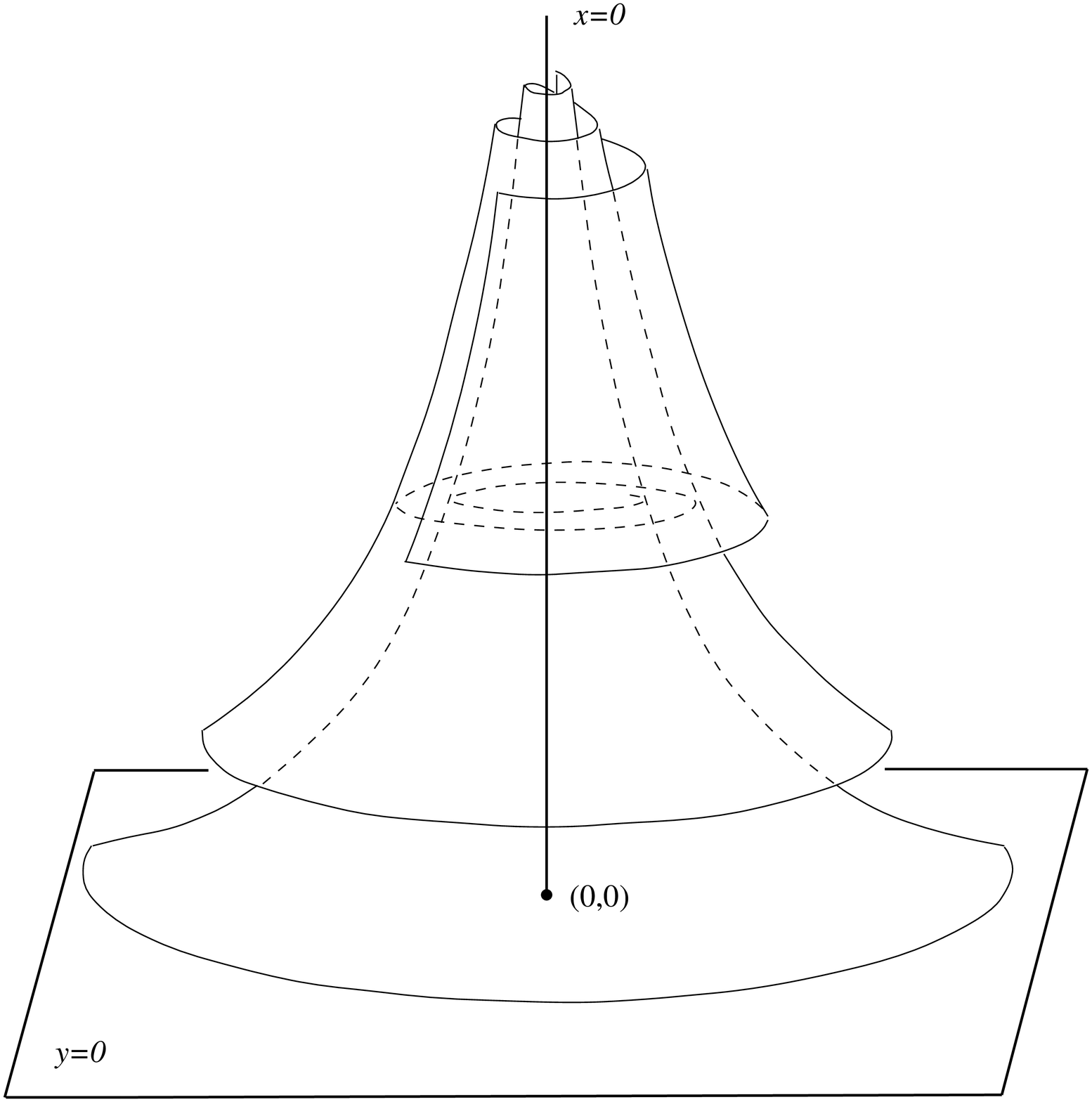}{\sl A hyperbolic singularity in $\CC^2$}{8cm}
{\bf Proof.}  Locally $X$ can be transformed into $\sigma _1 x\ \partial /\partial x + \sigma _2 y\ \partial /\partial y $
by Theorem 3.17. The integral curve passing through $(x_0,y_0)$ is parametrized by $T\mapsto (x _0e^{\sigma _1 T} , y_0e^{\sigma _2T})$. It follows that the punctured axes $\{x=0\} \sm \{ (0,0)\}$ and $\{y=0\} \sm \{ (0,0)\}$
are local separatrices. Since $\sigma _1/\sigma _2 \not \in {\Bbb R}$,
there exist sequences $\{ T_n \}$ and $\{ T'_n \}$ such that
$$\begin{array}{ll}
e^{\sigma _1 T_n} &=1 \ (n=1, 2, \ldots ) , \ \ \ \ \ e^{\sigma _2 T_n}\rightarrow 0 \ ({\rm as} \ n \rightarrow \infty ),\\
e^{\sigma _2 T'_n} &=1 \ (n=1, 2, \ldots ) , \ \ \ \ \ e^{\sigma _1 T'_n}\rightarrow 0 \ ({\rm as} \ n \rightarrow \infty ).\\
\end{array}$$
It follows that if $x_0 y_0 \neq 0$, the integral curve passing through $(x_0,y_0)$ accumulates on
$(x_0,0)$ and $(0,y_0)$, hence on the whole $\{x=0\}$ and $\{ y=0\}$ by Proposition 1.5. $\hfill \Box $\\ \\
{\bf 3.19 Relation with $G_\infty $} \ Let ${\cal F} \in {\cal A}_n',\ L_0 \cap$ sing$({\cal F}) = \{p_1 , \ldots , p_{n+1} \}$, and
$\lambda _j$ be the characteristic number of $p_j$, as in 2.15(3). It follows that $p_j$ is a hyperbolic
singularity of $X_1$ (hence of any vector field representing ${\cal F}$ near $p_j$) if and only if $\lambda _j\not \in {\Bbb R}$.
By Proposition 2.16, the last condition is equivalent to $|\nu _j|\neq 1$, where $\nu _j$
is the multiplier at 0 of the monodromy mapping $f_j \in G_\infty $ (see 2.14). We conclude that {\it $p_j$ is a hyperbolic singularity if and only if $f_j$ is a hyperbolic germ in} Bih$_0(\CC)$.

Note that if $p_j$ is hyperbolic, then one of the local separatrices through $p_j$ is
contained in the leaf at infinity itself; the other one is transversal to ${\cal L}_\infty$ (Fig. 13).\\ \\
\realfig{figthesis12}{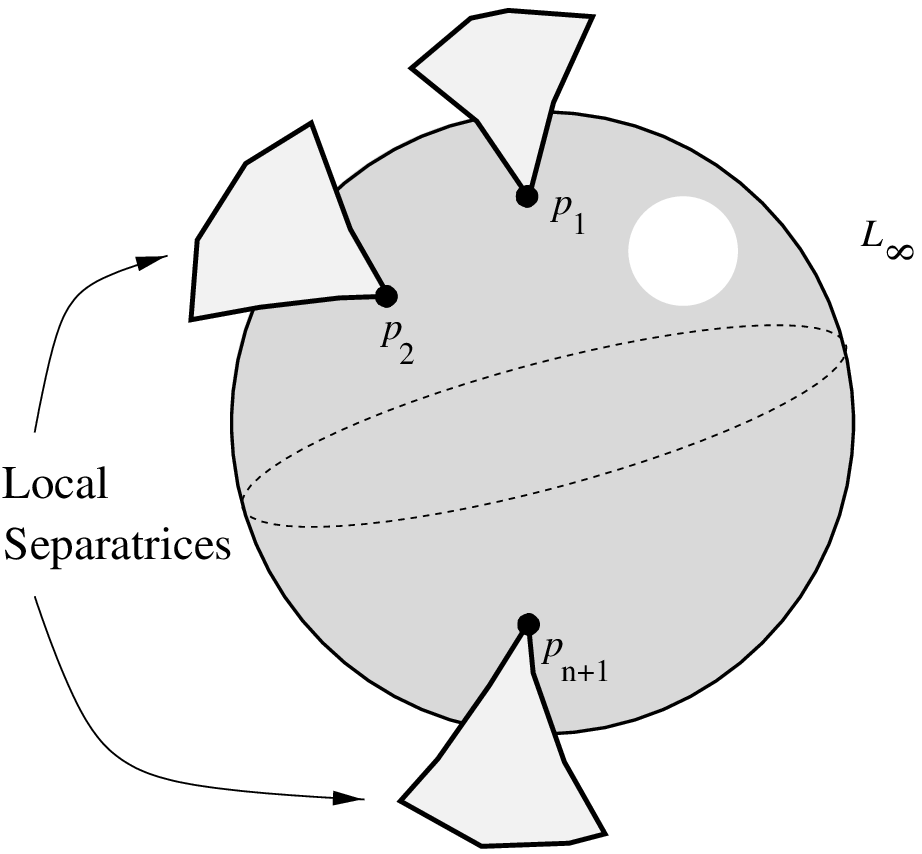}{}{7cm}
{\bf 3.20 Proposition} \ {\it A typical ${\cal F} \in {\cal A}_n,\ n\geq 2$, satisfies the following conditions:}\\
(i) \ \ ${\cal F}\in {\cal A}'_n,$\\
(ii) \ {\it $|\nu _j |\neq 1\ (1 \leq j\leq n)$; in other words, all generators of $G_\infty $ are hyperbolic,}\\
(iii) {\it $DG_\infty = \langle \nu_1 , \ldots , \nu _n \rangle$ is dense in $\CC$}.\\ \\
{\bf Proof.} Since the union of sets of measure zero has measure zero, it
suffices to prove that each condition is typical in ${\cal A}_n$.

The first condition is typical by 1.34. The second one is clearly typical by
2.15(3): If $|\nu _j|=1$, then $\lambda_j\in {\Bbb R}$, and this can be easily destroyed by perturbing the coefficients of $P$ and $Q$.

The third condition is more subtle. By 2.15(3) and 2.16, it suffices to prove that
for almost every $(\lambda_1, \ldots , \lambda _n ) \in \CC^n$, the
multiplicative subgroup generated by $\{ e^{2\pi i \lambda _1}, \ldots , e^{2\pi i\lambda _n } \}$ is dense in $\CC$. Evidently it is enough to prove this statement for $n=2$, since then we can
take the product of the resulting subset of $\CC^2$ by $\CC^{n-2}$ to obtain a
subset of full measure in $\CC^ n$.

We shall prove that for almost every $(\lambda _1, \lambda _2)\in \CC^2$,
the additive subgroup generated by $\{1,\lambda _1 , \lambda _2 \} $ is dense in $\CC$. Suppose that $(\lambda _1 , \lambda _2 )$ is chosen such that
\begin{enumerate}
\item[$(\ast )$]
No two vectors in $\{1, \lambda _1, \lambda _2 \}$ are $\RR$-dependent,
\item[$(\ast \ast )$]
If $1=a\lambda_1 +b\lambda _2$, then $b/a$ is irrational.
\end{enumerate}
Let $\Lambda $ be the lattice generated by $\{\lambda _1, \lambda _2 \}$, and consider the quotient torus
${\RR}^2 \stackrel {\pi } {\longrightarrow }{\Bbb R}^2/ \Lambda $.
Let $L:{\Bbb R}^2 \rightarrow {\Bbb R} ^2$ be a linear map with $L(\lambda _1 )=(1,0), \ L(\lambda _2 )=(0,1)$. Then
$L(1) = (a,b)$, and $L$ induces a homeomorphism $\tilde {L}:{\Bbb R}^2 /\Lambda \rightarrow {\Bbb R} ^2 /{\ZZ}^2$ such that the following diagram commutes:
\setcounter{equation}{0}
\begin{equation}
\begin{CD}
{\Bbb R}^2 @> L >> {\Bbb R}^2 \\
@VV{\pi}V @VV{\pi'}V \\
{\Bbb R}^2 / \Lambda @> {\tilde L} >> {\Bbb R}^2 / {\ZZ}^2 
\end{CD}
\end{equation}
Note that the slope of $L(1)$ is irrational by $(\ast \ast )$. So the sequence $\{ \pi' (L(n))\}_{n\geq 0} $ is dense in ${\Bbb R}^2 / {\ZZ}^2$. Pulling
back this sequence to the torus ${\Bbb R}^2 / \Lambda $ by $\tilde L$, it follows from diagram (1) that the sequence  $\{ \pi (n) \}_{n\geq 0} $ is dense in ${\Bbb R}^2 / \Lambda $. Hence if
$(\lambda _1, \lambda _2 )$ satisfies $(\ast )$ and $(\ast \ast )$, the subgroup generated
by $\{1, \lambda _1, \lambda _2 \}$ is dense in $\CC$. Finally, it is straightforward to
check that $(\ast )$ and $(\ast \ast )$ hold for almost every $(\lambda _1, \lambda _2 ) \in \CC^2$. $\hfill \Box $\\

Now the first version of the density theorem is quite easy to prove.\\ \\
{\bf 3.21 Theorem } \ {\it For a typical ${\cal F} \in {\cal A}_n$, all but at most $n+2$ non-singular leaves are dense in $\CC \PP^2$.}\\ \\
{\bf Proof.} Let ${\cal F} \in {\cal A}_n$ has properties (i), (ii), and (iii) of Proposition 3.20.
Since all singular points $\{p_1, \ldots , p_{n+1} \} = L_0 \cap$ sing$({\cal F})$ are of hyperbolic type,
there are exactly two local separatrices through each $p_j$, one of them lies in ${\cal L}_\infty $.
Denote by ${\cal L}_j$ the {\it global separatrix} through $p_j$ transversal to ${\cal L}_\infty $ (by the global separatrix we mean the continuation of the local
separatrix as a leaf). Note that we might have ${\cal L}_i = {\cal L}_j$ even if $p_i\neq p_j$.

Now let $p\in \CC \PP^2 \sm \{{\rm sing} ({\cal F}) \cup {\cal L}_1 \cup \ldots \cup {\cal L}_{n+1} \cup {\cal L}_\infty \}$.
By Corollary 2.13, ${\cal L}_p$ has a point of accumulation on $ L_0$. By the
choice of $p$ and Corollary 3.18, this point may be
chosen on ${\cal L}_\infty $ , hence by Proposition 1.5 the whole ${\cal L}_\infty $
is in the closure of ${\cal L}_p$. In particular, ${\cal L}_p$ intersects the section $\Sigma $
which is transversal to ${\cal L}_\infty $ at the base point $a$ of $\pi _1 ({\cal L} _\infty )$
(see 2.14). By Theorem 3.5, ${\cal L}_p \cap \Sigma $ is dense in some neighborhood $\Omega \subset \Sigma $ of $a$.

Now choose an open set $U\subset \CC \PP^2$ and a point $q \in U \sm \{ {\rm sing} ({\cal F}) \cup {\cal L}_1 \cup \ldots \cup {\cal L}_{n+1} \cup {\cal L}_\infty \}.$
By the above argument, ${\cal L}_q \cap \Sigma $ is dense in $\Omega $. Therefore
$q\in \overline {{\cal L}}_p$ by another application of Proposition 1.5. Since $U$ was arbitrary,
${\cal L}_p$ is dense in $\CC \PP^2$. $\hfill \Box $\\

The above proof showed that all non-singular leaves other than ${\cal L}_\infty$
and the global separatrices ${\cal L}_j$'s are dense. The condition $p \not \in \bigcup _{j=1} ^{n+1} {\cal L}_j$
was used where we were looking for a point in ${\cal L}_\infty \cap \overline {{\cal L}}_p$. If the only
accumulation points of ${\cal L}_p$ on $L_0$ are singular, it turns out that $\overline {{\cal L}}_p$ is
an {\it algebraic curve}. Since the algebraic leaves other than ${\cal L}_\infty $ are not typical for
elements of ${\cal A}_n$ (cf. Proposition 3.22 below), we will be able to prove a stronger version of the density theorem without the assumption $p \not \in \bigcup _{j=1}^{n+1} {\cal L}_j$ in the above proof.\\ \\
{\bf 3.22 Proposition} \ {\it For a typical ${\cal F} \in {\cal A}_n$,
the only algebraic leaf is ${\cal L}_\infty$.}\\ \\
{\bf Proof.} The argument below is an adaptation of an idea due to I. Petrovski\u\i \ and E. Landis \cite{P-L1}. Let ${\cal F}: \{Pdy -Qdx=0\} \in {\cal A}_n$.
Suppose that the algebraic curve $S_K:\{K=0 \}$ (see 1.8) with singular points of ${\cal F}$ deleted is a leaf of ${\cal F}$,
where $K=K(x,y)$ is an irreducible polynomial of degree $k$. Since $S_K$ is a leaf, we have
$$\frac {\partial K}{\partial x} (x,y) P(x,y)+ \frac {\partial K}{\partial y} (x,y) Q(x,y)=0$$
whenever $K(x,y)=0$. It follows that there exists a polynomial $\tilde {K}$ of degree at most
$(n-1)$ such that
\setcounter{equation}{0}
\begin{equation}
\frac {\partial K}{\partial x} P + \frac {\partial K}{\partial y} Q = K \tilde {K}.
\end{equation}
Conversely, if there exist polynomials $K, \tilde {K} , P,$ and $Q$ satisfying (1) with
$K$ irreducible and $P$ and $Q$ relatively prime, then $S_K \sm$sing$({\cal F})$ is an algebraic leaf of ${\cal F}: \{ Pdy -Qdx =0 \}.$

Let $E$ be the complex linear space of the coefficients of $K, \tilde {K}, P,$ and $Q$, which has dimension
$$\frac {(k+1)(k+2)}{2}+\frac {n(n+1)}{2}+ 2 \ \frac {(n+1)(n+2)}{2} = \frac {1}{2} \{ (k+1) (k+2) +(3n^2 +7n +4 )\}.$$
If we impose (1) on these coefficients, we obtain equalities that define an algebraic variety $S$ in $E$. Note that if $a\in S$, so does $\lambda a$ for every $\lambda \in
\CC^{\ast}$ by (1). Therefore, $S$ projects to an algebraic variety $S^{\ast}$
in $\CC \PP^d$, with $d={\rm dim}E-1$. Decompose $S^{\ast}$ as $\bigcup_{j=1}^m S_j^{\ast}$,
where each $S_j^{\ast}$ is irreducible. Let $E'\subset E$ be the subspace of
coefficients of $P$ and $Q$, which has dimension $n^2+3n+2$. The linear projection
$\pi :E\rightarrow E'$ induces a projection $\CC \PP^d \rightarrow \CC \PP^N$
(still denoted by $\pi$), where $N=({\rm dim}E')-1=n^2+3n+1$. Each $\pi (S_j^{\ast})$
is an algebraic variety in $\CC \PP^N$. Since there are SHFC's in ${\cal A}_n$
which do not have any algebraic leaf other than ${\cal L}_\infty$, we have
$\pi (S_j^{\ast})\neq \CC \PP^N$, so dim $\pi (S_j^{\ast}) \leq N-1$. Taking
the union for all $j=1,\cdots ,m$, we obtain the algebraic variety $\pi (S^{\ast})$
in $\CC \PP^N$, each irreducible component of which has dimension $\leq N-1$.
Each point in $\pi (S^{\ast}) \cap {\cal A}_n$ corresponds to a SHFC having an algebraic leaf other
than ${\cal L}_\infty$. It follows that the Lebesgue measure of $\pi (S^{\ast})$ in ${\cal A}_n$ is zero, and
we are done. $\hfill \Box $\\ \\
{\bf 3.23 Remark } By a much more difficult argument, using an index theorem of
Camacho and Sad and the concept of the Milnor number of a local branch of a
singular point, A. Lins Neto has shown that for $n\geq 2$ there exists an open and dense
subset of ${\cal D}_n$ (see 1.24) consisting of SHFC's which do not have any
algebraic leaf \cite{L}.\\ \\
{\bf 3.24 Proposition } {\it Let ${\cal F} \in {\cal A}_n'$ and all points in $L_0 \cap {\rm sing} ({\cal F})$ be
of hyperbolic type. Let ${\cal L}$ be a non-singular leaf of ${\cal F}$ such that $\overline {\cal L} \cap L_0$
consists of singular points only. Then ${\cal L}$ is an algebraic leaf.}\\ \\
{\bf Proof.} Let $\Omega =\CC \PP^2 \sm {\rm sing}({\cal F})$. First we show that
${\cal L}$ is closed in $\Omega $. Since $\overline {{\cal L}}\cap L_0 \subset L_0 \cap {\rm sing} ({\cal F})$,
it follows that ${\cal L}={\cal L}_j$, the global separatrix transversal to ${\cal L}_\infty$ through some singular point
$p_j$ on $L_0$. Let $p\in \overline {{\cal L}}$ be any
non-singular point. Then $L_0 \cap \overline{{\cal L}}_p \subset L_0 \cap \overline {{\cal L}} \subset L_0 \cap {\rm sing}({\cal L})$, which
shows that ${\cal L}_p$ also coincides with ${\cal L}_j$. Hence ${\cal L}_p ={\cal L}$ and $p\in {\cal L}$.

Next we show that $\overline {{\cal L}}$ is an analytic subvariety of $\CC \PP^2$. Let $p\in \Omega $. If
$p\in \overline {{\cal L}}$, then actually $p \in {\cal L}$. Suppose that $(U, \varphi )$ is a foliation chart
around $p$, $\Sigma $ is a section transversal to ${\cal L}_p= {\cal L}$ at $p$, $\Sigma '$ is another
section transversal to ${\cal L}$ at $p'$ near $p_j \in L_0 \cap {\rm sing} ({\cal F})$, and
$\gamma $ is any path in ${\cal L}$ joining $p$ to $p'$ (Fig. 14). Let $f_\gamma : \Sigma \rightarrow \Sigma '$
be the associated holonomy mapping. If there exists a sequence $p_n\in {\cal L} \cap \Sigma $ which
converges to $p$, then by considering the sequence $f_\gamma (p_n) \in {\cal L} \cap \Sigma '$
we conclude from Corollary 3.18 that ${\cal L}$ must have an accumulation
point on ${\cal L}_\infty $, which contradicts our assumption. Therefore, by choosing $U$
small enough, the only plaque of ${\cal L}$ in $U$ is the one which passes through $p$,   
and evidently there exists a holomorphic function $f:U\rightarrow \CC$
such that $f^{-1}(0) ={\cal L} \cap U$.

\realfig{figthesis13}{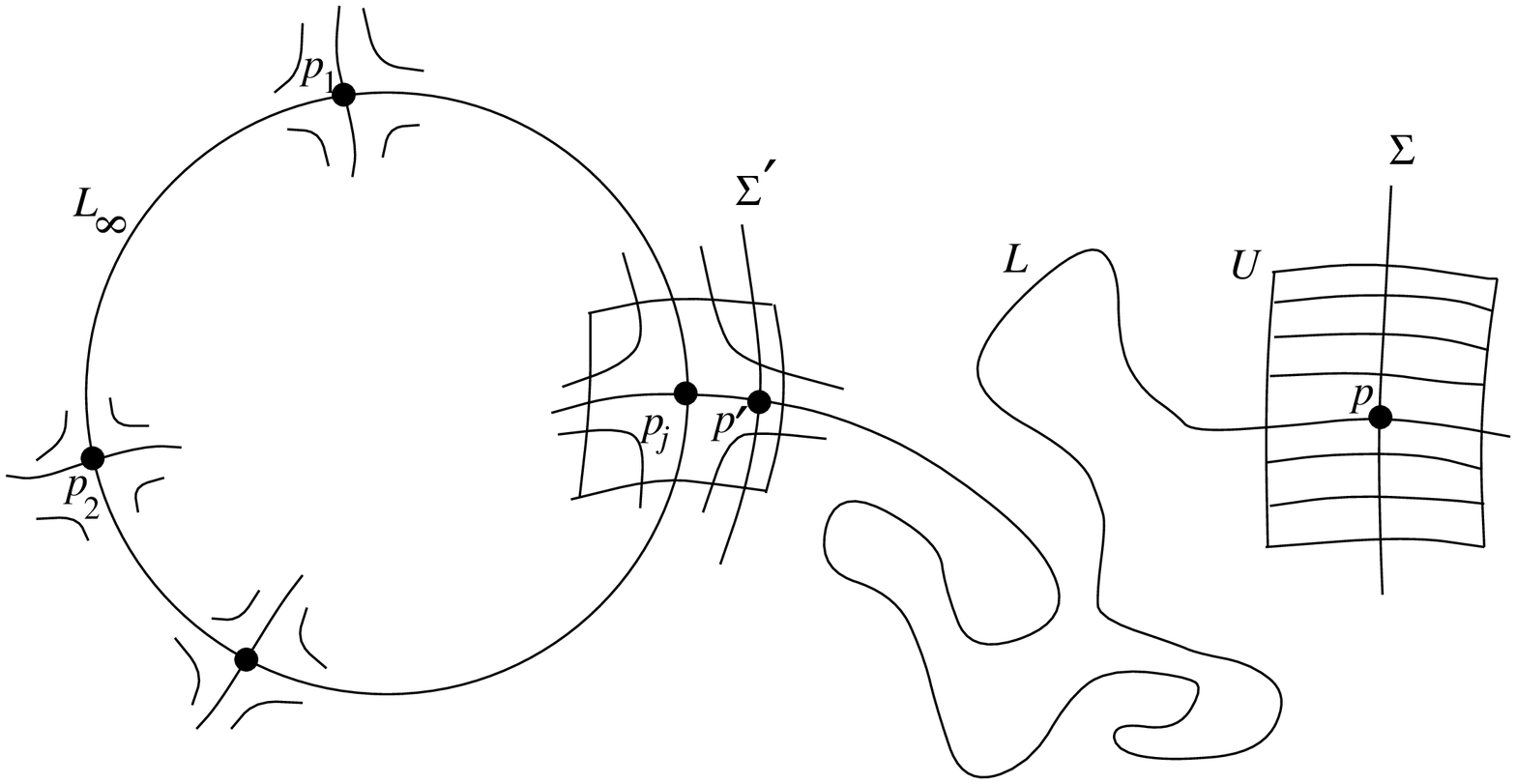}{}{10cm}


This means that $\ov{\cal L} \sm \mbox{sing}({\cal F})$ is a $1$-dimensional analytic subvariety of $\Omega$. Therefore, since dim(sing(${\cal F}))=0<1=$ dim($\ov{\cal L} \sm \mbox{sing}({\cal F})$), the well-known theorem of Remmert-Stein (see \cite{G-R}) shows that $\overline {{\cal L}}$ is an
analytic subvariety of $\CC \PP^2$. Finally, every analytic subvariety of
$\CC \PP^2$ is algebraic by the Chow's Theorem \cite{G-H}. $\hfill \Box $\\

Now, by Theorem 3.21, Proposition 3.22 and Proposition 3.24, we conclude
the density theorem of Khudai-Veronov.\\ \\
{\bf 3.25 Theorem }  {\it For a typical ${\cal F} \in {\cal A}_n$, all non-singular leaves, except the
leaf at infinity, are dense in $\CC \PP^2$}. $\hfill \Box $\\ \\
\vspace {4 mm}
{\large \bf Ergodicity of SHFC's on $\CC \PP^2$}\\
In what follows, we say that a set $A\subset \CC \PP^2$ has measure zero if for every chart $(U, \varphi)$ compatible with the standard smooth structure of $\CC \PP^2$, the set $\varphi(A\cap U) \subset \RR^4$ has Lebesgue measure zero. In other words, we consider the {\it Lebesgue measure class} on $\CC \PP^2$.\\ \\ 
{\bf 3.26 Definition}\ Let ${\cal F} \in {\cal A}_n$. The {\it saturation} $s(A)$ of
$A\subset \CC \PP^2$ is the set of all points $q$ such that $q\in {\cal L}_p$
for some $p\in A$. A set $A$ is called {\it saturated} if $s(A)=A$ up to a set of measure zero. ${\cal F}$ is
called {\it ergodic} if for every measurable saturated set $A$, either $A$
or $\CC \PP^2\sm A$ has measure zero.\\

Evidently if ${\cal F}$ is ergodic, then every non-singular leaf of ${\cal F}$
is either dense in $\CC \PP^2$ or its closure has measure zero. Also the following
observation is useful: Let ${\cal L}_p$ be a non-singular leaf and $\Sigma $
be a section transversal to ${\cal L}_p$ at $p$. Suppose that $A\subset \Sigma $
be a set of measure zero (with respect to the Lebesgue measure class on $\Sigma \simeq \DD)$.
Then $s(A)$ has measure zero in $\CC \PP^2$. This is because $s(A)$ can be
covered by a countable number of foliation charts $\{(U_i, \varphi _i)\}_{i=1} ^\infty$
and each $s(A) \cap U_i$ has measure zero. In particular, each single leaf has measure zero (take $A=\{p\}$).\\ \\
{\bf 3.27 Theorem } \ {\it Let ${\cal F} \in {\cal A}_n$ have properties (i),(ii), and (iii) of Proposition 3.20. Then ${\cal F}$ is ergodic. In particular, ergodicity is typical for elements of ${\cal A}_n$.}\\ \\
{\bf Proof.} Let $A$ be a measurable saturated subset of $\CC \PP^2$. Without loss of generality, we may assume that $A$ does not contain any global separatrix ${\cal L}_j$ through
$p_j \in L_0 \cap {\rm sing}({\cal F})$, for each individual leaf has measure zero. Let
$p\in A$ be non-singular. Then ${\cal L}_p$ must accumulate on ${\cal L}_\infty$, so it
has to intersect the section $\Sigma $ which is transversal to ${\cal L} _\infty $ at the base point $a$.
By Theorem 3.15, $G_\infty $ is ergodic. Since $A$ is saturated, $A\cap \Sigma $ is $G_\infty$-invariant, so there
is an open disk $U\subset \Sigma $ around $a$ such that either $A\cap U$ or $U\sm A$ has measure zero (with respect to the Lebesgue measure class on $\Sigma $ ). It is clear that
$$s(A\cap U)=A \sm {\rm sing} ({\cal F}) \ \ {\rm and}\  \ s(U \sm A) =\CC \PP^2 \sm \{A \cup {\rm sing } ({\cal F}) \}.$$
By the observation before the statement of the theorem, it follows that either $A$
or $\CC \PP^2 \sm A$ has measure zero. $\hfill \Box$
\newpage
\vspace*{4 mm}
\noindent
{\large {\bf Chapter 4}}\vspace{4 mm}\\
{\Large{\bf Non-Trivial Minimal Sets}}\\ \\ \\ \\ \\ \\ \\ \\ \\ \\ \\ \\
\thispagestyle{plain}
\noindent
This chapter deals with a somewhat different global aspect of SHFC's
on $\CC \PP^2$. As will be seen, the foliations under consideration
are essentially those which do not have any algebraic leaf. In particular, because of the absence of the leaf at infinity, we cannot utilize such powerful tools as the monodromy group $G_\infty$. Recall that a typical ${\cal F}\in {\cal A}_n$ has at least one algebraic leaf (i.e., ${\cal L}_\infty )$. Hence from the point of view of differential equations for which the decomposition into the ${\cal A}_n$ is more natural, the foliations we consider in this chapter almost never occur. However, from the point of view of foliation theory, for which the natural decomposition is by the ${\cal D}_n$, the property of having no algebraic leaf seems to be typical (see Remark 3.23). 

The study of limit sets of foliations and flows in the real domain has
proved to be of great significance in those theories. The classical theorem of Poincar\'e-Bendixson asserts that for every smooth real flow on the $2$-sphere, every trajectory accumulates either on a periodic orbit or a singular point (or both). It is natural to ask a similar question for SHFC's on $\CC \PP^2$. Here the analogue of a periodic orbit is a compact non-singular leaf and it is not difficult to prove that no such leaves could exist (Theorem 4.12). So we naturally arrive at the following question, apparently first asked by C. Camacho:\vspace{0.1in}\\
{\bf Question.} {\it Is there a non-singular leaf of a SHFC on $\CC \PP^2$ which does not accumulate on any singular point?}\vspace{0.1in}\\
Oddly enough, the question has remained open since the mid $80$'s. One can formulate it in a slightly different language, commonly used in foliation theory.\\ \\
{\bf 4.1 Definition} A {\it minimal set} for a SHFC on $\CC \PP^2$ is a compact
saturated non-empty subset of $\CC \PP^2$ which is minimal with respect to these three properties. A {\it non-trivial minimal set} is a minimal set which is not a singular point. Throughout this chapter, ${\cal M}$ will always denote a non-trivial minimal set.\\

Minimality shows that if $p\in {\cal M}$, then $\overline {{\cal L}}_p= {\cal M}$. It follows that the problem of finding a non-singular leaf which does not accumulate on any singular point is equivalent to finding a non-trivial minimal set. Therefore, we can reformulate the above question as\vspace{0.1in}\\
{\bf The Minimal Set Problem.} {\it Does there exist a SHFC on $\CC \PP^2$ which has a non-trivial minimal set?} \vspace{0.1in}\\
Such a non-trivial minimal set is an example of a {\it Riemann surface lamination}. By definition, a Riemann surface lamination (RSL) is a compact space which locally looks like the product of the unit disk and a compact metric space (usually a Cantor set). The transition maps between various charts are required to be holomorphic in the leaf direction and only continuous in the transverse direction. Clearly every compact Riemann surface is such a space, but they form the class of {\it trivial} RSL's. Although there are some basic results on uniformization of RSL's (see \cite{Candel}), the corresponding embedding problem is rather unexplored. A classical theorem asserts that every compact Riemann surface can be holomorphically embedded in $\CC \PP^3$. The Minimal Set Problem, as E. Ghys has suggested, could be viewed as a special case of the embedding problem for RSL's: ``Can a non-trivial Riemann surface lamination be holomorphically embedded in $\CC \PP^2$?''
     
In what follows we show some basic properties of non-trivial minimal sets. This theory, due to Camacho-Lins Neto-Sad \cite{C-L-S1}, was developed in part in the hope of arriving at a contradiction to the existence of non-trivial minimal sets.\\ \\
\vspace {4 mm}
\noindent
{\large \bf Uniqueness of Minimal Sets}\\
How many distinct non-trivial minimal sets, if any at all, can a SHFC on $\CC \PP^2$ have?\\ \\
{\bf 4.2 Theorem } {\it A SHFC on $\CC \PP^2$ has at most one non-trivial minimal set.}\\ \\
The proof of this nice fact is quite elementary, and is based on the study of
the distance between two non-singular leaves and the application of the Maximum Principle for real harmonic functions (cf. \cite{C-L-S1}).

To study the distance between two leaves, we have to choose a suitable Riemannian
metric on $\CC \PP^2$. Consider the Hermitian metric
\setcounter{equation}{0}
\begin{equation}
ds^2=\frac {|dx|^2+|dy|^2+|xdy-ydx|^2}{(1+|x|^2+|y|^2)^2}
\end{equation}
in the affine chart $(x,y)\in U_0$, which extends to a Hermitian metric on
the whole projective plane. It is called the {\it Fubini-Study metric} on $\CC \PP^2$ \cite{G-H}. We will denote by $d$ the Riemannian distance induced by this metric.

The associated $(1,1)$-form of the Fubini-Study metric is of the form
\begin{equation}
\begin{array}{rl}
\Omega & =\displaystyle{\frac{\sqrt{-1}}{2 \pi}\  \frac{dx \wedge d\ov{x}+dy \wedge d\ov{y}+(x dy-y dx)\wedge (\ov{x} d\ov{y}-\ov{y} d\ov{x})}{(1+|x|^2+|y|^2)^2}}\\
       & =\displaystyle{\frac{\sqrt{-1}}{2 \pi}\  \partial \ov{\partial} \log(1+|x|^2+|y|^2)}.
\end{array}
\end{equation}

\noindent
{\bf 4.3 Lemma} {\it For any $p_0=(x_0,y_0)$ and $p_1=(x_1,y_1)$ in the affine chart $(x,y)\in U_0$,
we have}
$$d(p_0,p_1) \leq \frac {|p_0-p_1|}{(1+\delta ^2(p_0,p_1))^{1/2} },$$
{\it where $\delta (p_0,p_1)$ is the minimum (Euclidean) distance form the origin of the line segment which joins $p_0$ and $p_1$ .}\\ \\
{\bf Proof.} By definition, $d(p_0,p_1)=\inf _\gamma \{\int_0^1 \| \gamma '(t) \| dt \}$,
where the infimum is taken over all piecewise smooth curves $\gamma :[0,1] \rightarrow \CC \PP^2$ with
$\gamma (0)=p_0, \gamma (1)=p_1$. In particular, when $\gamma (t)=(1-t)p_0+tp_1=:(x(t),y(t))$, one has
$$d^2 (p_0,p_1) \leq (\int^1_0 \| \gamma '(t) \| dt)^2 \leq \int ^1_0 \| \gamma ' (t) \| ^2 dt.$$
Now we estimate $\| \gamma ' \| ^2$:
$$\begin{array}{ll}\vspace{4 mm}
\| \gamma' (t) \|^2 & = \displaystyle { \frac {|x'(t)|^2+|y'(t)|^2+|x(t)y'(t)
-y(t)x'(t)|^2}{(1+|x(t)|^2+|y(t)|^2)^2}}\\ \vspace{4 mm}
 & \leq \displaystyle { \frac {|x'(t)|^2+|y'(t)|^2+(|x(t)|^2+|y(t)|^2)(|x'(t)|^2+|y'(t)|^2)}{(1+|x(t)|^2+|y(t)|^2)^2}}\\  \vspace{4 mm}
 & = \displaystyle {\frac{(|x'(t)|^2+|y'(t)|^2)(1+|x(t)|^2+|y(t)|^2)}{(1+|x(t)|^2+|y(t)|^2)^2}}\\  \vspace{4 mm}
 & = \displaystyle {\frac{|x_0-x_1|^2+|y_0-y_1|^2}{1+|x(t)|^2+|y(t)|^2}}\\  \vspace{3 mm}
 & \leq \displaystyle {\frac {|p_0-p_1|^2}{1+\delta^2(p_0,p_1)}},
\end{array}$$
and this completes the proof. $\hfill \Box$\\ \\
{\bf 4.4 Corollary} {\it Let E and F be two disjoint compact subsets of $\CC \PP^2$, and
$E':=E\cap U_0$ and $F':=F\cap U_0$ be both non-empty. If $\epsilon :=\inf \{|p-q|:(p,q)\in E' \times F' \}$, then
$\epsilon > 0$ and there exists a pair $(p,q) \in E' \times F' $ with $|p-q|=\epsilon $.}\\ \\
{\bf Proof.} Let $(p_n,q_n)\in E' \times F'$ be such that $|p_n-q_n| \rightarrow \epsilon $ as
$n\rightarrow \infty$. By taking subsequences, if necessary, we may assume that $p_n\rightarrow p\in E$ and
$q_n\rightarrow q \in F$. If $(p,q) \in E' \times F'$, we are done. Otherwise, if
$p \in E\sm E'$, one has $q\in F\sm F'$ since $\{|p_n-q_n| \}$
is bounded. Therefore $\delta (p_n,q_n)\rightarrow \infty $ as $n\rightarrow \infty $, so that
$d(p,q)=0$ by Lemma 4.3, which is a contradiction. $\hfill \Box $\\ \\
{\bf Proof of Theorem 4.2.} Let ${\cal M}_1$ and ${\cal M}_2$ be two non-trivial minimal sets
of ${\cal F}:\{Pdy-Qdx=0\}$. They are necessarily disjoint by minimality. Set
${\cal M}'_1:={\cal M}_1\cap U_0$, ${\cal M}'_2:={\cal M}_2\cap U_0$ and
$\epsilon :=\inf \{|p-q|:(p,q) \in {\cal M}_1' \times {\cal M}'_2 \}$. By
Corollary 4.4, there exists $(p,q)\in {\cal M}'_1\times {\cal M}'_2$, with
$|p-q|=\epsilon$. For simplicity, let $p=(0,0)$ and $q=(0, y_0)$, with $|y_0|=\epsilon$.
It follows from the definition of $\epsilon $ that the $y$-axis is normal to ${\cal L}_p$ and
${\cal L}_q$ at $p$ and $q$. Since $p$ and $q$ are not singular, we can parametrize the
leaves by the $x$-parameter in a disk $\DD(0,r)$ around the origin:
$${\cal L}_p:\ \ \ x\mapsto y_p(x),\ \ \ y'_p(x)= \frac {Q(x,y_p(x))}{P(x,y_p(x))},\ \ \ y_p(0)=0,$$
$${\cal L}_q:\ \ \ x\mapsto y_q(x),\ \ \ y'_q(x)= \frac {Q(x,y_q(x))}{P(x,y_q(x))},\ \ \ y_q(0)=y_0.$$
Define $h: \DD(0,r) \rightarrow {\Bbb R}$ by $h(x)=\log |y_p(x)-y_q(x)|.$
This is a harmonic function with a minimum at $x=0$. Therefore $h\equiv \log \epsilon $ on $\DD(0,r)$, so that locally ${\cal L}_q$ is just the translation of ${\cal L}_p$ by
$y_0$. By analytic continuation, this is true globally, i.e., ${\cal L}_q\cap
U_0=({\cal L}_p\cap U_0)+(0,y_0)$.
By Corollary 2.12, there exists a sequence $q_n \in {\cal L}_q \cap U_0$
tending to infinity. The sequence $p_n:=q_n-(0,y_0)\in {\cal L}_p \cap U_0$
also converges to infinity, so that $d(p_n,q_n)\rightarrow 0$ by Lemma 4.3. This shows that
$\overline {{\cal L}}_p \cap \overline {{\cal L}}_q={\cal M}_1 \cap  {\cal M}_2 \neq \emptyset$,
which is a contradiction. $\hfill \Box $\\ \\
\vspace {4 mm}
\noindent
{\large \bf Hyperbolicity in Minimal Sets}\\
Our next goal is to determine the type of a leaf (as a Riemann surface) in the non-trivial minimal set. First we need some preliminaries.\\ \\
{\bf 4.5 Uniformization of Riemann Surfaces}\ According to the celebrated Uniformization Theorem of Koebe-Poincar\'e-Riemann, every simply-connected Riemann surface
is biholomorphic to one of the three standard models: The Riemann sphere $\overline{\CC}$,
the complex plane $\CC$, or the unit disk $\DD$ \cite{Ah}. Since every Riemann
surface has a holomorphic universal covering, it follows that {\it every Riemann surface
can be covered holomorphically by $\overline{\CC}$, $\CC$, or $\DD$, in which case it is called
spherical, Euclidean, or hyperbolic, respectively}. Elementary considerations on the group of
covering transformations acting on the universal covering shows that every spherical Riemann surface is biholomorphic to $\overline{\CC}$, and the Euclidean
ones are biholomorphic to $\CC$, a $2$-torus, or the cylinder $( \simeq \CC^\ast )$. The last two surfaces are biholomorphic to quotients
of $\CC$ under the action of the groups generated by
$$ \{z\mapsto z+a, z\mapsto z+b, \ {\rm for \ some } \ a,b\in \CC^\ast \ {\rm with} \ a/b\not \in {\Bbb R}\},$$
and
$$ \{z\mapsto z+a, \ {\rm for\ some} \ a\in \CC^\ast \},$$
respectively. All other Riemann surfaces are thus hyperbolic.\\ \\
{\bf 4.6 Invariant Metrics on Hyperbolic Surfaces} \ The Riemannian metric
\setcounter{equation}{0}
\begin{equation}
\rho_\DD:=\frac{4}{(1-|z|^2)^2} \ |dz|^2,
\end{equation}
called the {\it Poincar\'e} or {\it hyperbolic metric} on the unit disk, is the unique (up to
multiplication by a non-zero constant) conformal metric on $\DD$ which is invariant under all biholomorphisms $\varphi :\DD \rightarrow \DD$. All biholomorphisms of $\DD$ are therefore isometries with respect to the distance induced
by $\rho_\DD$. The unit disk equipped with this metric is a complete metric space \cite{M1}.

Now let $X$ be a hyperbolic Riemann surface, with the covering map $\pi :\DD
\rightarrow X$. Since the Poincar\'e metric $\rho_\DD$ is invariant under
all the covering transformations, $\rho_{\DD}$ induces a Poincar\'e metric 
$\rho_X$ on $X$ which is conformal and invariant under all biholomorphisms $X\rightarrow X$. By the very definition of $\rho_X$, the projection $\pi $ is a local isometry.

It is a direct consequence of the Schwarz Lemma that if $\varphi :X\rightarrow Y$ is a holomorphic mapping between hyperbolic Riemann surfaces, then $\varphi $
decreases the Poincar\'e distance, i.e., for every $x,y\in X$,
$$d_Y(\varphi(x),\varphi(y))\leq d_X(x,y),$$
where $d_X$ and $d_Y$ are the Riemannian distances induced by $\rho_X$ and $\rho_Y$, respectively. If the equality holds for a pair $(x,y)$, then $\varphi $ will be a local isometry \cite{M1}.\\ \\
{\bf 4.7 Curvature of Conformal Metrics} \ Recall that the {\it Gaussian curvature} $\kappa $ of a conformal metric $ds^2=h^2|dz|^2$ on a Riemann surface 
is given by
\setcounter{equation}{0}
\begin{equation}
\kappa (z)= - \frac {(\Delta \log h)(z)}{h^2(z)},
\end{equation}
where, as usual, $z$ denotes a local chart on the surface. It follows that the {\it Gaussian curvature is a conformal invariant}, that is, if
$\varphi : X \rightarrow Y$ is a holomorphic map between Riemann surfaces, and
if $ds^2$ is a conformal metric on $Y$, then at any point $z\in X$ for which $\varphi'(z)\neq 0$, the curvature $\kappa '$ at $z$ of the
pull-back metric on $X$ is equal to the curvature $\kappa$ of $ds^2$ at $\varphi (z)$.

It follows from 4.6(1) and 4.7(1) that the Poincar\'e metric $\rho_\DD$ on the
unit disk has constant Gaussian curvature $-1$. The same is true for every hyperbolic surface
equipped with the Poincar\'e metric since the curvature is a conformal invariant.

The fact that hyperbolic Riemann surfaces admit a metric of strictly negative curvature is a
characteristic property, as can be seen by the following \\ \\
{\bf 4.8 Theorem} {\it Suppose that $X$ is a Riemann surface that has a conformal metric whose associated curvature $\kappa$ satisfies $\kappa (z) <\sigma <0$ for all $z\in X$ and some constant $\sigma$. Then X is hyperbolic.} $\hfill \Box $\\

In fact, this is a special case of a more general fact: A complex manifold $M$
which admits a distance for which every holomorphic mapping $\DD\rightarrow M$ is distance-decreasing
is hyperbolic in the sense of Kobayashi. For Riemann surfaces, the usual notion
of hyperbolicity is equivalent to the hyperbolicity in the sense of Kobayashi.
By Ahlfors' generalized version of the Schwarz Lemma, a Riemann surface $X$ which has a conformal metric of strictly negative curvature admits a distance for which every holomorphic mapping $\DD \rightarrow X$ is distance-decreasing (see \cite{Kob}, and also \cite{Kr} for a nice exposition in the case of domains in $\CC$).\\ \\
{\bf 4.9 A Hermitian Metric on $\CC \PP^2$} \ We now construct a Hermitian metric on $\CC\PP^2\sm $sing$({\cal F})$ which induces a conformal metric of negative Gaussian curvature on each non-singular leaf of a given SHFC ${\cal F}$. The metric is a modification of the Fubini-Study metric 4.2(1).

Suppose that ${\cal F}: \{ \omega =Pdy-Qdx=0 \} \in {\cal D}_n$, and let $R=yP-xQ$. Consider
the following Hermitian metric on $U_0 \sm {\rm sing}({\cal F})$:
\setcounter{equation}{0}
\begin{equation}
\rho :=(1+|x|^2+|y|^2)^{n-1} \frac { |dx|^2+|dy|^2+|xdy-ydx|^2}{|P(x,y)|^2+|Q(x,y)|^2+|R(x,y)|^2}.
\end{equation}
$\rho $ extends to a Hermitian metric on $\CC \PP^2 \sm {\rm sing}({\cal F})$. To
see this, let us for example compute the extended metric on the affine chart $(u,v)\in U_1$
(cf. 1.7(3)):
$$\begin{array}{ll}  \vspace{4 mm}
(\phi_{10}^\ast \rho )(u,v) & \displaystyle{=(1+\frac {1}{|u|^2} +\frac {|v|^2}{|u|^2})^{n-1}\frac {|u|^{-4}(|du|^2+|udv-vdu|^2+|dv|^2)}
{|P(\frac{1}{u} , \frac{v}{u})|^2 +|Q(\frac{1}{u} , \frac{v}{u})|^2+|R(\frac{1}{u} , \frac{v}{u})|^2}}\\
 & =(1+|u|^2+|v|^2)^{n-1} {\displaystyle \frac {|du|^2 +|dv|^2 +|udv-vdu|^2}{|\tilde {P} (u,v)|^2 +|\tilde {Q} (u,v)|^2 +|\tilde {R} (u,v)|^2}},
\end{array}$$
where $\tilde {P}(u,v)=u^{n+1} \displaystyle{P(\frac{1}{u} , \frac{v}{u})}$, $\tilde {Q}(u,v) =u^{n+1} \displaystyle{Q(\frac{1}{u} , \frac{v}{u})}$, and $\tilde {R} (u,v) =u^{n+1} \displaystyle{R(\frac{1}{u} , \frac{v}{u})}$ are {\it polynomials} in $u,v$.

Now let $p\in U_0$ be a non-singular point of ${\cal F}$, and let $T \stackrel {\eta }{\longmapsto } (x(T), y(T))$
be a local parametrization of ${\cal L}_p$ near $p$ with $\eta (0)=p$. By (1) above, the
induced conformal metric on ${\cal L}_p$ is
$$ds^2 = (1+|x(T)|^2 +|y(T)|^2)^{n-1} \displaystyle {\frac{|x'(T)|^2+|y'(T)|^2+|x(T)y'(T)-y(T)x'(T)|^2}{|P(\eta (T))|^2+|Q(\eta (T))|^2+|R(\eta (T))|^2} |dT|^2}$$
$$\begin{array}{ll} \vspace{3 mm}
 & =(1+|x(T)|^2 +|y(T)|^2)^{n-1} |dT|^2  \\
 & =:h^2(T)|dT|^2.
\end{array}$$
The Gaussian curvature of ${\cal L}_p$ at $p$ is
$$\kappa (p) = -\frac {(\Delta \log h) (0)}{h^2(0)}$$
by 4.7(1) and conformal invariance of the curvature. Computation gives
\begin{equation}
\begin{array}{ll}\vspace{2 mm}
\kappa (p) &= \displaystyle {\frac {-2}{(1+|p|^2)^{n-1}} \left ( \frac {\partial }{\partial T} \frac {\partial } {\partial \overline {T}} \log h^2 (T) \right ) \rule[-2.5 mm]{.1 mm}{8 mm} _{\ T=0}}\\ \vspace{2 mm}
 &= \displaystyle {\frac {-2(n-1)}{(1+|p|^2)^{n-1}} \left ( \frac {\partial }{\partial T} \frac {\partial } {\partial \overline {T}} \ {\rm log} (1+|x(T)|^2 +|y(T)|^2) \right )\rule[-2.5 mm]{.1 mm}{8 mm}_{\ T=0}}\\ \vspace{2 mm}
 &= \displaystyle {\frac {-2(n-1)}{(1+|p|^2)^{n-1}} \left ( \frac {|x'(T)|^2+|y'(T)|^2+|x(T)y'(T)-y(T)x'(T)|^2}{(1+|x(T)|^2+|y(T)|^2)^2} \right )\rule[-2.5 mm]{.1 mm}{8 mm}_{\ T=0}}\\ \vspace{2 mm}
 &= \displaystyle {\frac {-2(n-1)}{(1+|p|^2)^{n+1}} (|P(p)|^2 +|Q(p)|^2 +|R(p)|^2),}
\end{array}
\end{equation}
which is strictly negative.

Now let ${\cal F}\in {\cal D}_n$ have a non-trivial minimal set ${\cal M}$, and $p\in {\cal M}.$
As (2) is a continuous function of $p$ which extends to the whole $\CC \PP^2 \sm {\rm sing} ({\cal F})$,
it follows that the Gaussian curvature of the induced metric on ${\cal L}_p$ is
bounded from above by a strictly negative constant. By Theorem 4.8, we have\\ \\
{\bf 4.10 Theorem} \ {\it Every leaf contained in the non-trivial minimal set is a hyperbolic Riemann surface.} $\hfill \Box $\\ \\
{\bf 4.11 Example}\ We can use the preceding result to show that no SHFC ${\cal F}$
of geometric degree 1 can have a non-trivial minimal set.
Let $\cal L$ be a leaf of $\cal F$ contained in the non-trivial minimal set $\cal M$.
Choose $p\in {\cal L}$ and let $T\mapsto \eta (T)$ be a local parametrization
of $\cal L$ near $p$, with $\eta (0)=p$. The germ of $\eta$ can be analytically
continued over the whole plane $\CC$. To see this, observe that ${\cal F}$ is
induced by a holomorphic vector field {\it over the whole} $\CC \PP^2$ (see 1.11(5)).
The vector field is the projection of a linear vector field on $\CC^3$.
Since every integral curve of a linear vector field is parametrized by the whole $\CC$, the same is true for $\cal L$.
It follows that the result of analytic continuation of $\eta$
is a single-valued function. Now $\cal L$ is hyperbolic by Theorem
4.10. Let $\DD\stackrel{\pi}{\longrightarrow} \cal L$ be the covering map, and lift $\eta$
to the universal covering to obtain a holomorphic map $\tilde{\eta}: \CC \to \DD$ with $\pi \circ \tilde{\eta}=\eta$.
By the Liouville's Theorem, $\tilde{\eta}$ will be constant, and this is a contradiction.\\ \\
\vspace {4 mm}
\noindent
{\large \bf Algebraic Leaves and Minimal Sets}\\
The next theorem answers a basic question which is of special interest in the
case of any foliated manifold (cf. \cite{C-L-S1}).\\ \\
{\bf 4.12 Theorem} {\it No SHFC on $\CC \PP^2$ has a compact non-singular leaf.}\\ \\
{\bf Proof.} Let ${\cal L}$ be a compact non-singular leaf of ${\cal F}: \{Pdy-Qdx=0 \}$. By the Chow's Theorem \cite{G-H}, ${\cal L}$ is a (smooth) algebraic curve in $\CC \PP^2$, so the intersection of ${\cal L}$ with the 
curve $S_P:\{ P=0 \} $ is a finite set $\{p_1, \ldots , p_k \}$.

Consider the 1-forms 
$$\alpha=\frac{\partial}{\partial y}\left( \frac{Q}{P}\right ) dx\ \ \ \ \ \ \ \ {\rm and}\ \ \ \ \ \ \ \ \beta=-\partial \log (1+|x|^2+|y|^2)$$
in the affine chart $(x,y)\in U_0$. An easy computation shows that the 1-form
$\tau=\alpha +\beta$ is well-defined on $\CC \PP^2\sm S_P$. For example,
in the affine chart $(u,v)\in U_1$ it is given by  
$$\frac{\partial}{\partial v} \left( \frac{\tilde R}{\tilde P} \right) du-\partial \log (1+|u|^2+|v|^2),$$ 
where $\tilde P$ and $\tilde R$ are polynomials in $u,v$ defined in 4.9. The restriction $\alpha|_{\cal L}$ has poles at the finite set $\{p_1, \ldots , p_k \}$ where $\cal L$ has a vertical tangent line. Without loss of generality we assume that all the $p_j$ are in the affine chart $U_0$. Furthermore, it is easy to compute the residue of $\alpha|_{\cal L}$ at $p_j$: If $p_j=(x_j,y_j)$ and if $y\mapsto x_j+\sum_{i=m_j}^{\infty}a_i(y-y_j)^i$ is the local parametrization of $\cal L$ near $p_j$ with $a_{m_j}\neq 0$, then 
$${\rm Res}[\alpha|_{\cal L}; p_j]=1-m_j.$$
Now consider small disks $D_j \subset \cal L$ around each $p_j$ and integrate $d\tau$ over ${\cal L}'={\cal L}\sm (\bigcup D_j)$:
$$\begin{array}{rl}
\displaystyle{ \int}_{\cal L'}d\tau & =\sum \displaystyle{\int}_{\partial D_j} \tau\\
                    & =\sum \displaystyle{\int}_{\partial D_j} \alpha +\sum \displaystyle{\int} _{\partial D_j} \beta \\
                    & =(2\pi \sqrt{-1}) \sum (1-m_j) + \sum \displaystyle{\int}_{\partial D_j} \beta.
\end{array}$$
On the other hand, $d\tau|_{\cal L'}= d\alpha|_{\cal L'}+d\beta|_{\cal L'}=d\beta|_{\cal L'}=\bar{\partial} \beta |_{\cal L'}=(2\pi \sqrt{-1}) \Omega |_{\cal L'}$, where $\Omega$ is the standard
area form 4.2(2) coming from the Fubini-Study metric. Therefore
$$(2 \pi \sqrt{-1})\ {\rm area}({\cal L'})=(2\pi \sqrt{-1}) \sum (1-m_j)+\sum \int_{\partial D_j}\beta.$$
Letting $D_j$ shrink to $p_j$, we get 
$${\rm area}({\cal L})= \sum (1-m_j) \leq 0,$$
which is a contradiction. $\hfill \Box $\\ 

As a result, every algebraic leaf of a SHFC ${\cal F}$ must have some singular
points of ${\cal F}$ in its closure. Note that every singularity of an algebraic
leaf is indeed a singular point of ${\cal F}$ as well.\\ \\
{\bf 4.13 Proposition} \ {\it Let ${\cal M}$ be a non-trivial minimal set for
a SHFC ${\cal F}$ on $\CC \PP^2$. Then ${\cal M}$ intersects every algebraic
curve in $\CC \PP^2$}.\\ \\
{\bf Proof.} Let $S_K:\{ K=0\} $ be an algebraic curve in $\CC \PP^2$ of degree
$k$ (cf. 1.8). For every triple $(a,b,c)$ of positive real numbers, define
$$\varphi (x,y)=\varphi _{a,b,c}(x,y):=\frac{|K(x,y)|^2}{(a+b|x|^2+c|y|^2)^k}$$
as a non-negative real-analytic function on the affine chart $(x,y)\in U_0$. Since $K$ has
degree $k$, $\varphi$ can be extended to a real analytic function on the entire
$\CC \PP^2$, with $S_K=\varphi ^{-1}(0)$.

Suppose that ${\cal M}\cap S_K=\emptyset$. Then $\varphi$ attains a positive minimum on ${\cal M}$, i.e., there exists $p_0\in {\cal M}$ such that $\varphi(p)\geq
\varphi (p_0)>0$ for all $p\in {\cal M}$. Define $\psi :{\cal M}\rightarrow
{\Bbb R}$ by $\psi(p)= \log \varphi(p)$. Clearly $\psi (p)\geq \log \varphi(p_0) > -\infty$ for all $p\in {\cal M}$. On the other hand, $\psi$ is superharmonic along the non-singular leaf ${\cal L}_{p_0}$. To see this, let $T\stackrel{\eta}{\longmapsto}(x(T),y(T))$ be a local parametrization of ${\cal L}_{p_0}$ near $p_0$, with $\eta (0)=p_0$. Then $\psi (\eta (T))>-\infty $, and
$$\begin{array}{rl}\vspace{3 mm}
\Delta \psi (\eta (T))= & \displaystyle{4 \frac{\partial}{\partial T} \frac{\partial}{\partial \overline{T}} \psi (\eta (T))}\\ \vspace{3 mm}
= & -4k(a+b|x(T)|^2+c|y(T)|^2)^{-2}\\ \vspace{3 mm}
  & (ab|x'(T)|^2+ ac|y'(T)|^2+ bc|x(T)y'(T)-y(T)x'(T)|^2)
\end{array}$$
which is negative. Now $\psi (\eta (T))$ has a minimum at $T=0$. It follows that
$\varphi$ is constant on ${\cal L}_{p_0}$, hence on ${\cal M}$ since $\overline{\cal L}_{p_0}={\cal M}$.
Therefore, for any triple $(a,b,c)$ of positive real numbers, there exists $\alpha >0$
such that
$$|K(x,y)|^{2/k}=\alpha (a+b|x|^2+c|y|^2)$$
for all $(x,y)\in {\cal M}$. Clearly this is a contradiction.\ $\hfill \Box$\\ \\
{\bf 4.14 Corollary} \ {\it No SHFC on $\CC \PP^2$ admitting an algebraic leaf
can have a non-trivial minimal set}.\\ \\
{\bf Proof.} Let ${\cal L}$ be an algebraic leaf of a SHFC ${\cal F}$. By Theorem
4.12, $\overline{\cal L}$ necessarily contains a singular point of $\cal F$, say
$q$. If $\cal M$ is a non-trivial minimal set of $\cal F$, then there exists
$p\in {\cal M}\cap {\cal L}$ by Proposition 4.13. As $p$ is non-singular, $\overline{\cal L}=
\overline{\cal L}_p=\cal M$, so $q\in {\cal M}\cap {\rm sing}(\cal F)$, which is a contradiction.\ $\hfill \Box$\\

Therefore, in order to find an $\cal F$ with a non-trivial minimal set, we must
look for $\cal F$ in the sub-class of ${\cal D}_n$ consisting of SHFC's which do not
admit any algebraic leaf. This sub-class is open and dense in ${\cal D}_n$ (see Remark 3.23).

Note that the above corollary gives another proof of the fact that no SHFC of
geometric degree 1 can have a non-trivial minimal set, since one can easily see that
every SHFC of geometric degree 1 has a projective line as a leaf.

We conclude with few important remarks.\\ \\
{\bf 4.15 Remarks}

(a) It is shown in \cite{C-L-S1} that each leaf $\cal L$ contained in the non-trivial
minimal set $\cal M$ has {\it exponential growth}. This means that if we fix
some Riemannian metric on $\cal L$ and some $p\in {\cal L}$, then
$$\liminf_{r\rightarrow +\infty}\frac{\log ({\rm area}(B_r(p)))}{r}>0,$$
where $B_r(p)$ denotes the open ball in $\cal L$ of radius $r$ centered at $p$.

(b) Recently C. Bonatti, R. Langevin, and R. Moussu \cite{B-L-M} have shown that for
any non-trivial minimal set $\cal M$, there exists a leaf ${\cal L}\subset {\cal M}$ such
that the monodromy group $G({\cal L})$ contains a hyperbolic germ in Bih$_0(\CC)$. The
real version of this theorem is a famous 1965 result of R. Sacksteder \cite{Sa}: An exceptional minimal set of a transversely orientable codimension
one $C^2$ foliation on a compact manifold contains a leaf with a hyperbolic monodromy mapping.

(c) Here is a related result due to A. Candel and X. G\'omez-Mont \cite{C-G} (see also the paper by A. Glutsyuk in \cite{I6} for a generalization): Let $\cal F$ be a SHFC with no algebraic leaves and all singular points of hyperbolic type. Then every leaf of $\cal F$ is a hyperbolic Riemann
surface. In fact, a non-hyperbolic leaf gives rise to a non-trivial invariant
transverse measure for $\cal F$. The support of this measure cannot intersect 
the leaves outside of the (possible) minimal set since then it has to be supported on the (global) separatrices by hyperbolicity of the singular points, and this means that these separatrices (with singular points added) are compact,
hence algebraic, which is a contradiction. On the other hand, no invariant 
transverse measure can live on the minimal set by a result of \cite{C-L-S1}. 

\newpage
\begin{center}
{\sc Symbol Index}
\end{center}
\vspace*{0.15in}
The number in each item denotes the section where the item appears for the first time.\\ \\
\begin{tabular}{lll}
$[A]$ & {\small germ at 0 of the set $A$ in $\CC$} & 3.6\vspace{2 mm}\\
${\cal A}_n$ & {\small class of all singular foliations on $\CC \PP^2$ having affine degree $n$} & 1.27\vspace{2 mm}\\
${\cal A}_n'$ & {\small class of all singular foliations in ${\cal A}_n$ having Petrovski\u\i-Landis property} & 1.34\vspace{2 mm}\\
$B(n)$ & {\small standard line bundle on $\CC \PP^2$ with the first Chern class $n$} & 1.18\vspace{2 mm}\\
${\rm Bih}_0(\CC)$ & {\small group of germs at 0 of biholomorphisms of $\CC$ fixing the origin} & 2.5\vspace{2 mm}\\
$c_1(B)$ & {\small first Chern class of the line bundle $B$} & 1.18\vspace{2 mm}\\
$\DD,\DD(z,r)$ & {\small open unit disk, open disk of radius $r$ centered at $z$} & 1.4\vspace{2 mm}\\
$DG$ & {\small tangent group of a subgroup $G$ of ${\rm Bih}_0(\CC)$} & 3.4\vspace{2 mm}\\
${\cal D}_n$ & {\small class of all singular foliations on $\CC \PP^2$ having geometric degree $n$} & 1.23\vspace{2 mm}\\
$\doteq$ & {\small equal up to a set of measure zero} & 3.6\vspace{2 mm}\\
${\cal F}$ & {\small a (singular) holomorphic foliation} & 1.1\vspace{2 mm}\\
${\cal F}_\omega$ & {\small singular foliation induced by a polynomial 1-form $\omega $ on $\CC \PP^2$} & 1.6\vspace{2 mm}\\  
${\cal F}_X$ & {\small singular foliation induced by a holomorphic vector field $X$} & 1.4\vspace{2 mm}\\
$f_\gamma$ & {\small holonomy (monodromy) mapping associated with a curve (loop) $\gamma $} & 2.1\vspace{2 mm}\\
$\phi_{ij}$ & {\small change of coordinates for standard atlas of $\CC \PP^2$} & 1.7\vspace{2 mm}\\
$G_\infty$ & {\small monodromy group of the leaf at infinity} & 2.14\vspace{2 mm}\\
$G({\cal L})$ & {\small monodromy group of the leaf ${\cal L}$} & 2.6\vspace{2 mm}\\
$\kappa$ & {\small Gaussian curvature of a given conformal metric on a Riemann surface} & 4.7\vspace{2 mm}\\
$\lambda _j$ & {\small characteristic number of a singular point on the line at infinity} & 2.15\vspace{2 mm}\\
${\cal L}_ \infty $ & {\small the leaf at infinity} & 1.26\vspace{2 mm}\\
${\cal L}_p$ & {\small the leaf through $p$ of given (singular) foliation} & 1.1\vspace{2 mm}\\
$L_0,L_1,L_2$ & {\small lines at infinity with respect to the affine charts $U_0,U_1,U_2$ in $\CC \PP^2$} & 1.7\vspace{2 mm}\\
$m$ & {\small Lebesgue measure on $\CC$} & 3.8\vspace{2 mm}\\
${\cal M}$ & {\small non-trivial minimal set of a singular foliation on $\CC \PP^2$} & 4.1\vspace{2 mm}\\
$\nu _j$ & {\small multiplier at 0 of a generator of $G_\infty $} & 2.16\vspace{2 mm}
\end{tabular}
\newpage
\noindent
\begin{tabular}{lll}
${\cal O}$ & {\small sheaf of holomorphic functions on a complex manifold} & 1.18\vspace{2 mm}\\
${\cal O^\ast }$ & {\small sheaf of non-vanishing holomorphic functions on a complex manifold} & 1.17\vspace{2 mm}\\
$PG$ & {\small pseudo-group obtained from a finitely-generated subgroup $G$ of Bih$_0(\CC)$} & 2.8\vspace{2 mm}\\
$\rho _X$ & {\small Poincar\'e metric on a hyperbolic Riemann surface $X$} & 4.6\vspace{2 mm}\\
$\cal S$ & {\small space of all SHFC's on $\CC \PP^2$} & 1.24\vspace{2 mm}\\
$s(A)$ & {\small saturation of a set $A$ in $\CC \PP^2$} & 3.26\vspace{2 mm}\\
SHFC & {\small Singular Holomorphic Foliation by Curves} & 1.6\vspace{2 mm}\\
sing$({\cal F})$ & {\small singular set of a singular foliation ${\cal F}$} & 1.3\vspace{2 mm}\\
$S_P$ & {\small algebraic curve in $\CC \PP^2$ defined by $P=0$} & 1.8\vspace{2 mm}\\
$U_0,U_1,U_2$ & {\small affine charts for standard atlas of $\CC \PP^2$} & 1.7
\end{tabular}
\vspace*{2 cm}

\newpage
\begin{thebibliography}{*****}
\bibitem [\bf{Ah}]{Ah} L. Ahlfors, {\it Conformal Invariants: Topics in Geometric Function Theory}, McGraw-Hill, New York, 1973.
\bibitem [\bf{An}]{An} D. Anosov, V. Arnold (eds.), {\it Dynamical Systems (I)}, Encyclopedia of Mathematical Sciences, Springer-Verlag, 1988.
\bibitem [\bf{Ar}]{Ar} V. Arnold, {\it Geometrical Methods in the Theory of Ordinary Differential Equations}, 2nd ed., Springer-Verlag, 1988.
\bibitem [\bf{AA}]{A-A} V. Arnold, A. Avez, {\it Ergodic Problems of Classical Mechanics}, Benjamin, Inc., New York, 1968.
\bibitem [\bf{BLM}]{B-L-M} C. Bonatti, R. Langevin, R. Moussu, {\it Feuilletages de $\CC \PP(n)$: De l'Holonomie Hyperbolique pour les Minimaux Exceptionnels}, Publ. Math. IHES, {\bf 75} (1992) 124-134.
\bibitem [\bf{C}]{C} C. Camacho, {\it Problems on Limit Sets of Foliations on Complex Projective Spaces}, Proceedings of the International Congress of Mathematicians, Kyoto, Japan, 1990, 1235-1239.
\bibitem [\bf{CKP}]{C-K-P} C. Camacho, N. Kuiper, J. Palis, {\it The Topology of Holomorphic Flows with Singularity}, Publ. Math. IHES, {\bf 48 }(1978) 5-38.
\bibitem [\bf{CL}]{C-L} C. Camacho, A. Lins Neto, {\it Geometric Theory of Foliations}, Birkhauser Boston, 1985.
\bibitem [\bf{CLS1}]{C-L-S1} C. Camacho, A. Lins Neto, P. Sad, {\it Minimal Sets of Foliations on Complex Projective Spaces,} Publ. Math. IHES, {\bf 68} (1988) 187-203.
\bibitem [\bf{CLS2}]{C-L-S2} C. Camacho, A. Lins Neto, P. Sad, {\it Foliations with Algebraic Limit Sets}, Annals of Math.,{\bf 136} (1992) 429-446.
\bibitem [\bf{CS}]{C-S} C. Camacho, P. Sad, {\it Invariant Varieties Through Singularities of Holomorphic Vector Fields}, Annals of Math., {\bf 115} (1982) 579-595.
\bibitem [\bf{Ca}]{Candel} A. Candel, {\it Uniformization of Surface Laminations}, Ann. Scient. \'Ec. Norm. Sup., {\bf 26} (1993) 489-516.
\bibitem [\bf{CaG}]{C-G} A. Candel, X. G\'omez-Mont, {\it Uniformization of the Leaves of a Rational Vector Field}, Ann. Inst. Fourier, {\bf 45} (1995) 1123-1133.
\bibitem [\bf{CoL}]{Co-L} E. Coddington, N. Levinson, {\it Theory of Ordinary Differential Equations}, McGraw-Hill, New York, 1965.
\bibitem [\bf{D}]{D} H. Dulac, {\it Sur les Cycles Limites}, Bull. Soc. Math. France, {\bf 51} (1923) 45-188.
\bibitem [\bf{GO}]{G-O} X. G\'omez-Mont, L. Ortiz Bobadilla, {\it Systemas Dinamicos Holomorfos en Superficias,} Sociedad Matematica Mexicana, 1989.
\bibitem [\bf{GL}]{G-L} X. G\'omez-Mont, I. Luengo, {\it Germs of Holomorphic Vector Fields in $\CC^3$ without a Separatrix}, Invent. Math., {\bf 109} (1992) 211-219.
\bibitem [\bf{GH}]{G-H} P. Griffiths, J. Harris, {\it Principles of Algebraic Geometry}, Wiley-Interscience, New York, 1978.
\bibitem [\bf{GR}]{G-R} R. Gunning, H. Rossi, {\it Analytic Functions of Several Complex Variables}, Prentice-Hall, 1965.
\bibitem [\bf{I1}]{I2} Yu. Il'yashenko, {\it Global and Local Aspects of the Theory of Complex Differential Equations}, Proceedings of the International Congress of Mathematicians, Helsinki, Finland, 1978, 821-826.
\bibitem [\bf{I2}]{I3} Yu. Il'yashenko, {\it Singular Points and Limit Cycles of Differential Equations in the Real and Complex Planes}, Preprint, Pushchino: Sci. Res. Comput. Cent., Acad. Sci. USSR, 1982 (in Russian).
\bibitem [\bf{I3}]{I4} Yu. Il'yashenko, {\it Topology of Phase Portraits of Analytic Differential Equations in the Complex Projective Plane},
Trudy Sem. Petrovski\u\i, {\bf 4} (1978) 84-136 (in Russian), English Translation in Sel. Math. Sov., {\bf 5} (1986) 140-199.
\bibitem [\bf{I4}]{I1} Yu. Il'yashenko, {\it Finiteness Theorems for Limit Cycles}, Russ. Math. Surv., {\bf 42} (1987) 223.
\bibitem [\bf{I5}]{I5} Yu. Il'yashenko (ed.), {\it Nonlinear Stokes Phenomena}, AMS Publications, Providence RI, 1992.
\bibitem [\bf{I6}]{I6} Yu. Il'yashenko (ed.), {\it Differential Equations with Real and Complex Time}, Proceedings of the Steklov Institute of Mathematics, vol. 213, 1996. 
\bibitem [\bf{IY}]{I-Y} Yu. Il'yashenko, S. Yakovenko, {\it Analytic Differential Equations}, ICTP Lecture Notes of School on Dynamical Systems, 1991.
\bibitem [\bf{J}]{J} J.P. Jouanolu, {\it Equations de Pfaff Algebriques}, Lecture Notes in Mathematics 708, Springer-Verlag, 1979.
\bibitem [\bf{Kob}]{Kob} S. Kobayashi, {\it Hyperbolic Manifolds and Holomorphic Mappings}, Marcel Dekker, New York, 1970.
\bibitem [\bf{Kod}]{Kod} K. Kodaira, {\it Complex Manifolds and Deformation of Complex Structures}, Springer-Verlag, 1986.
\bibitem [\bf{Kr}]{Kr} S. Krantz, {\it Complex Analysis: The Geometric Viewpoint}, The Carus Mathematical Monographs No. 23, MAA, 1990.
\bibitem [\bf{L}]{L} A. Lins Neto, {\it Algebraic Solutions of Polynomial Differential Equations and Foliations in Dimension Two}, in Holomorphic Dynamics, Lecture Notes in Mathematics 1345, Springer-Verlag, 1988.
\bibitem [\bf{M}]{M1} J. Milnor, {\it Dynamics in One Complex Variable}, SUNY Stony Brook IMS Preprint, 1990/5.
\bibitem [\bf{N}]{N} V. Naishul, {\it Topological Invariants of Analytic and Area-Preserving Mappings and Their Application to Analytic Differential Equations in $\CC^2$ and $\CC \PP^2$}, Trans. Moscow Math. Soc., Issue {\bf 2} (1983) 234-250.
\bibitem [\bf{P}]{P} R. Perez-Marco, {\it Solution Compl\`ete au Probl\`eme de Siegel de Lin\'earisation d'une Application
Holomorphe au Voisinage d'un Point Fixe (d'apr\`es J.C. Yoccoz)},\ S\'eminaire Bourbaki, {\bf 753} (1991-1992) 1-31.
\bibitem [\bf{PL1}]{P-L1} I. Petrovski\u\i, E. Landis, {\it On the Number of Limit Cycles of the Equation
$dy/dx=P(x,y) / Q(x,y)$ where $P$ and $Q$ are Polynomials of Degree 2, } Amer.
Math. Society Trans., Series 2, {\bf 10} (1963) 125-176.
\bibitem [\bf{PL2}]{P-L2} I. Petrovski\u\i, E. Landis, {\it On the Number of Limit Cycles of the Equation
$dy/dx=P(x,y) / Q(x,y)$ where $P$ and $Q$ are Polynomials}, Amer.
Math. Society Trans., Series 2, {\bf 14} (1964) 181-200.
\bibitem [\bf{Re}]{Re} G. Reeb, {\it Sur Certaines Propri\'et\'es Topologiques des Vari\'et\'es Feuillet\'ees}, Actualit\'es Sci.
Indust., Hermann Paris, {\bf 1183} (1952) 91-154.
\bibitem [\bf{Ru}]{Ru} W. Rudin, {\it Real and Complex Analysis}, 3rd Ed., McGraw-Hill, New York, 1986.
\bibitem [\bf{Sa}]{Sa} R. Sacksteder, {\it Foliations and Pseudo-Groups}, Amer. J. Math., {\bf 87} (1965) 79-102.
\bibitem [\bf{Sh}]{Sh} A. Shcherbakov, {\it Topological and Analytical Conjugacy of Non-abelian Groups of Germs of Conformal Mappings}, Trudy Sem. Petrovski\u\i, {\bf 10} (1984) 170-196 (in Russian).
\bibitem [\bf{W}]{W} H. Whitney, {\it Complex Analytic Varieties}, Addison-Wesley, Reading MA, 1972.
\end{thebibliography}
\end{document}